\documentclass[11pt]{amsart}
\usepackage{amsmath,amsxtra,amssymb}
\numberwithin{equation}{section}
\newtheorem{lemma}{Lemma}[section]
\newtheorem{prop}[lemma]{Proposition}
\newtheorem{theorem}[lemma]{Theorem}
\newtheorem{cor}[lemma]{Corollary}

\newtheorem{prob}[lemma]{Problem}
\newtheorem{rem}[lemma]{Remark}

\newcommand{\re}{\begin{rem}\rm}
  \newcommand{\mar}{\end{rem}}
\newtheorem{exam}[lemma]{Example}

\newcommand{\kla}{\left ( }
\newcommand{\mer}{\right ) }

\renewcommand{\for}{\begin{eqnarray*}}

\newcommand{\mel}{\end{eqnarray*}}

\newcommand{\kl}{\pl \le \pl}
\newcommand{\gl}{\pl \ge \pl}

\newcommand{\lel}{\pl = \pl}

\newcommand{\ez}{{\mathbb E}}

\newcommand{\nz}{{\mathbb N}}
\newcommand{\nen}{n \in \nz}

\newcommand{\rz}{{\mathbb R}}
\newcommand{\zz}{{\mathbb Z}}

\newcommand{\Mz}{{\mathbb M}}

\newcommand{\cz}{{\mathbb C}}
\newcommand{\kz}{{\rm  I\! K}}

\newcommand{\ten}{\otimes}

\newcommand{\wet}{\stackrel{\wedge}{\otimes}}

\DeclareMathOperator{\gr}{gr}

\DeclareMathOperator{\Pe}{Perm}

\newcommand{\p}{\hspace{.05cm}}
\newcommand{\pl}{\hspace{.1cm}}
\newcommand{\pll}{\hspace{.3cm}}

\newcommand{\hz}{\vspace{0.5cm}}
\newcommand{\lz}{\vspace{0.0cm}}

\newcommand{\qd}{\end{proof}\vspace{0.5ex}}

\newcommand{\Om}{\Omega}
\newcommand{\om}{\omega}

\newcommand{\al}{\alpha}
\newcommand{\si}{\sigma}
\newcommand{\Si}{\Sigma}

\newcommand{\la}{\lambda}
\newcommand{\eps}{\varepsilon}

\newcommand{\E}{{\mathcal E}}
\newcommand{\A}{{\mathcal A}}

\newcommand{\K}{{\mathcal K}}

\newcommand{\Ma}{{\mathbb M}}

\newcommand{\M}{{\mathcal M}}

\newcommand{\N}{{\mathcal N}}

\newcommand{\g}{\gamma}

\newcommand{\U}{{\mathcal U}}

\newcommand{\noo}{\left \|}
\newcommand{\rrm}{\right \|}

\newcommand{\bet}{\left |}
\newcommand{\rag}{\right |}

\newcommand{\intt}{\int\limits}
\newcommand{\summ}{\sum\limits}

\newcommand{\prodd}{\prod\nolimits}

\newcommand{\lb}{\langle \langle }
\newcommand{\rb}{\rangle \rangle }

\newcommand{\8}{\infty}

\allowdisplaybreaks

\begin{document}

\title[Operator spaces and Araki-Woods
factors]{Operator spaces and Araki-Woods factors\\
-- A quantum probabilistic approach --}
\author[M. Junge]{M. Junge$^*$}
\address{Department of Mathematics\\
University of Illinois, Urbana, IL 61801, USA} \email[Marius
Junge]{junge@math.uiuc.edu}

\thanks{${}^*$ The author is partially supported by the
National Science Foundation Foundation DMS-0301116 and DMS
05-56120}

\begin{abstract} We show that the operator Hilbert space OH
introduced by Pisier embeds into the predual of the hyerfinite
III$_1$ factor. The main new tool is a Khintchine type inequality
for the generators of the CAR algebra with respect to a quasi-free
state. Our approach  yields a Khintchine type inequality for the
$q$-gaussian variables for all values $-1\le q\le 1$. These
results are closely related to recent results of Pisier and
Shlyakhtenko in the free case.
\end{abstract}  \maketitle

\vspace{-0.2cm}

 {\bf Content:}
 \begin{enumerate}
 \item[0.] Introduction and notation \hfill 1
 \item[1.] Preliminaries \hfill 5
 \item[2.] The algebraic central limit theorem \hfill 6
 \item[3.] CAR, CCR and $q$-commutation relations \hfill 12
 \item[4.] Limit Distributions \hfill 20
 \item[5.] A uniqueness result \hfill 27
 \item[6.] An inequality for sums of independent copies \hfill 33
 \item[7.] Khintchine type inequalities \hfill 44
 \item[8.] Applications to operator spaces \hfill 49
 \end{enumerate}

\vspace{0.2cm}

\setcounter{section}{-1}
\section{Introduction and Notation}

Probabilistic methods play an important role in the theory of
operator algebras and Banach spaces. It is not surprising that a
quantized theory of Banach spaces will require tools from quantum
probability. This connection between noncommutative probability
and the recent theory of operator spaces (sometimes called
quantized Banach spaces) is well-established through the work of
Haagerup, Pisier \cite{HP} and the general theory of Khintchine
type inequalities by Lust-Piquard \cite{LP}, Lust-Piquard and
Pisier \cite{LPP}. The importance of type III von Neumann algebras
in this line of research was discovered only recently through the
work of  Pisier/Shlyahtenko \cite{PS} on Grothendieck's theorem
for operator spaces and in \cite{Joh}. Both papers make essential
use of the theory of free probability. It is well-known that in
free probability theory some probabilistic estimates, classically
only valid for $p<\infty$, hold even for $p=\infty$. For example
this holds for Biane/Speicher's \cite{BS} work on stochastic
process and Voiculescu's inequality for sums of free random
variables \cite{Voi2}). In this paper we prove norm estimates for
the sum of independent copies in noncommutative $L_1$ spaces in a
quite general setting. This includes free random variables as  in
\cite{Joh} and also classical commuting or anti-commuting random
variables. Using a central limit procedure, similar as in
\cite{Joh}, we derive Khintchine type inequalities for the
classical Araki-Wood factors. Although our results are motivated
by the theory of operator spaces, the techniques used in the proof
are (quantum) probabilistic in nature.

Let us fix the notation required to state the main inequality of
this paper. Let us assume that we have an inclusion of von
Neumann algebras $\Mz\subset A,B\subset \N$ and that $E:\M\to \Mz$
is a faithful normal conditional expectation. Then $A,B$ are
called \emph{independent} over $\Mz$ (better over $E$) if
 \[ E(ab)\lel E(a)E(b) \]
holds for all $a\in A$ and $b\in B$. Let us now assume that
$\Mz\subset M_1,...,M_n\subset\N$ are von Neumann subalgebras.
Then $(\M_i)$ is called (increasingly) independent if $M_i$ is
independent of the von Neumann algebra $\M_{i-1}$ generated by
$M_1,...,M_{i-1}$ over $\Mz$. This definition  is due to
Sch\"{u}rmann.

Our main inequality is an inequality for independent copies,
allowing matrix valued coefficients. Let us therefore assume that
$\Mz\subset M$ and $\al_i:M\to \N$ are faithful homomorphisms such
that $\al_i|_{\Mz}=id$. We say that $(\Mz,M,\al_1,...,\al_n,\N,E)$
is a system of \emph{independent top-subsymmetric copies} if the
$(M_i)$'s are independent and
 \[ E(\al_{i_1}(a_1)\cdots \al_{i_m}(a_m)) \lel
  E(\al_{j_1}(a_1)\cdots \al_{j_m}(a_m)) \]
holds for all $a_1,...,a_m\in M$, and all functions
$\i:\{1,..,m\}\to \{1,...,n\}$,  $\j:\{1,..,m\}\to \{1,...,n\}$
such that $\i|_{\{1,..,m\}\setminus A}=\j|_{\{1,..,m\}\setminus
A}$, $A$ has at most $2$ elements, such that $k\in A$ implies
 \[ i_k\lel \max\{i_1,...,i_m\}\quad ,\quad j_k \lel
 \max\{j_1,...,j_m\}\pl .\]
This means we are allowed to exchange at most two  top values. In
our applications we have often have much stronger assumptions, for
example for free or independent copies. We say that such a system
is \emph{conditioned} if there is a faithful normal state $\phi$
on $\N$ such that $\phi\circ E\lel \phi$ and
 \[ \si_t^{\phi}(\al_i(M))\subset \al_i(M)\]
holds for all $i=1,...,n$. This allows us to extends the maps
$\al_i$ to all $L_p$ spaces. Our main inequality is an estimate
for sums of independent copies in $L_1$:

\begin{theorem}\label{sumind} Let $(\Mz,M,\al_1,...,\al_n,\N,E)$ be a system of
independent, conditioned top-sub- symmetric copies. Then
 \[ \ez \|\summ_{k=1}^n \eps_k \al_k(x)\|_1
 \sim \inf_{x=x_1+x_2+x_3} n\|x_1\|_1 +
 \sqrt{n}\|E(x_2^*x_2)^{1/2}\|_1+ \sqrt{n}\|E(x_3x_3^*)^{1/2}\|_1
 \pl .\]
\end{theorem}
Here  $(\eps_k)_{k=1}^n $ is a sequence of independent Bernoulli
variables with Prob$(\eps_k=\pm 1)=\frac12$. We will use the
symbol $a\sim b$ if there exists an absolute constant $c>0$ such
that $c^{-1}a\le b\le ca$ (of course independent of $x$ in the
theorem above). For non-tracial von Neumann algebras which occur
in the context of free probability, we work with Haagerup's $L_1$
spaces and the `natural' extension of $\al_k$ and $E$ to these
spaces (see e.g. \cite{JD}, \cite{JX}). We may replace $\eps_k$ by
$\eps_k\ten v_k$ where $v_k$ are unitaries. This is important in
the context of Speicher's \cite{Sp} interpolation technique for
$q$-commutation relations.

An essential ingredient in our proof of Theorem \ref{sumind} is
the noncommutative Khintchine inequality due to Lust-Piquard and
Pisier  \cite{LPP}:
 \[ \ez \|\summ_{k=1}^n \eps_k x_k\|_1
 \sim \inf_{x_k=c_k+d_k} tr\big((\sum_{k=1}^n c_k^*c_k)^{1/2}\big)+
 tr\big( (\sum_{k=1}^n d_kd_k^*)^{1/2}\big) \pl .\]
Our passage from three terms above to two terms uses the central
limit theorem. Given a sequence $(x_k)$ of classical independent
copies, the central limit theorem tells us that
$n^{-1/2}\sum_{k=1}^n x_k$ converges to a gaussian variable.
Central limit theorems in quantum probability have a long and
impressive tradition starting with  the work of Cushen/Hudson and
\cite{Cu-Hud}, see also Hudsen \cite{Hud0}, von Waldenfels
\cite{vWa}, Giri/von Waldenfels \cite{GiW}, Cockroft, Gudder and
Hudsen \cite{cock}, Quaegebeur\cite{Quaeg} and Hegerfeldt
\cite{Heg} among many others. Our interest in the central limit
theorem is two-fold. First, we consider limits of the form
 \[ u_n(x) \lel \sqrt{\frac{T}{n}} \summ_{k=1}^n u_k \ten  \pi_k(x) \]
where $\pi_k:N\to N^{\ten_n}$ is the homomorphism which sends $N$
in the $k$-th component and the $u_k\in \Mz_{2^n}$ are unitaries
satisfying the CAR-relations $u_ku_j=-u_ju_k$. Let $\psi$ be a
state on $N$, and $\tau_n$ be the normalized trace on
$\Mz_{2^n}$. Using the classical combinatorial approach, we see
that
 \begin{align}\label{comb-1}
 \lim_n (\tau_n \ten \psi^{\ten_n})(u_n(x_1)\cdots u_n(x_m))
 &=  \summ_{\si=\{\{i_1,j_1\},...,\{i_{\frac{m}{2}},j_{\frac{m}{2}}\}\}
 \in P_2(m)}
 (-1)^{I(\si)} \prod_{l=1}^{\frac{m}{2}} (T\psi)(x_{i_l}x_{j_l}) \pl
 .
 \end{align}
Here $P_2(m)$ stands for the set of pair partitions of
$\{1,...,m\}$ and $I(\si)$ for the number of inversions. Using
Speicher's trick it is not difficult to replace $(-1)$ by any $q$,
$-1\le q\le 1$. For $q=1$ we see that algebraically the formal
limit object $u_{\infty}(x)$ satisfies
 \[ [u_{\infty}(x),u_{\infty}(y)]\lel T\psi([x,y]) \pl .\]
A suitable change of variables yields  the classical commutation
relations. One problem in our paper is to associate with the limit
object $u_{\infty}(x)$ a selfadjoint operator affiliated to a
suitable von Neumann algebra. This is achieved using ultraproducts
of von Neumann algebras and the fundamental work of Raynaud
\cite{Ra} (see section 4). The reader might object that generators
for the CAR algebra are already bounded and hence there is no need
for this ultraproduct procedure. This is where the second interest
in the central limit procedure becomes apparent. We also want to
guarantee that the norm estimates from Theorem \ref{sumind} hold
for the limit object in $L_1$. This requires extra knowledge on
the action of the modular group and adds technical difficulties to
this paper. Moreover, our approach works uniformly for all $-1\le
q\le 1$. In fact we find a simultaneous realization of all the
$q$-commutation relations (which seems to be new). Note that for
$q=1$ the limit objects are classical gaussian variables which are
indeed unbounded. In section 2 we establish the combinatorial
aspect of the central limit theorem. The connection to the
classical CAR and CCR relations and the Araki-Woods factors is
established in section 3. In section 5 we prove  a noncommutative
version of the Hamburger moment problem which allows us to
identify von Neumann algebras using combinatorial information as
in \eqref{comb-1}. The proof of Theorem \ref{sumind} is contained
in section 6. Combining these results in section 7  we obtain  a
Khintchine type inequality for quasi-free states. More precisely,
we consider the CAR-algebra $\A$ generated by a sequence $(a_k)$
sequence satisfying the canonical anti-commutation relations
 \[ a_ka_j+a_ja_k\lel 0 \quad ,\quad a_ka_j^*+a_j^*a_k\lel
 \delta_{kj} \pl .\]
A quasi free state $\phi_{\mu}$ is characterized by
 \begin{equation}\label{carrk}  \phi_{\mu}(a_{i_r}^*\ten \cdots a_{i_1}^*\ten a_{j_1}\ten
 \cdots \ten a_{j_s}) \lel \delta_{r,s} \prod_{l=1}^r \delta_{i_l,j_l}\mu_{i_l}
 \end{equation}
for all increasing sequences  $i_1<i_2<\cdots <i_r$ and
$j_1<i_2<\cdots <j_s$. Let us also recall the usual notation
$x.\phi(y)\lel \phi(xy)$. The Khintchine inequality for quasi-free
states reads as follows:

\begin{theorem}\label{kh-carr} Let $(x_k)\subset \Mz_m$ be matrices. Then
 \begin{align*}
 \|\summ_k  x_k.tr \ten a_k.\phi_{\mu}\|_{(\Mz_m\ten \A)^*}
 &\sim \inf_{x_k\lel c_k+d_k}
  tr\bigg((\summ_k  \mu_k c_k^*c_k)^{1/2}\bigg) +
  tr\bigg((\summ_k  \frac{\mu_k^2}{1-\mu_k} d_kd_k^*)^{1/2}\bigg) \pl .
 \end{align*}
\end{theorem}

Analogous results hold for all $-1\le q\le 1$. The case $q=0$ in
Theorem \ref{kh-carr} in this form follows from \cite{Pl} but main
ideas are already contained in \cite{PS} and somehow independently
in \cite{Joh}. Using techniques from \cite{Joh}, we deduce from
Theorem \ref{kh-carr} formulas for
 \[ \|\summ_k c_{kj} \pl a_k.\phi \ten a_j.\phi \|_{(\A\ten \A)^*} \pl
 .\]
This leads to a four term infimum. We refer to section 8 for the
precise formula which is slightly involved. At the time of this
writing a proof seems  impossible without using free probability.
We conclude the paper with an application to the theory of
operator spaces. Only in section 8 we require the reader to have
some background in operator spaces theory. Indeed, the reader
familiar with the theory of operator spaces recognizes immediately
that Theorem \ref{kh-carr} is a result on the operator space
structure of $\overline{\rm span}(a_k.\phi_{\mu})$. Let us recall
the definition of the column space $C=\overline{\rm
span}\{e_{k1}:k\in \nz\}$  and the row space $R=\overline{\rm
span}\{e_{1k}:k\in \nz\}$ as  subspaces of the compact operators
$\K(\ell_2)$. An operator space comes either with a concrete
embedding $\iota:X\to B(H)$ or with a sequence  $(\|\pl \|_m)$ of
norms on $(M_m(X))$ satisfying Ruan's axioms (see \cite{ER} or
\cite{Po} for more information on operator spaces). Ruan's theorem
tells us  that if these axioms are satisfied, then there is an
embedding $\iota:X\to B(H)$ such that
 \[ \|[x_{kl}]\|_{\Mz_m(X)} \lel \|[\iota(x_{kl})]\|_{\Mz_m(B(H))}
 \pl .\]
The theory of operator spaces is closed under taking subspaces,
dual spaces and quotient spaces. For example the sequence of
matrix norms on $X/Y$ defined by
 \[ \|[x_{kl}+Y]\|_{\Mz_m(X/Y)} \lel \inf_{y_{kl}\in \Mz_m(Y)} \|[x_{kl}-y_{kl}]\|_{\Mz_m(Y)} \pl \]
satisfies Ruan's axioms. The isomorphism in this category are
those maps which control uniformly all matrix norms, i.e. maps
$u:X\to Y$ such that $\|u\|_{cb}\|u^{-1}\|_{cb}$ is finite. Here
the cb-norm is defined as
 \[ \|u\|_{cb} \lel \sup_m \|id_{\Mz_m} \ten u: \Mz_m(X)\to
 \Mz_m(Y)\|  \pl. \]
The following application to the theory of operator spaces is
proved in section 8.

\begin{cor}\label{emb0} Let $Q$ be a quotient of $R\oplus C$. Then $Q$ embeds
completely isomorphically into the predual of the hyperfinite {\rm
III}$_{\la}$ factor, $0<\la\le 1$. In particular, the operator
space OH completely embeds into the predual of the {\rm III}$_\la$
factor.
\end{cor}

Combining this application with the results in \cite{PS}, we
obtain an analogue of Grothendieck's characterization of Hilbert
spaces (see Corollary \ref{caracG}). The connection between
quotients of $R\oplus C$ and the classical Araki-Wood factors is
very tight. Both objects are essentially determined by the graph
of an unbounded operator on $\ell_2$, indeed it suffices  to
consider diagonal operators on $\ell_2$. Retrospectively, it seems
that the Araki-Wood factors constructed from CAR-relations are a
perfect fit for quotients of $R\oplus C$-an interesting connection
to operator space theory. We refer to \cite{JX-ros} for Khintchine
inequalities for the CAR generators in $L_p$. At the moment the
result holds for $p=1$ and $p>1$, but the known constants $c_p$
are not bounded as $p\to 1$. It is an open problem how to close
this gap. At the end of section 8 we construct a completely
isomorphic embeddings of $OH_n$ in the Brown algebra $B_m^{nc}$,
the universal algebra of coefficients of an $m\times m$ unitary.
Brown algebras have been introduced in the context of K-theory by
\cite{Bro} and can be understood as noncommutative analogues of
the full $C^*$-algebra of the free group. For an embedding of the
infinite dimensional $OH$ we have to use a suitable direct limit
of the Brown algebras. These result could be considered as a first
step towards constructing a `concrete' embedding of the operator
space $OH$, a problem which remains open.

I owe special thanks to Quanhua Xu for many clarifying discussions
on this topic and for insisting on a formulation with the CAR
generators. I would also like to thank Eric Ricard for fruitful
discussions and Anthony Yew for reading several versions of this
paper.

\section{Preliminaries}

We use standard notation in operator algebras
\cite{Tak,TK-II,TK-III}, \cite{kar-I,kar-II} and
\cite{Strat,Strat-Z}. We refer to \cite{ER} and \cite{Po} for
standard references on operator spaces. The relevant information
on the projective tensor product for noncommutative $L_1$ spaces
can also be found in \cite{Joh}. Let us recall very briefly some
facts about Haagerup's $L_p$ spaces associated to a von Neumann
algebra $N$. Let $\phi$ be a normal faithful weight with modular
group $\si_t^{\phi}$. We obtain a representation $\pi$ of $N$ on
the Hilbert space $L_2(N,\phi)$, the completion of $N$ with
respect to the scalar product $(x,y)=\phi(x^*y)$. On the crossed
product $N\rtimes_{\si_t^{\phi}}\rz$ there exists a unique trace
$\tau$ such that the dual action $\theta_s$ satisfies
$\tau(\theta_s(x))=e^{-s}\tau(x)$ for all $x$ in the positive cone
of $N\rtimes_{\si_t^{\phi}}\rz$. Haagerup's space
$L_p(N)=L_p(N,\phi)$ is defined as the space of $\tau$ measurable
operators $x$ affiliated to  $N\rtimes_{\si_t^{\phi}}\rz$ such
that $\theta_s(x)=e^{-s/p}x$. $L_p(N)$ is a left and right $N$
module. On the space $L_1(N)$ Haagerup  defines a unique linear
functional $tr:L_1(N)\to \cz$ such that
 \[ tr(Dy)\lel tr(yD) \]
holds for all $y\in N$ and $D\in L_1(N)$.  This functional defines
a one-to-one (linear) correspondence between the predual $N_*$
and $L_1(N)$ given by $\phi_D(y)\lel tr(xy)$. We call the operator
$D$ the \emph{density} of $\phi$. The norm in $L_p(N)$ is defined
by
 \[ \noo x\rrm_p \lel [tr(|x|^p)]^{\frac1p} \pl .\]
It is important to note for $x\in L_p(N)$ the polar decomposition
$x=u|x|$ satisfies $u\in N$ and $|x|\in L_p(N)$. $L_p(N)$ is a
Banach space for $1\le p<\infty$ and $p$-normed for $0<p\le 1$.
We refer to \cite{Te} for proofs and more information. Let us
consider a further semifinite von Neumann algebra $M$ with trace
$\tau$. Then $\tau\ten \phi$ is a normal faithful weight on the
tensor product $M\bar{\ten}N$. Using $\si_t^{\tau\ten \phi}=1\ten
\si_t^{\phi}$ we may identify
$(M\bar{\ten}N)\rtimes_{\si_t^{\tau\ten \phi}}\rz$ and
$M\bar{\ten}(N\rtimes_{\si_t^{\phi}}\rz)$. This implies that for
$a\in L_p(M)$ and $x\in L_p(N)$ we have $x\ten a\in
L_p(M\bar{\ten}N)$. Moreover, the tracial functional satisfies
 \[ tr_{M\ten N} \lel \tau\ten tr_N \pl .\]
An easy application of Kaplansky's density theorem shows that
$L_p(M)\ten L_p(N)$ is a dense subspace of $L_p(M\bar{\ten}N)$.
For more information on tensor products of noncommutative $L_p$
spaces in the non-tracial case we refer to \cite{JF}.

The natural operator space structure on $L_1(N)$ is inherited from
$N_*^{op}$. Let us recall that $N^{op}$ is the von Neumann algebra
$N$ equipped with the reverse multiplication $x\circ y=yx$. Since
$N$ and $N^{op}$ coincide as Banach spaces, the same holds true
for $N_*$ and $N^{op}_*$. In particular, the mapping
$\iota:L_1(N)\to N^{op}_*$, $\iota(x)(y)=tr(xy)$ is well-defined.
As explained in \cite{JR-ap} and \cite{Joh}, the advantage of
this embedding is the equality
 \begin{equation}\label{wetL1}  L_1(M)\wet L_1(N)\lel L_1(M\bar{\ten}N)
 \end{equation}
which holds for all von Neumann algebras $M$. We will also
frequently use the following results of Effros and Ruan (see
\cite{ER})
 \begin{equation}\label{erdual}  (M_*\wet N_*)^*\lel M\bar{\ten}N \pl .
 \end{equation}
Let us consider the special case  $M=B(\ell_2)$. Then, we see by
complementation that
 \begin{align*}
  &\noo \summ_k e_{k1}\ten x_k\rrm_{S_1\wet L_1(N)}
  \lel \noo \summ_k e_{k1}\ten
  x_k\rrm_{L_1(B(\ell_2)\bar{\ten}N)} \lel
  \noo (\summ_k x_k^*x_k)^{\frac12} \rrm_{L_1(N)} \pl .
 \end{align*}
On the other hand, we deduce from \eqref{wetL1} that
 \begin{align*}
 &\noo \summ_k e_{1k}\ten x_k\rrm_{R\wet L_1(N)}
 \lel \sup\left\{ \bet\langle \summ_k e_{1k}\ten x_k, u\rangle\rag \pl:\pl
 \|u:R\to N^{op}\|_{cb}\le 1\right\} \\
 &= \sup\{ | \summ_k tr(x_ky_k)| \pl:\pl
 \|\summ_k y_k^*\circ y_k\| \le 1\} \lel
  \sup\{ | \summ_k tr(y_kx_k)| \pl:\pl
 \|\summ_k y_ky_k^*\| \le 1\} \\
 &=
  \|(\summ_k x_k^*x_k)^{\frac12}\|_{L_1(N)} \pl.
 \end{align*}
This means that the \emph{column space  in
$S_1=L_1(B(\ell_2,tr))$ carries the operator space structure of
the row space in $B(\ell_2)$}.

At the end of these preliminaries we recall the definition  of the
operator space $OH$. This space is obtained from $C$ and $R$ by
complex interpolation. On the matrix level we have
 \[ \Mz_m(OH)\lel [\Mz_m(C),\Mz_m(R)]_{\frac12}  \pl .\]
Let $(e_k)$ be the canonical basis of $OH$. Then
 \[ \noo \summ_k x_k\ten e_k\rrm_{\Mz_m(OH)}
 \lel \noo \summ_k x_k\ten \bar{x}_k\rrm^{\frac12}_{\Mz_m \ten
 \bar{\Mz}_m} \pl. \]
Here $\overline{B(H)}=B(\bar{H})$. In many aspects the operator
space $OH$ is the appropriate analogue of $\ell_2$ in the category
of operator spaces. We refer to \cite{Po} for more details and
information. For general Hilbert spaces we use the notations
$H^c=B(\cz,H)$, $H^r=B(H,\cz)$ and $H^{oh}=[H^c,H^r]_{\frac12}$.

\section{The algebraic central limit theorem}

Due to the work of Hudson/Parthasarathy \cite{HP}, von Waldenfels
\cite{vWa} and Hegerfeldt \cite{Heg} it is well-known that the
central limit theorem applies for CCR and CAR relations. More
recently Speicher developed an approach which applies to the
$q$-commutation relations as well. We will apply Speicher's
approach using the data given by a von Neumann algebra $N$ and a
weight. We will first start with a normal faithful state $\psi$ on
$N$. For fixed $\nen$, we consider the $n$-fold tensor product
$N_n=N^{\ten_n}$ and the faithful normal state
$\phi_n=\psi^{\ten_n}$ on $N_n$. Let us denote by $\pi_k:N\to
N^{\ten_n}$ the homomorphism which sends $N$ into the $k$-th copy
 \[ \pi_k(x)\lel 1\ten \cdots \ten \underbrace{x}_{\scriptsize \mbox{k-th position}} \ten 1\ten \cdots \ten 1 \pl
 .\]
As an additional reservoir we consider a finite von Neumann
algebra $\Mz_{2^n}$ with normalized trace $\tau_n$ and selfadjont
contractions $v_1(n),...,v_n(n)$. Let $T$ be an additional
scaling factor. Then we consider the new random variable
 \[ u_n(x)\lel u_{n,T}(x) \lel \sqrt{\frac{T}{n}} \summ_{k=1}^n v_k(n) \ten \pi_k(x)
 \pl .\]
We need some further notation to formulate the central limit
procedure. Let $A\subset \{1,...,m\}$ and $x_1,...,x_m$ be
elements in $N$. Then we use
 \[ \psi_{A}[x_1,...,x_m] \lel \psi(\prod_{i\in A}^{\to}x_i) \]
for the evaluation of $\psi$ of the ordered  product
$\stackrel{\to}{\prod}_{i\in A}x_i$ arising from those indices
which are in $A$. Let $\si=\{A_1,...,A_r\}$ be a partition of
$\{1,...,m\}$. Then the multiplicative extension of $\psi$ is
defined as
 \[ \psi_{\si}[x_1,....,x_m] \lel \prod_{j=1}^r
 \psi_{A_j}[x_1,...,x_m] \pl .\]
We denote by $P(m)$ the collection of all partitions of
$\{1,...,m\}$ and by $P_2(m)$ the set of pair partitions. For an
element $(k_1,....,k_m)\in \{1,...,n\}^m$, we  use the notation
$(k_1,....,k_m)\le \si$ if
\[ k_i\lel k_l \Leftrightarrow \exists_{\p 1\le j\le r} \pl
 \{i,l\}\subset A_j   \pl.\]
A rather general form of the central limit theorem may be
formulated as follows. (The special cases discussed below are
well-known). Let us recall that for an ultrafilter $\U$ on an
index set $I$ the limit $\lim_{i,\U} a_i=a$ holds if for every
$\eps>0$ the set $\{i: |a_i-a|<\eps\}\in \U$.

\begin{lemma}\label{comb00} Let $\psi$ be a state. Let $(v_k(n))_{k=1,..n}$ be
contractions  such that the singleton condition is satisfied:
 \begin{equation} \label{sing}  \tau_n(v_{k_1}(n)\cdots v_{k_m}(n))=0
 \end{equation}
holds whenever one of the variables $k_1,...,k_m$ occurs only
once. Let $\U$ be a free ultrafilter on $\nz$ and
 \begin{equation} \label{gen-part}
  \beta(\si) \lel  \lim_{n,\U} n^{-\frac m2}
  \summ_{(k_1,....,k_m)\le \si} \tau_n(v_{k_1}(n)\cdots v_{k_m}(n)) \pl.
  \end{equation}
Let $x_1,...,x_m\in N$. Then
 \[ \lim_{n,\U} \tau_n\ten \psi^{\ten_n}(u_{n,T}(x_1)\cdots u_{n,T}(x_m))
  \lel \sum_{\si\in P_2(m)}   \beta(\si)  (T\psi)_{\si}[x_1,...,x_m] \pl
  .\]
\end{lemma}

\begin{proof} Let $\si=\{A_1,...,A_r\}$ be a partition of
$\{1,...,m\}$ and $\nen$. Now, we may develop  the terms
 \begin{align*}
 &\tau_n\ten \psi^{\ten_n}(u_{n,T}(x_1)\cdots u_{n,T}(x_m))\\
 &\pl = \kla \frac{T}{n}\mer^{\frac{m}{2}} \summ_{k_1,....,k_m=1}^n   \tau_n(v_{k_1}(n)\cdots v_{k_m}(n))
 \psi^{\ten_n}(\pi_{k_1}(x_1)\cdots \pi_{k_m}(x_m)) \\
 & \pl = \kla \frac{T}{n}\mer^{\frac{m}{2}}  \summ_{\si\in P(m)} \summ_{(k_1,....,k_m)\le \si}
 \tau_n(v_{k_1}(n)\cdots v_{k_m}(n))
 \psi_{\si}[x_1,...,x_m] \pl .
 \end{align*}
Let us fix a partition $\si$ and perform the limit for $n$ along
the ultrafilter $\U$. The singleton condition \eqref{sing} implies
that only partitions where all the subsets have cardinality bigger
than $2$ provide a non-trivial contribution. In particular, it
suffices to consider partitions with cardinality $|\si|\le m/2$.
Since the $v_k(n)$'s are assumed to be contractions, we deduce
from the Cauchy-Schwarz inequality that
 \begin{align}
 |\tau_n(v_{k_1}(n)\cdots v_{k_m}(n))|\! &\le \!  \tau_n(v_{k_1}(n)
 v_{k_1}(n)^*)^{\frac12} \! \prod_{j=2}^{m-1}
 \|v_{k_j}(n)\|_{\infty} \tau_n(v_{k_m}(n)^*
 v_{k_m}(n)^*)^{\frac12} \le  1 \label{cs00} \pl .
 \end{align}
Since there are at most $n (n-1) \cdots (n-r+1)$ many tuples
$(k_1,...,k_m)$ such that $(k_1,...,k_m)\le \si$, we deduce for
$r<\frac{m}{2}$ that
 \[ \lim_{n,\U} n^{-\frac{m}{2}} \summ_{(k_1,....,k_m)\le \si}
 \tau_n(v_{k_1}\cdots v_{k_m}) \lel 0 \pl .\]
Therefore only the pair partitions with $|\si|=\frac{m}{2}$
provide a non-trivial contribution $\beta(\si)$.
 \qd

\begin{rem} {\rm The additional parameter $T$ will be used later
for strictly semifinite weights. We refer to \cite{AB} for
representation of arbitrary selfadjoint random variables
represented with pair partitions.}
\end{rem}

In the next section it will be important to analyze a growth
condition for central limits.

\begin{cor}\label{growth} For $x\in N$  we use the length function
 \[ |x| \lel \max\{(T\psi(x^*x))^{\frac12},(T\psi(xx^*))^{\frac12},\noo
x\rrm_{\infty}\} \pl .\] Let $x_1,...,x_m\in N$. Then
 \[ |\tau_n\ten \psi^{\ten_n}(u_{n,T}(x_1)\cdots u_{n,T}(x_m))| \kl m^{\frac{m}{2}} \prod_{i=1}^m |x_i| \pl .\]
holds for all $\nen$.
\end{cor}

\begin{proof} According to the proof of Lemma  \ref{comb00}, we
have
 \begin{align*}
 &|\tau_n\ten \psi^{\ten_n}(u_n(x_1)\cdots u_n(x_m))| \kl \!\!
 \summ_{\si \in P_{ns}(m)} n^{-\frac m2} \!\!\!\!\!\!\!\!\summ_{(k_1,...,k_m)\le \si, |\si|\le \frac{m}{2}}
\!\! \!\!  \!\!  \!\!   |\tau_n(v_{k_1}(n)\cdots v_{k_m}(n))|
          (T\psi)_{\si}[x_1,...,x_m] \pl .
\end{align*}
Here $P_{ns}(m)$ stand for partitions not containing singletons.
For fixed $\si$, we deduce from \eqref{cs00} and the fact that
there are $n\cdots (n-|\si|+1)$ many tuples satisfying
$(k_1,...,k_m)\le \si$ that
 \[ n^{-\frac m2} \summ_{(k_1,...,k_m)\le \si} |\tau_n(v_{k_1}(n)\cdots
 v_{k_m}(n))|\kl \frac{n\cdots (n-|\si|+1)}{n^{\frac m2}} \kl 1 \pl. \]
The Cauchy-Schwarz inequality implies
 \[ |(T\psi)(y_1\cdots y_r)|\kl ((T\psi)(y_1y_1^*)^{\frac12}
 ((T\psi)(y_r^*\cdots y_2^*y_2\cdots y_r)^{\frac12} \kl
 \prod_{i=1}^r|y_i| \pl .\]
This implies $|(T\psi)_{A_j}[x_1,...,x_m]|\kl \prod_{i \in A_j}
|x_i|$ for $|A_j|\gl 2$. Thus we have
 \[ |(T\psi)_{\si}[x_1,...,x_m]|\kl \prod_{i=1}^m |x_i| \pl \]
for partitions without singletons. Combining these estimates we
get
 \[ |\tau_n\ten \psi^{\ten_n}(u_n(x_1)\cdots u_n(x_m))|\kl
 \summ_{\si \in P(m),|\si|\le \frac{m}{2}} \prod_{i=1}^m |x_i| \pl
 .\]
Finally, we note that there are not more than  $m^{\frac{m}{2}}$
partitions $\si$ with $|\si|\le \frac m2$. \qd

Before we discuss concrete examples let us indicate how these
estimates remain valid in the context of weights. We recall that a
function $\psi:N_+\to [0,\infty]$ is called an n.s.f. weight if
 \begin{enumerate}
  \item[n.)] $\psi(\sup_i x_i)\lel \sup_i \psi(x_i)$ holds for
  every increasing net $(x_i)$ of positive elements;
  \item[s.)] For every $0\le x$ one has $\psi(x)=\sup_{0\le y\le
  x,\psi(y)<\infty} \psi(y)$.
  \item[f.)] $\psi(x)=0$ implies $x=0$.
  \end{enumerate}
For a weight $\psi$ we may define $n_{\psi}=\{x\in N:
\psi(x^*x)<\infty\}$. We use the notation $n^*_{\psi}=\{x\in N:
\psi(xx^*)<\infty \}$ and $(x,y)_{\psi}=\psi(x^*y)$. The
completion of $n_{\psi}$ with respect to the inner product norm
is denoted by $L_2(N,\psi)$. It is well-known that $\psi$ extends
to a linear functional on $n_{\psi}^*n_{\psi}$, see \cite{Strat}
for details.

\begin{exam}\label{exssw}{\rm  Let $(f_j)$ be a net of mutually orthogonal projections which sum up to the identity, i.e.
$\sum_j f_j=1$. Let  $(\psi_j)$ be a family of  normal faithful
states on $f_jNf_j$.  Then $\psi(x)=\sum_j \psi(f_jxf_j)$ defines
an n.s.f. weight.}
\end{exam}

The example in \ref{exssw} is the prototype of a strictly
semifinite normal faithful weight. A weight $\psi$ is called
strictly semifinite if there exists an increasing net $(e_i)$ of
projections in the centralizer $N_{\psi}$ converging strongly to
$1$ and such that $\psi(e_i)<\infty$ holds for all $i$. We recall
that an element $x$ belongs to the centralizer $N_{\psi}$ if the
modular group $\si_t^{\psi}$ with respect to $\psi$ satisfies
$\si_t^{\psi}(x)=x$ for all $t\in \rz$. In the example above the
index set is the collection $P_{<\aleph_0}(I)$ of all finite
subsets of $I$ and the increasing family of projections is given
by $e_J=\sum_{j\in J} f_j$. We will use the following well-known
lemma for a strictly semifinite weight with invariant projections
$(e_i)$. The proof follows immediately from the fact that $xe_i$
and $e_ix$ converges in $L_2(N,\psi)$ to $x$. To show the
convergence of $\lim_i xe_i=x$, we use standard modular theory,
$S(x)=x^*$, $S=J\Delta^{\frac12}$ and observe that
 \begin{equation}\label{modtr} x^*e_i \lel (e_ix)^*\lel J\Delta^{\frac12}(e_ix)\lel Je_i\Delta^{\frac12}x\lel
 Je_iJx^*\end{equation}
holds for all analytic elements $x\in n_{\psi}\cap n^*_{\psi}$
and all $i$.

\begin{lemma}\label{mod-approx} Let $\psi$ be a strictly semifinite faithful normal
weight and $(e_j)$ as above. Then
 \[ \psi(x_1\cdots x_m) \lel \lim_i \psi(e_ix_1e_i\cdots
 e_ix_me_i) \]
holds for all $x_1,...,x_m\in n_{\psi}\cap n_{\psi}^*$.
\end{lemma}

%\begin{proof}
%Since $x_m\in n_{\psi}\subset L_2(n,\psi) $, we know that
%$L_2-\lim_j f_jx_k\cdots x_m=x_k\cdots x_m$ holds for all
%$k=1,..,m$. Let $\Delta$ be the generator of the modular group. We
%note that $\si_t^{\psi}(f_j)=f_j$ implies that $f_j$ is in the
%domain of $\Delta^{z}$ for every $z$ and $\Delta^z(f_j)=f_j$. We
%denote by $J$ the partial isometry such that
%$J\Delta^{\frac12}(a)=a^*$. With $f_j^*=f_j$ we deduce that
% \begin{equation}\label{modtr} a^*f_j \lel (f_ja)^*\lel J\Delta^{\frac12}(f_ja)\lel Jf_j\Delta^{\frac12}a\lel
% Jf_jJa^*\end{equation}
%holds for all analytic elements $a$. By density, we deduce that
%$af_j=Jf_jJa$ holds for all $a\in L_2(N,\psi)$. In particular,
%$L_2-\lim_j x_mf_j=x_m$. Let $\pi:N\to B(L_2(N,\psi))$ be the GNS
%representation. Using the strong convergence of $Jf_jJ$ and $f_j$
%to $1$, we deduce inductively that
% \begin{align*}
% &\lim_j f_jx_2\cdots f_jx_mf_j \lel  \lim_j\pi(f_jx_2\cdots
% f_jx_{m-1}f_j)(x_mf_j)  \lel \lim_j\pi(f_jx_2\cdots f_jx_{m-1}f_j)x_m \\
% &= \lim_j\pi(f_jx_2\cdots f_jx_{m-1})(f_jx_m) \lel \lim_j\pi(f_jx_2\cdots
% f_jx_{m-1})(x_m)  \\
% &=
% \lim_j\pi(f_jx_2\cdots f_jx_{m-2})(
% f_jx_{m-1}x_m) \lel  \cdots \lel x_2\cdots x_m \pl .
% \end{align*}
%Here all the limits are with respect to the $L_2(N,\psi)$-norm.
%Finally, $L_2-\lim_j x_1^*f_j=x_1^*$ and hence we get
%\begin{align*}
% \psi(x_1\cdots x_m) &= (x_1^*,x_2\cdots x_m)_{\psi}
% \lel \lim_j (x_1^*f_j,\pi(f_jx_2\cdots f_j)x_mf_j) \lel
%  \lim_j \psi(f_jx_1f_j \cdots f_jx_mf_j) \pl . \qedhere
% \end{align*}
% \qd

In this purely algebraic setting it is very convenient to use the
algebraic notion of a tensor algebra
 \[ A(V)=\sum_{n=0}^\infty V^{\ten_n} \]
of a vector space $V$. The formal multiplication is obtained by
linear extension of
\[ (v_1\ten \cdots \ten v_r)\times (v_{r+1}\ten \cdots \ten v_{m})\lel v_1\ten \cdots\ten v_r\ten v_{r+1}\ten \cdots \ten v_m\pl .\]
The words of length $0$, $V^{0}=\cz$ (or $V^{0}=\rz$), correspond
to the multiples of the identity.  In our context we apply this to
$V=n_{\psi}\cap n_{\psi}^*$.

\begin{cor}\label{weight}  Let $\psi$ be a strictly semifinite weight and $(e_i)_{i\in I}$
as above.  Let $\U'$ be a free ultrafilter on $I$ and
$x_1,...,x_m\in n_{\psi}\cap n_{\psi}^*$. Then the function
 \[ \phi(x_1\ten \cdots \ten x_m) \lel \lim_{j,\U'} \lim_{n,\U}
 \tau_n\ten
 (\frac{\psi}{\psi(e_j)})^{\ten_n}(u_{n,\psi(e_j)}(e_jx_1e_j)\cdots
 u_{n,\psi(e_j)}(e_jx_me_j)) \]
extends to a  linear functional on the tensor algebra
$A(n_{\psi}\cap n_{\psi}^*)$ and satisfies
 \begin{equation} \label{limpp}  \phi(x_1\ten \cdots \ten x_m) \lel \summ_{\si\in P_2(m)} \beta(\si)
  \psi_{\si}[x_1,...,x_m] \end{equation}
and for $|x|=\max\{\noo
x\rrm_{\infty},\psi(x^*x)^{\frac12},\psi(xx^*)^{\frac12}\}$
 \begin{equation} \label{limppp}  |\phi(x_1\ten \cdots \ten x_m)| \kl m^{\frac{m}{2}} \prod_{i=1}^m
 |x_i|\pl . \end{equation}
\end{cor}

\begin{proof} For fixed $j$ we use $T_j=\psi(e_j)$. Then
$\psi_j(x)= \frac{\psi(x)}{T_j}$ is a normal faithful  state on
$e_jNe_j$. Note that $T_j\psi_j$ is the restriction of $\psi$ to
$e_jNe_j$ and hence independent of $j$. According to Lemma
\ref{comb00} and Corollary \ref{growth} we know  that
\eqref{limpp} and \eqref{limppp} hold for all elements in
$\bigcup_j e_jNe_j$. We apply Lemma \ref{mod-approx} and obtain
that
 \[ \psi_{\si}[x_1,...,x_m]\lel \lim_j
 \psi_{\si}[e_jx_1e_j,...,e_jx_me_j]\]
for all $x_j\in n_{\psi}\cap n_{\psi}^*$. Similarly, we deduce
from $|x|= \max\{\|x\|,\psi(x^*x)^{1/2},\psi(xx^*)^{1/2}\}$ that
 \begin{align*}
  |x| &=  \lim_j
 \max\{\|e_jxe_j\|,T_j\psi_j((e_jxe_j)^*(e_jxe_j))^{\frac12},T_j\psi_j((e_jxe_j)(e_jxe_j)^*)^{\frac12}\}
 \pl .\end{align*}
Thus the assertion holds for arbitrary elements $x_1,..,x_m\in
n_{\psi}\cap n_{\psi}^*$.\qd
 We will now follow Speicher's approach for the $q$-commutation relations.
Let $(\Om,\mu)$ be a probability space and $s_{ij}:\Om \to
\{1,-1\}$ be random variables such that $s_{ij}=s_{ji}$ and such
that the family $(s_{ij})_{i<j}$ is independent. For fixed $\om
\in \Om$ we define the Pauli-matrices
 \[ v_{i,i}(\om) \lel \kla \begin{array}{cc}  0&1\\1&0 \end{array} \mer
 \quad \mbox{and for $i<j$}\quad
  v_{i,j}(\om) \lel \kla \begin{array}{cc}  1&0\\0 & s_{i,j}(\om)  \end{array} \mer
  \pl .\]
For fixed $\nen$, we use
 \[ v_j(\om)\lel v_{1,j}(\om)\ten \cdots \ten  v_{j,j}(\om)\ten  1
 \ten  \cdots \ten 1 \in \Mz_{2^n}\pl .\]
This yields selfadjoint unitaries such that
 \[ v_i(\om)v_j(\om)\lel s_{i,j}(\om) \pl v_j(\om)v_i(\om) \pl .\]
In the special case of a constant random variable  $s_{i,j}=-1$ we
obtain Clifford matrices. For a pair partition
$\si=\{A_1,...,A_{m/2}\}$ we write $(i,j)\in I(\si)$ if $i$ and
$j$ are part of an inversion of $\si$. This means $\{i,l\}\subset
A_{j_1}$ and $\{r,j\}\subset  A_{j_2}$ with $r<i<j<l$.

\begin{lemma} \label{comb01} The random variables $v_i\in
L_{\infty}(\Om,\Mz_{2^n})$ satisfy the singleton condition
\eqref{sing} with respect to $\tau_n(x)\lel \int
2^{-n}tr(x(\om))d\mu(\om)$. Moreover, for every pair partition
$\si\in P_2(m)$:
 \[  \beta(\si) \lel \lim_{n,\U} n^{-\frac{m}{2}}
 \summ_{(k_1,....,k_m)\le\si}
 \prod_{(i,j)\in I(\si)}\ez (s_{k_i,k_j}) \pl .\]
\end{lemma}

\begin{proof} Let $(k_1,...,k_m)$ be such that one index  occurs
only once. We fix $\om\in \Om $ and consider
 \[ 2^{-n}tr(v_{k_1}(\om)\cdots v_{k_m}(\om)) \pl.\]
Using the commutation relation
$v_{k_i}(\om)v_{k_j}(\om)=s_{k_i,k_j}(\om)v_{k_j}(\om)v_{k_i}(\om)$
and $v_{k_i}(\om)^2=1$ we may assume that all the variables occur
only once and $k_j$ is the largest index with this property. Then
we find some $v$ such that  $v_{k_1}(\om)\cdots
v_{k_m}(\om)=(-1)^{\pm 1} v\ten v_{k_j,k_j}(\om)\ten 1$. This
implies that $2^{-n}tr(v\ten v_{k_j,k_j}(\om)\ten 1)\lel 2^{-n}
tr(v)\pl tr(v_{k_j,k_j}(\om))=0$. For the second assertion, we
consider a pair partition $\si=\{A_1,....,A_{m/2}\}$. Let us
consider $(k_1,...,k_m)\le \si$ and $\om\in \Om$. We may assume
that $A_1=\{1,j\}$. Now, we commute $v_{k_j}(\om)$ with all the
unitaries $v_{k_{j-1}}(\om)$,...,$v_{k_{2}}(\om)$. This yields a
factor $\prod_{i=2}^{j-1}s_{k_i,k_j}(\om)$. If there is no
inversion between $i<j$ the index $k_i$ occurs twice. Therefore,
we get
 \[ \prodd_{i=2}^{j-1}s_{k_i,k_j}(\om) \lel  \prodd_{i,(i,j)\in I(\si)}
 s_{k_i,k_j}(\om)\pl .\]
After this procedure we use
$v_{k_1}(\om)v_{k_j}(\om)=v_{k_j}(\om)^2=1$. Thus we have
eliminated all the inversions with $j$ at the cost of the factor
$\prod_{i,(i,j)\in I(\si)}s_{k_i,k_j}(\om)$. Continuing by
induction, we find
 \[ 2^{-n}tr(v_{k_1}(\om)\cdots v_{k_m}(\om)) \lel \prod_{(i,j)\in I(\si)}
  s_{k_i,k_j}(\om)  \pl .\]
Taking  expectations yields the assertion by independence. \qd

The following result is due to Speicher:
\begin{cor}\label{prob}{\rm (}Speicher{\rm )}  Let ${\rm Prob}(s_{i,j}=1)=p$ and  ${\rm
Prob}(s_{i,j}=-1)=1-p$. Then
 \[ \beta(\si) \lel (2p-1)^{I(\si)} \pl ,\]
holds for all $\sin \in P_2(m)$. Here $I(\si)$ is the number of
inversions of $\si$.
\end{cor}

\begin{proof} This follows immediately from Lemma \ref{comb01}, $\ez s(k_i,k_j)=\al-(1-\al)\lel 2\al-1$
 and
 \begin{align*}
  \lim_n n^{-\frac{m}{2}} \summ_{(k_1,...,k_m)\le \si}
  &=  \lim_n n^{-\frac{m}{2}} n(n-1)\cdots (n-\frac{m}{2}) \lel 1
  \pl .\qedhere
  \end{align*}
\qd

Let us investigate the combinatorics including all values of
$-1\le q\le 1$ simultaneously. We use the random variable
 \[ s^q(t)\lel 1_{[-1,q]}-1_{[q,1]} \]
on $[-1,1]$ with respect to the Haar measure $\frac{dt}{2}$. Note
that $\ez s^q\lel \frac{1}{2}(2-2(1-q))=q$. Then $v_{j,q}$ is
constructed above using independent copies $(s^q_{i,j})_{i,j}$.

\begin{cor}\label{-one} Let $u_{n,q}(x)\lel \sqrt{\frac{T}{n}}\sum_{j=1}^n
v_{j,q}\ten \pi_j(x)$. Let $q_1,...,q_m\in [-1,1]$ and $\rho$ be
the partition given by sets where the $q_i$'s coincide. Then
 \begin{align*}
 \lim_n \tau_n\ten \kla\frac{\psi}{T}^{\ten_n}\mer (u_{n,q_1}(x_1)\cdots u_{n,q_m}(x_m))
 \lel \summ_{\si\in  P_2(m), \si \le \rho} t(\si,q_1,...,q_m)
 \psi_{\si}[x_1,...,x_m] \pl ,
 \end{align*}
where $t(\si)\lel ((m/2)!)^{-1}\sum_{\gamma \in \Pe(m/2)}
\prod_{(i,j)\in
I(\si)}[1_{i\le_{\gamma}j}q_j+1_{j\le_{\gamma}i}q_i]$. Here we use
the lexicographic order $(A_1,...,A_{|\si|})$ for the sets in
$\si$, and for a inversion $(i,j)\in I(\si)$ with $A_s=\{l,j\}$,
$A_t=\{i,k\}$ with $l<i<j<k$ we write $i\le_{\gamma}j$ if
$\gamma(t)<\gamma(s)$.
\end{cor}

\begin{proof} It suffices to analyze
 \[ \lim_n n^{-\frac{m}{2}}\sum_{(k_1,...,k_m)\le \si}
 \tau_n(v_{k_1,q_1}\cdots v_{k_m,q_m}) \pl . \]
As in Lemma \ref{comb01}, we find
 \[ \tau_n(v_{k_1,q_1}\cdots v_{k_m,q_m})
 \lel \prod_{(i,j)\in I(\si)} \eps(k_{i,q_i},k_{j,q_j})
 \tau_n(v_{k_1,q_{\pi(1)}}v_{k_1,q_{\pi(2)}}\cdots
 v_{k_{m/2},q_{\pi(m-1)}}v_{k_{m/2},q_{\pi(m)}}) \]
for a suitable permutation $\pi$. Here $I(\si)$ is the collection
of inversions of $\si$. The sign $\eps(k_{i,q_i},k_{j,q_j})$ is
calculated as follows. If $k_i<k_j$, then we know that
  \[
  \kla\begin{array}{cc} 0&1\\1&0\end{array}\mer
  \kla\begin{array}{cc} 1&0\\0&s_{k_i,k_j}^q\end{array}\mer
  \lel s_{k_i,k_j}^q \kla\begin{array}{cc} 1&0\\0&s_{k_i,k_j}^q\end{array}\mer
  \kla\begin{array}{cc} 0&1\\1&0\end{array}\mer \pl .\]
Thus we find
 \[ \eps(k_{i,q_i},k_{j,q_j})\lel 1_{k_i<k_j}q_j+1_{k_j<k_i}q_i
 \pl .\]
Thus the sign produces $q_i$ or $q_j$ depending on which color is
considered the bigger one. This leads to the description of
$t(\si,q_1,...,q_m)$ for $\si\le \rho$. If $\si\not\le \rho$, we
will find terms of the form
 \[  \frac{1}{n} |\summ_{j=k_{m/2-1}}^n [\frac12+\frac12 \ez
 s^{q_{m-1}}s^{q_{m}}]^{j-k_{m/2-1}}]|
 \kl \frac{1}{n|1-\al|} \]
 where $\al=\frac12+\frac12 \ez s^{q_{m-1}}s^{q_{m}}$ is in
$(-1,1)$. These terms disappear in the limit.\qd

\section{CAR, CCR  and $q$-commutation relations}

The CCR and CAR relations and their quasi-free states are
motivated by quantum mechanics. In this section we review
algebraic properties of the family of $q$-commutation relation and
the type of the associated von Neumann algebra. We will follow the
central limit approach from the previous section and use the data
provided by a weight on $N=L_{\infty}(\Om,\mu;\Ma_2)$. Let us
start with the discrete case. We consider a sequence of positive
numbers $(\mu_j)$ such that $0<\mu_j<1$ and the weight
 \[ \psi(x) \lel \sum_j[(1-\mu_j) x_{11}(j)+ \mu_j x_{22}(j)] \]
on $\ell_{\infty}(\nz,\Ma_2)$. Obviously  $\psi$ is a strictly
semifinite weight because the projections $e_j= \delta_{j}\ten
1_{\Mz_2}$ are invariant under the modular group. Here
$(\delta_j)$ denotes the sequence of unit vectors in
$\ell_{\infty}$ (where $\delta_j(k)=\delta_{jk}$ is given by the
Kronecker symbol).

For $-1\le q\le 1$, we may define a linear functional $\phi_q$ on
$A(n_\psi\cap n_{\psi}^*)$ by the formula
 \begin{equation}\label{q-rel} \phi_q(x_1\ten \cdots \ten x_m) \lel \summ_{\si\in P_2(m)}
 q^{I(\si)} \psi_{\si}[x_1,...,x_m] \pl .
 \end{equation}
We denote by $I_q$ the ideal generated by
 \[ \{z\in A(n_{\psi}\cap n_{\psi}^*) \pl:\pl \forall
 a,b\in A(n_{\psi}\cap n^*_{\psi}) \pl: \phi_q(a\ten z\ten b) \lel
 0 \}\]
Then, we may define $B_q=A(n_{\psi}\cap n_{\psi}^*)/I_q$. We
denote by $\pi_q$ the quotient homomorphism and observe that
$\phi_q$ induces a functional $\hat{\phi}_q$ on $B_q$ such that
$\hat{\phi}_q\pi_q=\phi_q$.

\begin{lemma}\label{carrel} Let $q=-1$ and $\al=0$. Let
$b_j\lel \delta_j\ten \kla \begin{array}{cc}
0&1\\0&0\end{array}\mer $. Then
 \[  \pi_{-1}(b_j\ten b_k+b_k\ten b_j)\lel 0
 \quad ,\quad \pi_{-1}(b_j\ten b_k^*+b_k^*\ten b_j)\lel
  \delta_{kj}  \]
and
 \begin{equation}\label{-1}
  \phi_{-1}(b_{i_r}^*\ten \cdots b_{i_1}^*\ten b_{j_1}\ten
 \cdots \ten b_{j_s}) \lel \delta_{r,s} \prod_{l=1}^r \delta_{i_l,j_l}\mu_{i_l}
  \end{equation}
for all increasing sequences  $i_1<i_2<\cdots <i_r$ and
$j_1<i_2<\cdots <j_s$.
\end{lemma}

\begin{proof} We first consider $y_1,...,y_m\in n_{\psi}\cap
n_{\psi}^*$ and deduce from Corollary \ref{prob} that
 \begin{equation}\label{car00}  \phi_{-1}(y_1\ten \cdots \ten y_m)
 \lel \summ_{\si\in P_2(m)} (-1)^{I(\si)} \psi_{\si}[y_1,...,y_m] \pl
 . \end{equation}
Now, we want to calculate the commutator between elements $y$ and
$z$:
 \begin{align*}
 &\phi_{-1}(y_1\ten \cdots \ten y_{r-1}\ten y\ten z \ten y_{r+1}\ten \cdots  \ten
 y_m) \\
 &=
 \summ_{\si\in P_2(m+1), \{r,r+1\}\in \si} (-1)^{I(\si)} \psi(yz)
 \psi_{\si}[y_1,\cdots
 y_{r-1},y_{r+1},...,y_{m}] \\
 &\pll  + \summ_{\si\in P_2(m+1), \{r,r+1\}\notin \si} (-1)^{I(\si)}
 \psi_{\si}[y_1,...,y_{r-1},y,z,y_{r+1}...,y_m] \pl.
 \end{align*}
For the second term we define a bijection $\rho$ on the partitions
$\si$ by exchanging $r$ and $r+1$. Since $\{r,r+1\}$ is not
contained in $\si$, we see that $\psi_{\si}[y_1,...,y,z,...y_m]=
\psi_{\rho(\si)}[y_1,...,y,z,...,y_m]$. However, $\rho$ changes
the number of inversion by one. This implies
 \begin{align*}
  \summ_{\si\in P_2(m), \{r,r+1\}\notin \si} (-1)^{I(\si)}
 \psi_{\si}[y_1,...,y,z,...,y_m] &=
 - \summ_{\si\in P_2(m), \{r,r+1\}\notin \si} (-1)^{I(\si)}
 \psi_{\si}[y_1,...,z,y,...,y_m] \pl.
 \end{align*}
By cancellation we obtain
 \begin{align*}
 &\phi_{-1}(y_1\ten \cdots \ten y_{r-1}\!\!\ten y\ten z \ten y_{r+1}\ten \cdots  \ten
 y_m) + \phi_{-1}(y_1\ten \cdots \ten y_{r-1}\!\!\ten z\ten y \ten
y_{r+1}\ten \cdots  \ten
 y_m) \\
  &= \psi(yz+zy) \phi_{-1}(y_1\ten \cdots y_{r-1}\ten y_{r+1}\ten
  \cdots y_{m}) \pl .
   \end{align*}
This shows
 \begin{equation}\label{carr}
  y\ten z+ z\ten y-\psi(yz+zy)1 \in I_{-1} \pl .
  \end{equation}
However, it is easily checked  that $\psi(b_jb_k+b_kb_j)=0$ for
all $j,k$ and  $\psi(b_jb_k^*+b_k^*b_j)=0$ for  $j\neq k$. For
$j=k$, we note that
  \[ \psi(b_k^*b_k+b_kb_k^*)\lel (1-\mu_k)+\mu_k\lel 1 \pl .\]
This completes the proof for the anticommutation relations.
Equation \eqref{-1} is a direct consequence of \eqref{car00}.
Indeed, given a pair partition $\si=\{A_1,....,A_{\frac{m}{2}}\}$,
$m=r+s$, we see that only  for $\{i_1,...,i_r\}=\{j_1,...,j_s\}$,
we have a chance to match up all the coefficients. Since the
subsequences  are assumed to be increasing, we deduce that only
the inversion-free partition $\{\{1,m\},\{2,m-1\},
,....,\{m/2,m/2+1\}\}$ leads to a nontrivial contribution.
 \qd

\begin{rem}\label{cont} {\rm In the continuous case we consider
weights of the form
 \[ \psi(x)\lel \intt [f_1(\om) x_{11}(\om)+ f_2(\om)x_{22}(\om)]
 \pl d\mu(\om)  \pl .\]
We assume in addition  that $f_1(\om)+f_2(\om)=1$. We use the
elements $b(g)(\om)\lel \kla
\begin{array}{cc}0& g(\om)\\0&0
\end{array}\mer$ and deduce from \eqref{carr} that
 \[ \pi_{-1}(b(g)\ten b(h)+b(h)\ten b(g))=0\pll ,\pll
  \pi_{-1}(b(g)\ten b(h)^*+b(h)^*\ten b(g))= (\intt \bar{h}g \mu ) \]
and
  \[ \phi_{-1}(b(g_r)^*\ten \cdots \ten b(g_1)^*\ten b(h_1)\cdots
  b(h_s))
  \lel \delta_{rs}  \det(\intt \bar{g_i}h_j f_2d\mu)  \pl.\]
}
\end{rem}

\begin{rem}\label{typecar} {\rm The type of the associated
factor satisfying these commutation relations can be determined
from the work of Araki-Woods \cite{AW}. Indeed, let us first
consider a selfadjoint operator $\rho$  on a real Hilbert space
$K$ such that $0\le \rho\le 1$. Let $H=K+iK$ be the
complexification and $\{u(f):f\in H\}$ a field satisfying the CAR
relations. We repeat the same construction for $K\oplus K$ and use
the well-known map
 \[ u_{\rho}(f)\lel u((1-\rho)^{\frac12}f,0)+
  u(0,\rho^{\frac12}f)^* \pl .\]
This is equivalent to \cite[(12.22)]{AW} (see also \cite{Rid}).
Indeed, one can check that
 \[ u_{\rho}(f)u_{\rho}(g)+u_{\rho}(g)u_{\rho}(f)\lel 0 \quad
 \mbox{and} \quad u_{\rho}(f)u_{\rho}(g)^*+u_{\rho}(g)^*u_{\rho}(f)\lel (f,g)1 \pl.  \]
Moreover, the moments with respect to the vacuum state $\xi$ of
the full system is given by
 \[ (\xi,u_{\rho}(f_r)^*\cdots u_{\rho}(f_1)^*u_{\rho}(g_1)\cdots
 u_{\rho}(g_s)\xi) \lel \delta_{rs} \det(e_j,\rho(g_i)) \pl.  \]
The corresponding von Neumann algebra $R(\rho)$ is defined by
 \[ R(\rho)=\{u_{\rho}(f):f\in K\}''\]
where the commutant is  taken in the GNS-representation with
respect to the vacuum state $\xi$. For discrete $\rho$ we see that
$R(\rho)$ is the ITPFI factor with $R(\rho)\cong
\ten_{\nen}(\Mz_2,\phi_n)$ where
$\phi_n(x)=(1-\la_n)x_{11}+\la_nx_{22}$ is determined by the
spectral values $\{\la_n:\nen\}$ of $\rho$. (This can also be
checked using the standard matrices satisfying the CAR relation
and by considering the GNS-construction with respect to the tensor
product state $\phi=\ten_{\nen}\phi_n$.) It follows from
\cite[section12]{AW} that for $\rho$ with  continuous spectrum,
$R(\rho)$ is the hyperfinite III$_1$ factor. If there exists a
$0<\la<1$ such that for every $\nen$ there exists $k_n,l_n$ such
that $\la_n\lel \la^{k_n}/(\la^{k_n}+\la^{l_n})$, then $R(\rho)$
is a type III$_{\la}$ factor (see \cite[section8]{AW}).}\end{rem}

The CCR relations can be obtained in a similar manner.

\begin{lemma}\label{ccrrel} Let $N=\ell_{\infty}(\Mz_2)$ and
$(\mu_k)$ as in Lemma \ref{carrel}. Let $X_k= \delta_k\ten \kla
\begin{array}{cc}0&1\\1&0\end{array}\mer$ and $Y_k=
\delta_k\ten \kla
\begin{array}{cc}0&i\\-i&0\end{array}\mer$. Then
 \[ \pi_1(X_k\ten X_j-X_j\ten X_k)\lel \delta_{kj} \pll \mbox{and}\pl
 \quad \pi_1(Y_k\ten Y_j-Y_j\ten Y_k)=\delta_{kj} \pl \]
and
 \[  \pi_1(X_k\ten Y_j-Y_j\ten X_k)\lel  \delta_{kj} \pl 2i(2\mu_k-1)
  \pl .\]
Moreover,
 \[ \phi_{1}(e^{zX_k}e^{wY_j}) \lel  e^{izw(2\mu_k-1)\delta_{kj}} e^{\frac{z^2+w^2}{2}} \pl. \]
\end{lemma}

\begin{proof} Let us first consider arbitrary elements
$y_1,...,y_{r-1}$, $y_{r+1},...,y_{m+1}$ and $y,z\in N$. Now, we
proceed as in Lemma \ref{carrel}. The only  difference is that if
we replace the partition $\si$ containing $\{i,r\}$ and
$\{r+1,j\}$ by the one containing $\{i,r+1\}$ and $\{j,r\}$ this
does not produce a sign change. Therefore
 \begin{align*}
 &\phi_1(y_1\ten \cdots y_{r-1}\ten y\ten z\ten y_{r+1}\ten
 \cdots y_m)- \phi_1(y_1\ten \cdots y_{r-1}\ten z\ten y\ten y_{r+1}\ten
 \cdots y_m)\\
 &= \psi(yz-zy)\phi_1(y_1\ten \cdots \ten y_{r-1}\ten y_{r+1}\ten
 \cdots y_{m+2}) \pl .
 \end{align*}
This implies
 \begin{equation}\label{ccr0} [\pi_1(y),\pi_1(z)]\lel \psi([y,z])
 \pl. \end{equation}
From
 \[ \psi([X_k,X_j]) \lel \delta_{kj}\pll,\pll
 \psi([Y_k,Y_j]) \lel \delta_{kj} \pll \mbox{and}\pll
 \psi([X_k,Y_j]) \lel  \delta_{kj} 2i (2\mu_k-1) \pl \]
we deduce the first three equations. For the last assertion
 we use the power series expansion. For a fixed even $k$ we
have
 \[ \phi_1(y^{\ten_k})\lel \summ_{\si \in P_2(k)}
 \psi(y^{2})^{\frac{k}{2}}  \pl. \]
Let $g$ be a $N(0,1)$ distributed normal random variable and $w\in
\cz$ such that $w^2=\psi(y^2)$. Then
 \[ \ez (wg)^k\lel w^k \summ_{\si \in P_2(k)} 1\lel
 \phi_1(y^{\ten_k}) \pl.\]
This implies that
 \begin{align}\label{ccrmom}
 \phi_1(e^{y})&= \summ_{k=0}^{\infty}
 \frac{\phi_1(y^{\ten_k})}{k!} \lel \summ_{k=0}^{\infty} \ez
 \frac{(wg)^k}{k!} \lel
  \ez e^{wg} \lel e^{\frac{w^2}{2}} \lel e^{\frac{\psi(y^2)}{2}}
 \pl .
 \end{align}
If $k\neq j$ then $\pi_1(X_k)$ and $\pi_1(Y_j)$ commute. This
implies with \eqref{ccrmom}
 \begin{align*}
  \phi_1(e^{zX_k}e^{wY_j})&= \phi_1(e^{zX_k+wY_j})\lel
  e^{\frac{\psi((zX_k+wY_j)^2)}{2}} \lel e^{\frac{z^2+w^2}{2}} \pl
  . \end{align*}
For $k=j$, we use the Baker-Cambell-Hausdorff formula (see e.g.
\cite[Section 3.3]{OL}) and obtain from \eqref{ccr0} that
 \begin{equation}\label{Ba-ca-ha}  e^{t(zX_k+wY_k)}   \lel e^{-\frac{t^2zw\psi([X_k,Y_k])}{2}}
 e^{tzX_k}e^{twY_k} \pl .\end{equation}
holds as a formal power series in $t$ with values in $B_1$. On the
other hand, we deduce from Corollary \ref{growth} that
 \[ \phi_1(e^{tzX_k}e^{twY_k}) \lel \summ_{n,m} t^{n+m}
 \frac{z^nw^m}{n!m!} \phi_1(X_k^{\ten_n}Y_k^{\ten_m}) \]
is absolutely  converging because
$\phi_1(X_k^{\ten_n}Y_k^{\ten_m})\kl |X_k|^n|Y_k|^m$. Similarly,
we deduce that the map $t\mapsto \phi_1(e^{tzX_k}e^{twY_k})$ is an
absolutely converging power series. This implies
 \begin{equation}\label{CBH}  \phi_1(e^{t(zX_k+wY_k)})  \lel e^{-\frac{t^2zw\psi([X_k,Y_k])}{2}}
 \phi_1(e^{tzX_k}e^{twY_k}) \pl \end{equation}
for all $t\in \rz$. For $t=1$ we deduce the assertion from
$\psi([X_k,Y_k])=2i(2\mu_k-1)$. \qd

\begin{rem}\label{contccr} {\rm In the continuous case we consider
the weight $\psi(x)=\int [f_1x_{11}+f_2x_{22}]d\mu$ on
$N=L_{\infty}(\Om,\mu;\Mz_2)$. The real Hilbert space $H_{\rz}$
is given by functions $f:\Om \to \cz$. Then, we define the real
linear embedding
 \[ j(f)\lel \kla \begin{array}{cc} 0& f\cr  \bar{f}&0
 \end{array}\mer \pl . \]
We obtain that
 \[ [\pi_1(j(f)),\pi_1(j(g))] \lel \psi([j(f),j(g)])\lel \intt (f\bar{g}-g\bar{f})(f_1-f_2)
 d\mu
 \pl . \]
By  the Baker-Cambell-Hausdorff formula \eqref{Ba-ca-ha} (valid as
a for formal power series in the algebra $B_1$) this implies that
 \begin{equation}\label{bch-key}
  \phi_1(a \ten e^{zj(f)}e^{wj(g)}\ten b)
 \lel \phi_1(a\ten e^{-zw \intt (f\bar{g}-g\bar{f})(f_1-f_2)d\mu}
e^{zj(f)+wj(g)}\ten b)
  \end{equation}
Here $a,b$ is either $1$ or a finite tensor $a=x_1\ten \cdots \ten
x_m$, $b=y_1\ten \cdots \ten y_m$. Note that in both cases the
absolute convergence of the corresponding power series follows
from the growth condition. Hence we may evaluate them at $t=1$.
Moreover, we deduce from $\psi(j(f)^2)=\int |f|^2d\mu$ that
 \[ \phi_1(e^{zj(f)})\lel e^{\frac{z^2}{2}\int |f|^2d\mu} \pl .\]
} \end{rem}

The reader familiar with the classical CCR relations will have
observed that the formulas are not the usual ones. Let us first
recall these usual CCR relations and state some results on the
corresponding type of the underlying von Neumann algebra. We will
follow the representation in \cite{AW1}. Let $K$ be a real Hilbert
space. The Weyl representation is given by a collection of
unitaries $\{U(f)|f\in K\}$, $\{V(f)|f\in K\}$ such that
 \begin{equation}\label{ccr1}
  U(f+g)\lel U(f)U(g) \quad,\quad V(f+g)\lel V(f)V(g)
  \end{equation}
and
 \begin{equation}\label{ccr2} U(f)V(g)\lel e^{-i(f,g)}V(g)U(f) \pl .
 \end{equation}
The corresponding quasi free state $\phi$ is determined by the
relation
 \begin{equation}\label{ccrstate}
    \phi(U(f)V(g))\lel e^{-i\frac{(f,g)}{2}-\frac{\noo f\rrm^2+\noo
  g\rrm^2}{4}}\pl .
 \end{equation}
Given real subspaces $K_1$ and $K_2\subset K$ Araki introduced
 \[ R(K_1,K_2/K)\lel \{U(f)V(g)\pl:\pl f\in K_1, g\in
 K_2\}'' \pl .\]
Here the bicommutant is taken in the GNS representation of
$\{U(f),V(g):f,g\in K\}$ with respect to $\phi$. In
\cite[section12]{AW} we can find very precise information about
the type of $R(K_1,K_2/K)$. We shall assume in addition that
$K=K_1\oplus K_1^{\perp}$ and that $K_2$ is the graph of an
unbounded operator $B$. If $B^*B$ has discrete spectrum, then
$R(K_1,K_2/K)$ may be constructed as an ITPFI factor. Let $\rho$
be the operator  such that $B^*B=4\rho(1+\rho)$. We refer to
\cite[section 12, section 8]{AW} for the following result.

\begin{theorem}[Araki-Woods]\label{typccr} If $\rho$ has a continuous spectrum
then $R(K_1,K_2/K)$ is the hyperfinite {\rm III}$_1$ factor. If
the spectrum $\{\la_n; n\in \nz\}$  of $\rho$ is discrete then
$R(K_1,K_2/K)$ is isomorphic to $\ten_{\nen} (B(\ell_2),\phi_n)$
where $\phi_n(x)=tr(D_nx)$  and $D_n$ has the eigenvalues
$\{\mu_n^k(1-\mu_n): k=0,1,...\}$ and
$\mu_n=\frac{\la_n}{1+\la_n}$. Moreover, if  $0<\lambda<1$ such
that  $\mu_n=\lambda^{k_n}$ for some $k_n\in \nz$. Then
$R(K_1,K_2/K)$ is the type {\rm III}$_{\lambda}$ factor.
\end{theorem}

We will translate our algebraic relations into the usual CCR
relations. For a measurable function $f:\Om\to \rz$ we define
 \[
  U(tf) \lel   e^{it\pi_1(j(2^{-\frac12}f))} \quad \mbox{and} \quad V(tg)\lel e^{it\pi_1(j(2^{-\frac12}ig))}
 \pl \]
as formal power series in $B_1$. The relations $U(f+g)=U(f)U(g)$,
$V(f+g)=V(f)V(g)$ are obvious. We will  show in Theorem
\ref{limitth} that this may be interpreted as unitary operators on
a Hilbert space. According to Remark \ref{contccr} we have
 \begin{align*}
  \phi_1(a\ten U(tf)V(tg)\ten b) &= e^{-it^2\int fg (f_2-f_1)d\mu}
  \phi_1(a\ten V(tg)U(tf) \ten b)\pl .
 \end{align*}
Similarly, we have
 \[ \phi_1(U(f)V(g))\lel e^{-\frac{i}{2}\int fg (f_2-f_1)d\mu} e^{-\frac{1}{4}(\intt
 f^2+g^2 d\mu)} \pl .\]
We will now match this with the usual  CCR-relations. We introduce
$K=L_2(\mu;\rz)\oplus L_2(\mu;\rz)$ and $K_1\lel \{(f,0): f\in
L_2(\mu;\rz)\}$. We recall the assumption $f_1+f_2=1$. Thus the
operator  $A(f)\lel (f_2-f_1)f$ is a contraction on $L_2(\mu;\rz)$
and  we may define
 \begin{equation} \label{ktwo}  K_2\lel \{(Ag,\sqrt{1-A^2}(g))\pl:\pl g \in L_2(\mu,\rz)\} \pl
 . \end{equation}
Then $K_2$ is the graph of  the operator $B=A^{-1}\sqrt{1-A^2}$.
The map $g\mapsto (Ag,\sqrt{1-A^2}(g))$ is an isometry and
 \[ ((f,0),(Ag,\sqrt{1-A^2}(g)) \lel (f,Ag)\lel \intt fg
 (f_2-f_1)d\mu \pl .\]
In view of Theorem \ref{typccr}, we  now solve the equation
$B^*B=4(\rho(1+\rho))$ for $B=A^{-1}\sqrt{1-A^2}$ and find
$2\rho=1+|A^{-1}|$. As an application, let us state the following
result.

\begin{cor}\label{typeIII1} Let $K_1$ and $K_2$ as above. If the
operator $(1+|A|^{-1})(f)=(1+|f_2-f_1|^{-1})f$ has continuous
spectrum, then $R(K_1,K_2/K)$ is the hyperfinite III$_1$ factor.
\end{cor}

%\begin{rem}{\rm For a general weight $\psi$ it is still possible
%to understand the type of CCR-relations given by the family
%$e^{i\pi_1(x)}$ where  $x$ runs through the selfadjoint elements
%$x\in n_{\psi}$ which are orthogonal to the centralizer (with
%respect to the scalar product on $L_2(N,\psi)$). Indeed, if
%$2\Delta/(\Delta-1)$ has continuous spectrum and $N_*$ is
%separable we find the tensor product of the hyperfinite $II_1$
%factor (coming from the centralizer $N_{\psi}$)  and the
%hyperfinite $III_1$ factor.}
%\end{rem}

At the end of this section, we will discuss the unifying approach
introduced by Shlyakhtenko which allows us to describe the
$q$-commutation relations for all $-1\le q\le 1$. Here we shall
assume that $U_t$ is a one parameter unitary group on a real
Hilbert space $K$.  On the complexification $H=K+iK$ one may
write $U_t=A^{it}$ for a suitable positive nonsingular generator
$A$. The crucial information is contained in the new scalar
product
 \[ (x,y)_U\lel (2A(1+A)^{-1}x,y) \]
on $K+iK$. We assume the scalar product to be antilinear in the
first component. We denote by $H_U$ the completion of $H$ with
respect to this scalar product. Indeed, the following properties
characterize $(K,H_U,A)$.
 \begin{enumerate}
 \item[i)] $H\supset K$ is a complex Hilbert space and $(U_t)$ a one parameter
 group such that $U_t=A^{it}$ and $U_t(K)\subset K$,
 \item[ii)] $K+iK$ is dense in $H$ and $K\cap iK=\{0\}$,
 \item[iii)] The restriction of the real part of the scalar product on
 $H$ induces the scalar product on $K$.
 \item[iv)]  ${\rm Im}(x,y)_H\lel (i\frac{1-A^{-1}}{1+A^{-1}}x,y)_K$ for
 all $x,y\in K$.
 \end{enumerate}
Moreover, on a dense domain one has
 \begin{equation}\label{modinv}
  (x,y)_U\lel (y,A^{-1}x)_U \pl .
  \end{equation}
We refer to \cite{S} for more information on these conditions and
the fact that these properties characterize the inclusion
$K\subset H$. The algebraic information of the factor
$\Gamma_q(K,U_t)$ is contained in a generating family of bounded
operators $\{s_q(h):h\in K\}$ and the vacuum vector $\xi$
satisfying
 \begin{equation} \label{q-com} (\xi,s_q(h_1)\cdots s_q(h_m)\xi)\lel
 \summ_{\si=\{\{i_1,j_1\},...,\{i_{\frac{m}{2}},j_{\frac{m}{2}}\}\}
 \in P_2(m)}
 q^{I(\si)} \prod_{l=1}^{\frac{m}{2}} (h_{i_l},h_{j_l})_U
 \end{equation}
for even $m$. For odd $m$ we have $0$.

\begin{exam} \label{modgr} {\rm
Given a strictly semifinite weight $\psi$ on a von Neumann algebra
$N$, we may consider the algebra $B_q(n_{\psi}\cap n_{\psi}^*)$
and the scalar product given by
  \[ (x,y)_{\psi} \lel \psi(x^*y) \pl .\]
Then $\phi_q$ defines a functional on $A(n_{\psi}\cap n_{\psi}^*)$
satisfying
 \[ \phi_q(x_1\ten \cdots \ten x_m) \lel \summ_{\si \in P_2(m)}
 q^{I(\si)} \psi_{\si}[x_1,....,x_m] \]
which coincides with \eqref{q-com} for selfadjoint elements.
Indeed, we consider the real Hilbert space  $K$ of selfadjoint
elements in $n_{\psi}\cap n_{\psi}^*$ and the complex Hilbert
space $L_2(N,\psi)$. On $K$ we shall use the real scalar product
 \[ (x,y)_K\lel \frac{\psi(xy)+\psi(yx)}{2} \pl .\]
The modular group $\si_t^{\psi}$ is a one parameter group of
unitaries on $K$. From \eqref{modinv} we see that $A=\Delta^{-1}$
is a generator such that $\psi(xy)\lel \psi(yA^{-1}x)$ holds for
analytic elements. Thus we consider $U_t=\si_{-t}$ which leaves
$K$ invariant. For analytic selfadjoint elements  $x,y$ we have $
(x,y)_K = 1/2\psi((1+\Delta^{-1})(x)y)$. This implies
 \begin{align*}
  &(i\frac{1-A^{-1}}{1+A^{-1}}(x),y)_K \lel  (i\frac{1-\Delta}{1+\Delta}(x),y)_K
  \lel   (i\frac{\Delta^{-1}-1}{1+\Delta^{-1}}(x),y)_K \\
  &=
   \frac12 \psi(
  (1+\Delta^{-1})(i\frac{\Delta^{-1}-1}{1+\Delta^{-1}}(x))y) \lel
   \frac{i}{2} [\psi(\Delta^{-1}(x)y)-\psi(xy)] \\
   &=   \frac{i}{2} [\psi(yx)-\psi(xy)]
   \lel {\rm Im}(\psi(xy)) \pl .
  \end{align*}
This yields
\[ (x,y)_U\lel \psi(xy) \pl \]
for selfadjoint elements. Thus the formula \eqref{q-rel} extends
the combinatorial formula \eqref{q-com} for arbitrary elements in
$n_{\psi}\cap n_{\psi}^*$. In the literature these extensions are
obtained by considering the complex linear extension $\hat{s}_q$
(see \cite{S} and \cite{Hiai}) of $s_q$.}
\end{exam}

\begin{exam}\label{22-ext} {\rm Let us consider the special case
$N=L_\infty(\mu;\Mz_2)$ and
 \[ \psi(x)\lel \intt [f_1(\om)x_{11}(\om)+ f_2(\om)x_{22}(\om)]
 d\mu(\om) \pl. \]
For $f\in L_2(\mu;\cz)$ we define $j(f)\lel \kla \begin{array}{cc}
0& f\\ \bar{f}&0\end{array}\mer$. It is easily checked that $j$ is
an isometric embedding and
$\Delta^{-it}(j(f))=(\frac{f_2}{f_1})^{it}f$ leaves
$K=L_2(\mu;\cz)$ invariant. Since $j$ is only real linear, we
shall consider $K$ as a real Hilbert space. The subspace
$H=j(K)+ij(K)$ is spanned by elements of the form $j_2(f,g)=\kla
\begin{array}{cc} 0& f\\g& 0\end{array}\mer$. By restriction we find
\[ {\rm Re}\psi(j(f)j(g))\lel {\rm Re}\int
[f\bar{g}f_1+\bar{f}gf_2]d\mu \lel \int
\frac{f\bar{g}+\bar{f}g}{2}
 d\mu \pl .\]
For the imaginary part we have
 \begin{align*}
 {\rm Im}\psi(j(f)j(g)) &= {\rm Im}(\intt [f\bar{g}f_1+\bar{f}gf_2]d\mu)
 \lel \intt [\frac{f\bar{g}-\bar{f}g}{2i} f_1+
 \frac{\bar{f}g-f\bar{g}}{2i}f_2] d\mu \\
 &= \int
 \frac{i(f_2-f_1)f\bar{g}- i(f_2-f_1)\bar{f}g}{2} d\mu
 \lel \int
 \frac{i(f_2-f_1)f\bar{g}+\overline{i(f_2-f_1)f}g}{2} d\mu \\
 &= {\rm Re}\psi(j(i(f_2-f_1)f)j(g)) \pl .
 \end{align*}
This shows that on $H$ the operator $A$ is given by $j_2(f,g)=
j_2(f_2/f_1f,f_1/f_2 g)$ and $C=i\frac{1-A^{-1}}{1+A^{-1}}$ is
given by $C(j(f))=j(i(f_2-f_1)f)$.}
\end{exam}

We may consider Example \ref{22-ext} as the prototype for
arbitrary inclusions using a well-known argument of Araki
\cite{Ar}.

\begin{lemma}\label{arak-form} Every triple $(K,H_U,A)$ is unitarily equivalent to
\[ (K_0^{\rz},K_0^{\cz},id) \oplus \big(L_2(\Om,\mu;\cz),L_2(\Om,\mu;\sqrt{f_2}\cz) \oplus L_2(\Om,\mu;\sqrt{f_1}\cz),
(M_{f_2/f_1},M_{f_1/f_2})\big)
 \pl. \]
Here $U_t$ acts as the identity on the real Hilbert space
$K_0^{\rz}$ and $K_0^{\cz}=K_0^{\rz}+iK_0^{\rz}$ is the canonical
complexification. On the orthogonal complement
$K_1=L_2(\Om,\mu;\cz)$ the embedding is given by
$j(f)=(f,\bar{f})$ and $U_t(f)=(f_2/f_1)^{it}f$ holds for suitable
positive functions $f_1$, $f_2$ with $f_1+f_2=1$.
\end{lemma}

\begin{proof} We follow \cite[Section2]{S} and \cite{Ar}. On $K\oplus iK$ there
exists a selfadjoint operator $S$ such that $U_t=e^{itS}$. Being a
generator of a one parameter unitary group, we know that $B=iS$
is an unbounded operator on $K$. We may assume that $S$ is
injective. Indeed, there is a projection $e$ corresponding to the
kernel of $S$ and we might write $K=K_0+K_1$ such that
$eU_t=U_te=eU_te=e$ and $U_t$ acts non-trivially on $K_1$. On
$K_0$ the inclusion $K_0\subset H_U$ is given by the canonical
inclusion $K_0\subset K_0+iK_0$. Thus for the rest of the proof
we may assume $K=K_1$. Note that $B^{t}B=S^*S$ is a selfadjoint
real operator and thus we have a polar decomposition $B=V|S|$
(coming from $S=U|S|$ and $iU=V$) such that $V$ is a unitary on
$K$. Since $S$ is selfadjoint we have
 \[ 1\lel (S|S|^{-1})^2 \lel U^2 \lel (-iV)^2 \pl.  \]
This implies  $V^2=-1$. On the other hand we have $V^{t}V=1$ and
hence $V^{t}=-V$. Moreover, $S^*=S$ implies $B^{t}=-B$. Therefore
we get
 \[- V|S| \lel -B\lel  B^{t} \lel |S|V^{t}\lel-|S|V \pl .\]
Thus $V$ commutes with $S$. Using a maximal system of cyclic
vectors, one can construct $K_{+}$ such that $K_{+}+V(K_{+})=K$
and $|S|(K_+)\subset K_+$. Using a maximal system of cyclic
vectors again and functional calculus, we see that $|S|_{K_+}$ is
unitarily equivalent to a  multiplication operator $M_{f}$ on
$L_2(\Om,\mu;\rz)$ (see \cite{kar-I}). By orthogonality of $K_{+}$
and $K_{-}=V(K_{+})$, we deduce that $(K,B)$ is unitarily
equivalent to $(L_2(\mu;\cz),iM_{f})$ such that $K=K_{+}\oplus
K_{-}$ is real isomorphic to $L_2(\mu;\cz)$ and $U_t=e^{itf}$.
Now, we may use our standard model and set
$f_1=e^{-f}(1+e^{-f})^{-1}$, $f_2=(1+e^{-f})^{-1}$, and
 \[ \psi(x)\lel \intt [f_1 x_{11} + f_2x_{22}] d\mu    \]
According to Example \ref{22-ext} the corresponding action is
given by $(f_2/f_1)^{it}\lel e^{itf}=U_t$. Moreover,
$M_{i(f_2-f_1)}$ also provides the correct imaginary part for
$i\frac{1-A^{-1}}{1+A^{-1}}$. \qd

As a direct application of \cite[Theorem 3.3]{Hiai} we obtain the
following information on  types.

\begin{cor}\label{hia} Let $K=L_2(\mu;\cz)$ and $U_t(f)=(f_2/f_1)^{it}f$.
Let $-1<q<1$ and  $\Gamma_q(K,U_t)$ the $q$-deformed Araki-Woods
factor. Let $G$ be the closed subgroup generated by the spectrum
of $f_1/f_2$. Then $\Gamma_q(K,U_t)$ is of type {\rm II}$_1$ if
$G=\{1\}$, of type {\rm III}$_{\la}$ if $G=\{\la^n:n\in \zz\}$,
and of type {\rm III}$_1$ if $G=\rz_+$.
\end{cor}

 \section{Limit distributions}

In the previous section, we have seen an algebraic construction of
the CAR and CCR relations from the central limit theorem. In this
part we will discuss some analytic properties which will be used
to provide matrix models. We will first discuss ultraproducts of
noncommutative $L_p$ spaces.

Given an index set $I$ and a family  $(X_n)_{n\in I}$ of Banach
spaces, we will use the notation $\prod_b X_n$ for the subspace of
bounded sequences in $\prod X_n$.
%For an ultrafilter $\U$ on $\nz$
%and a bounded family  $(a_n)$ of complex numbers, we write
%$\lim_{n,\U} a_n=z$ if
% \[ \{ n\in \nz\pl:\pl |a_n-z|<\eps\}\in \U\]
%for every $\eps>0$.
The ultraproduct $\prod_{n,\U}X_n$ is defined
as the Banach space $\prod_b X_n/J_\U$, where $J_\U$ is the closed
subspace of $\prod_bX_n$ given by
 \[ J_\U\lel \{(x_n)\pl|\pl \lim_{n,\U} \|x_n\| \lel 0 \} \pl . \]
We will follow \cite{Ra} and use the notation
 \[ (x_n)^{\bullet} \lel (x_n)+J_\U \pl .\]
Let $(N_n)$ be a sequence of von Neumann algebras. Then, we may
consider the ultraproduct $\prod_{n,\U} L_1(N_n)$. Following the
natural duality for operator spaces (see \cite{4-author}), we
find an isometric embedding
 \[ \iota:\prodd_{n,\U} L_1(N_n)\to (\prodd_b N_n^{op})^* \quad \mbox{given
 by}\quad
  \iota(d_n)(x_n)\lel \lim_{n,\U} tr_{N_n}(d_nx_n) \pl .\]
Note  $\iota$ is a completely isometric embedding and the range
$V=\iota(\prod_{n,\U} L_1(N_n))$ is left and right invariant by
multiplication with elements $x\in \prod_b N_n$. Therefore (see
\cite{Tak}), we find a central projection $z_{\U}\in (\prod_n
N_n^{op})^{**}$ such that $V=z_\U(\prod_n N_n^{op})^{*}$. We
obtain a von Neumann algebra
 \[ \tilde{N}_\U \lel z_\U(\prodd_n N_n^{op})^{**} \pl .\]
The von Neumann algebra $\tilde{N}_\U$ is the dual of the
ultraproduct of preduals. The map $\iota$ provides a completely
isometric isomorphism between $\prod_{n,\U} L_1(N_n)$ and
$L_1(\tilde{N}_\U)\cong (\tilde{N}_\U^{op})_*$. This argument
(except for the additional complication with $op$ needed for the
operator space level) is due to Groh \cite{Groh}. For our
applications  $\tilde{N}_\U$ is still too big because, in general,
it is not $\si$-finite. Therefore, we will in addition assume that
$(\phi_n)$ is a sequence of faithful normal states such that
$\phi_n$ is defined on $N_n$. Let $(D_n)$ be the corresponding
sequence of densities associated to $(\phi_n)$. \emph{This
notation will be fixed in the sequel}. Clearly, $(D_n)^{\bullet}$
is an element in $\prod_{n,\U} L_1(N_n)$.  We  denote by $e_\U$
its support (in $\tilde{N}_\U$). Then, we define the \emph{von
Neumann algebra ultraproduct with respect to $(\phi_n)$} by
 \[ N_\U\lel e_{\U}\tilde{N}e_\U \pl. \]
If all the states $\phi_n$ are tracial, then $N_\U$ is the
well-known von Neumann algebra ultraproduct $\prod_n N_n/J_\U$.
Here $J_\U$ is the ideal of elements $(x_n)$ such that
$\lim_{n,\U} \phi_n(x_n^*x_n)=0$. If the $\phi_n$ are only states,
then  $J_\U\subset \tilde{N}_\U(1-e_\U)$ is only a left ideal.
Let us note that we have a $^*$-homomorphism $\pi:\prod_b N_n\to
\tilde{N}_\U$ which has strongly dense range. It will be
convenient to use the notation $(x_n)^{\bullet}=\pi((x_n))$. Let
us denote by $D_{\U}\in L_1(\tilde{N}_\U)$ the density of the
state $\iota((D_n)^{\bullet})$. The following result is due to
Raynaud \cite{Ra}.

\begin{lemma}\label{ra} Let $0<p<\infty$. There is an (completely) isometric embedding $I_p:L_p(N_\U)\to \prod_{n,\U}
L_p(N_n)$ satisfying
 \[ I_p(D_{\U}^{\frac{1-\theta}{p}}e_{\U}(x_n)^{\bullet} e_\U
 D_{\U}^{\frac{\theta}{p}}) \lel
 (D_n^{\frac{1-\theta}{p}}x_nD_n^{\frac{\theta}{p}})^{\bullet} \pl .\]
Moreover, the subspace $W=\{e_\U(x_n)^{\bullet}e_\U\pl|\pl
(x_n)\in \prod_b N_n\}$ is strongly dense in $N_\U$. If
$\si_t^{\phi_n}$ and $\si_t^{\phi_\U}$ denote the modular group of
$\phi_n$ and $\phi_\U$, respectively, then
 \begin{equation} \label{modular1}  \si_t^{\phi_\U}(e_\U(x_n)^{\bullet}e_\U) \lel
 e_\U (\si_t^{\phi_n}(x_n))^{\bullet}e_\U \pl .
 \end{equation}
\end{lemma}

The classical central limit theorem provides unbounded random
variables. If we work with unbounded operators in a non-tracial
setting, we have to multiply unbounded operators without the help
of the algebra of $\tau$-measurable operators. The context of
$L_p$ spaces yields a very convenient substitute in our setting.
Let $D$ be the density of a normal faithful state $\phi$ on a von
Neumann algebra $N$ and $x$ be a selfadjoint operator affiliated
to $N$. We will say that $xD^{\frac1p}\in L_p(N)$ if there exists
$y\in L_p(N)$ such that for any increasing sequence of intervals
$(I_k)$ with $\bigcup_k I_k=\rz$ we have
 \[ y\lel \lim_k 1_{I_{k}}(x)xD^{\frac1p} \pl. \]
The next Lemma is an easy application of Kosaki's interpolation
result.

\begin{lemma}\label{assoc} Let $x$ be a selfadjoint operator affiliated to $N$,
$\phi$ a normal faithful state and $\mu(A)\lel \phi(1_{A}(x))$
the induced measure on $\rz$. If $2\le p\le \infty$ and  the
identity function $f(t)=t$ belongs to $L_p(\mu)$, then
 \[ xD^{\frac1p}\in L_p(N) \pl .\]
In particular, let $p=2m$. If $\int_{\rz} t^{2m} d\mu(t)$ is
finite, then $xD^{\frac1p}$ for all $p\le m$.
\end{lemma}

\begin{proof} Let $M$ be the commutative von Neumann algebra generated by $1$
and the spectral projections $e_A=1_{A}(x)$ where $A$ ranges
through the borel-measurable sets. Functional calculus provides a
natural $^*$-homomorphism $\pi:L_\infty(\rz,\mu)\to N$. For a
bounded function $f\in L_{\infty}(\mu)$ we define
 \[ \pi_p(f)\lel \pi(f)D^{\frac1p} \pl. \]
Let us recall  Kosaki's interpolation result \cite{Kos}
 \[ L_p(N)D^{\frac{1}{2}-\frac{1}{p}}\lel
 [ND^{\frac{1}{2}},L_2(N)]_{\frac{2}{p}} \pl .\]
Here the inclusion map  $\iota$ of $N$ in $L_2(N)$ is given by the
map $\iota(x)= xD^{\frac12}$. We observe that
$\pi_p(f)D^{\frac12-\frac1p}\lel \iota(\pi(f))$. Thus, taking the
inclusions into account, we see that the family $(\pi_p)_{2\le
p\le \infty}$ is indeed induced by the `same' operator
$\pi_2(f)=\pi(f)D^{\frac12}$. Note also that $\pi_p$ is a isometry
for $p=\infty$ and $p=2$ (recall
$\|f\|_{L_2(\mu)}=[\phi(\pi(|f|^2))]^{1/2}$). By interpolation, we
deduce that $\pi_p$ extends to a contraction on $L_p(\mu)$.
However, by the dominated convergence theorem, we see that for
every increasing family of bounded intervals $(I_k)$ and $f\in
L_p(\mu)$ we have
 \[ L_p(\mu)-\lim_k 1_{I_k}f\lel f \pl .\]
By continuity, we deduce that $\lim_k \pi(1_{I_k}(x)x)D^{\frac1p}$
converges to $\pi_p(f)$. This proves the first assertion. For
$p=2m$, we simply note that $\intt t^{2m}d\mu(t)$ means that
$f(t)=t$ belongs to $L_{2m}(\mu)$. Since the inclusion map
$I_{2m,p}:L_{2m}(N)\to L_p(N)$ defined by
$I_{2m,p}(x)=xD^{\frac{1}{p}-\frac{1}{2m}}$ is continuous, we
obtain $xD^{\frac1p}\in L_p(N)$ for all $p\le m$.\qd

As usual we denote by $\si_t^{\phi}$ the modular group of a state
(or weight) $\phi$. An element $x\in N$ is called {\it analytic}
if the function $t\mapsto \si_t^{\phi}(x)$ extends to an analytic
function with values in $N$. In this case we use the notation
$\si_z(x)\in N$. We denote  by {\it $N_{a}\subset N$ the
subalgebra of analytic elements}. The map $\si_z:N_a\to N$ is a
homomorphism satisfying $\si_z(x)^*\lel \si_{\bar{z}}(x^*)$. We
denote by $N_{sa,a}$ the {\it real subalgebra of selfadjoint
analytic elements}.

\begin{prop}\label{limit1} Let $X$ be an index set and $|\pll |:X\to [0,\infty)$ be a
function. Let $(N_n,\phi_n)$ be a family  of von Neumann algebras
$N_n$ with normal faithful states $\phi_n$. For every $\nen$ let
$(u_n(x))_{x\in X}$ be a family of analytic selfadjoint elements
in $N_n$ such that
 \begin{enumerate}
  \item[i)] There exists a constant $C$ such that
  \[ |\phi_n(u_n(x_1)\cdots u_n(x_m))|\kl (Cm)^m |x_1|\cdots
  |x_m| \]
holds for all $\nen$ and $x_1,...,x_m\in X$.
  \item[ii)] There exists a $2< p <\infty$ such that for all $x\in
  X$ there exists a constant $c(x,k)$ such that
   \[ \sup_n \noo
   D_n^{\frac1p}\si_{\frac{i}{2}}^{\phi_n}(u_n(x)^k) \rrm_p
   \kl c(x,k)
  \pl .\]
 \end{enumerate}
Let $\al:\bigcup_{m\ge 0} X^m\to \prod_{n,\U} L_2(N_n)$ be defined
by $\al(\emptyset)\lel (D_n^{\frac12})^{\bullet}$ and
 \[ \al(x_1,...,x_m)\lel \big(u_n(x_1)\cdots
 u_n(x_m)D_n^{\frac12}\big)^{\bullet} \pl .\]
Then there exists a family of selfadjoint operators $(u(x))_{x\in
X}$ affiliated to $N_\U$ such that
 \begin{enumerate}
 \item[a)] If $|z|<\delta(x,x_1,...,x_m)$, then
 $\sum_{k=0}^\infty \frac{z^k}{k!}u(x)^k\al(x_1,...,x_m)$
 converges absolutely in the ultraproduct $\prod_{n,\U}L_2(N_n)$,
 and
  \item[b)]  $\frac{d}{i\, dt} e^{itu(x)}\al(x_1,....,x_m)|_{t=0}\lel
 \al(x,x_1,....,x_m)$, and
 \item[c)] $((D_n^{\frac12})^{\bullet},\al(x_1,....,x_m))\lel \lim_{n,\U}
 \phi_n(u_n(x_1)\cdots u_n(x_m))$.
 \end{enumerate}
\end{prop}

\begin{proof} For fixed $x\in X$, we see that
$U_t^x=(e^{itu_n(x)})^{\bullet}$ is a unitary group in
$\tilde{N}_\U$. Let us show that $\al(x_1,...,x_m)$ is in the
domain of the generator. Let us fix $\nen$ and define
$\al_n(x_1,...,x_m)=u_n(x_1)\cdots u_n(x_m)D_n^{\frac12}$.
According to our assumption
 \[ \noo \al_n(x_1,...,x_m)\rrm_2^2 \lel \phi_n(u_n(x_m)\cdots u_n(x_1)u_n(x_1)\cdots
 u_n(x_m))\kl (2Cm)^{2m} \prod_{i=1}^m |x_i|^2 \]
is uniformly bounded and hence $\al(x_1,...,x_m)$ is an element in
$L_2(\tilde{N}_{\U})\cong \prod_{n,\U}L_2(N_n)$. Similarly,  we
have
 \begin{align}
  &\noo u_n(x)^k \al_n(x_1,...,x_m)\rrm_2^2
   \lel  \phi_n(u_n(x_m)\cdots u_n(x_1)u_n(x)^{2k}u_n(x_1)\cdots
  u_n(x_m))  \nonumber \\
  &\le (2C(m+k))^{2(m+k)} \pl \prodd_{i=1}^m |x_i|^2  \pl |x|^{2k}
  \kl  (2Cm)^{2m} \prodd_{i=1}^m |x_i|^2 \pl (2Cmk)^{2k}|x|^{2k}
  \label{momentum}
  \pl .
 \end{align}
Let $\delta=(4eCm (1+|x|))^{-1}$, $c_m= (2Cm)^{m} \prod_{i=1}^m
|x_i|$, and   $|z|\le \delta$. Then we see  that
 \begin{align} \label{powerseries}
 \noo \summ_k \frac{z^k}{k!}u_n(x)^k \al_n(x_1,...,x_m)\rrm_2
  &\le  c_m \summ_{k=0}^\infty \frac{|z|^k}{k!}
 \frac{k^k \delta^{-k}}{(2e)^k} \kl \summ_{k\ge 0} 2^{-k} \pl .
 \end{align}
Since these estimates are uniformly in $n$, we deduce assertion a)
and that
 \[ \frac{d}{i\, dt}U_t(\al(x_1,..,x_m))|_{t=0}\lel \al(x,x_1,...,x_m)   \]
is well-defined. Therefore $\al(x_1,...,x_m)$ is in the domain of
the closed generator $u(x)$ of the unitary group $U_t^x$ (see
\cite[Theorem 5.6.36]{Kad}). Let us note in passing that also
$(D_n^{1/2})^{\bullet}\tilde{N}_\U$ is in the domain of $u(x)$.
We recall that $e_\U$ is the support of
$\phi_\U=\iota((D_n)^{\bullet})$. Now, we will use condition ii)
to show that for $|t|$ small enough, we have $U_t^xe_\U=e_\U
U_t^x= U_t^xe_\U$. We will first show that for $|t|\le \delta(x)$
we have
 \begin{equation} \label{eq1}  (1-e_\U)U_t^x\bigg(D_n^{\frac12}y_n\bigg)^{\bullet} \lel 0
 \end{equation}
for every bounded sequence $(y_n)\subset \prod_b N_n$. Again, we
fix $\nen$ and note
 \[ u_n(x)^kD_n^{\frac12} \lel
 D_n^{\frac12}\si_{\frac{i}{2}}^{\phi_n} (u_n(x)^k)
 \lel D_n^{\frac{1}{2}-\frac1p} D_n^{\frac1p} \si_{\frac{i}{2}}^{\phi_n}
 (u_n(x)^k) \pl .\]
By definition of the multiplication in $L_2(\tilde{N}_\U)\cong
\prod_{n,\U} L_2(N_n)$ (see \cite{Ra}), we deduce
 \begin{equation}\label{ddom} \bigg(u_n(x)^kD_n^{\frac12}\bigg)^{\bullet} \lel
 \bigg(D_n^{\frac{1}{2}-\frac{1}{p}}\bigg)^{\bullet}
 \bigg(D_n^{\frac1p}\si_{\frac{i}{2}}^{\phi_n}(u_n(x)^k)\bigg)^{\bullet} \pl
 .\end{equation}
According to assumption ii)
$\big(D_n^{\frac1p}\si_{\frac{i}{2}}^{\phi_n}(u_n(x))^k\big)^{\bullet}\in
\prod_{n,\U}L_p(N_n)$ is well-defined and henceforth
 \[ (1-e_\U) \bigg(u_n(x)^kD_n^{\frac12}y_n\bigg)^{\bullet} \lel
 (1-e_\U) \bigg(D_n^{\frac12-\frac1p}\bigg)^{\bullet}
 \bigg(D_n^{\frac1p}\si_{\frac{i}{2}}^{\phi_n}(u_n(x)^k)\bigg)^{\bullet}
 \bigg(y_n\bigg)^{\bullet} \lel 0 \pl .\]
However, for $|t|\le \delta(x)$, we deduce from absolute
convergence that
 \[ (1-e_\U)U_t^x\bigg(D_n^{\frac12}y_n\bigg)^{\bullet}
  \lel \summ_{k=0}^\infty \frac{t^k}{k!} (1-e_\U)
  \bigg(u_n(x)^kD_n^{\frac12}y_n\bigg)^{\bullet}
   \lel 0 \pl .\]
This shows \eqref{eq1}. By density we find $(1-e_\U) U_t^x
(D_n^{\frac12})^{\bullet}\tilde{N}_\U \lel 0$, i.e.
$(1-e_\U)U_te_\U=0$ for $|t|\le \delta$. Using adjoints, we deduce
$e_\U U_t^x\lel e_\U U_t^xe_\U \lel U_t^xe_\U$. For arbitrary $t$
we choose $m$ such that $|t|\le m\delta$ and observe
 \[ e_\U U_t^x\lel e_\U (U_{\frac{t}{m}}^x)^m \lel
 (U_{\frac{t}{m}}^x)^me_\U \lel U_t^xe_\U \pl .\]
Therefore $e_\U U_t^xe_\U$ is a strongly continuous unitary group
in $N_\U$ and the generator $u(x)$  satisfies the assertion. \qd

\begin{rem}\label{amom}{\rm
{\bf 1)} The proof of Proposition \ref{limit1} allows us to
replace the bounded operators $u_n(x)$ by unbounded selfadjoint
operators. In this case, we should require that
$\al_n(x_1,....,x_m)$ is in the domain of $u_n(x)$ and replace the
moment condition by
 \[ |(D_n^{\frac12},u_n(x_1)\cdots u_n(x_m)D_n^{\frac12})| \kl
 (Cm)^m \prod_{i=1}^m |x_i| \pl .\]
In order to formulate the modular condition, we recall that
$\si_t^{\phi_n}(x)=D_n^{it}xD_n^{-it}$ also defines an
automorphism of $L_p(N_n)$. Therefore, we shall require that
$u_n(x)D_n^{1/p}$ is analytic and that for
$\frac{1}{2}=\frac{1}{p}+\frac{1}{q}$
  \[ \sup_n \noo \si_{\frac{i}{q}}(u_n(x)^kD_n^{\frac1p})\rrm_p \kl
  c(x,k) \pl .\]
Then we have
 \[ u_n(x)^k D_n^{\frac12} \lel
 D_n^{\frac12-\frac1p}\si_{\frac{i}{q}}(u_n(x)^kD_n^{\frac1p}) \pl .\]
This allows us to deduce  \eqref{ddom} and complete the proof as
above.

{\bf 2)} For applications, it is often  more convenient to replace
condition  ii) in Proposition \ref{limit1} by  a moment condition.
Let us assume that $p=2m$ is an even integer and
$\frac12=\frac1p+\frac1q$.  We consider the real and imaginary
part
 \[ y_{n,k}  \lel \si_{\frac{i}{q}}^{\phi_n}(u_n(x)^k)+
 \si_{\frac{i}{q}}^{\phi_n}(u_n(x)^k)^* \quad ,\quad
  z_{n,k}  \lel \si_{\frac{i}{q}}^{\phi_n}(u_n(x)^k)-
  \si_{\frac{i}{q}}^{\phi_n}(u_n(x)^k)^* \pl .\]
Then, we deduce from Lemma \ref{assoc} that
 \[ \noo y_{n,k}D_n^{\frac1p}\rrm_p^p \kl \phi_n(y_{n,k}^p) \kl
 \summ_{s_1,...,s_p\in \{\emptyset,*\}}
 |\phi_n(\si_{\frac{i}{q}}(u_n(x)^k)^{s_1}\cdots
 \si_{\frac{i}{q}}(u_n(x)^k)^{s_p})| \pl .\]
A similar estimate holds for $z_{n,k}$. This shows that moment
estimate
 \begin{equation} \label{analyticmoment}  \sup_{n,s_1,...,s_p\in \{\emptyset,*\}} |\phi_n(\si_{\frac{i}{q}}(u_n(x)^k)^{s_1}\cdots
 \si_{\frac{i}{q}}^{\phi_n}(u_n(x)^k)^{s_p})| \kl c(x,k,p) \end{equation}
implies condition ii). Again this holds still in the context of
unbounded operators as long as we can justify the domain issues
for the  operators $\si_{\frac{i}{q}}(u_n(x)^k)$ and justify the
equation
\[ \si_{\frac{i}{q}}^{\phi_n}(u_n(x)^k)D_n^{\frac1p}=
\si_{\frac{i}{q}}^{\phi_n}\big( u_n(x)^kD_n^{\frac1p}\big)\pl. \]
}
\end{rem}

At the end of this section we show that the limit objects in the
algebraic central limit theorem may be realized as unbounded
operators.

\begin{theorem}\label{limitth}  Let $N$ be a von Neumann algebra and $\psi$ a
strictly semifinite weight with associated projections $(e_j)$.
Let $(v_{k}(n))_{k=1,...,n,\nen}$ be a family of selfadjoint
contractions such that $(v_{k}(n))_{k=1,..,n}\subset
(\Mz_{2^n},\tau_n)$ are contained in a finite von Neumann algebra
satisfying the singleton condition \eqref{sing}. Let $\U$ be an
ultrafilter and for a pair partition $\si\in P_2(m)$ the weight
is given by
 \[ \beta(\si) \lel\lim_{n,\U}  n^{-\frac m2}\summ_{(k_1,...,k_m)\le \si}
 \tau_n(v_{k_1}(n)\cdots v_{k_m}(n)) \pl .\]
Consider $X=N_{a}\cap n_{\psi}\cap n_{\psi}^*$ and
$|x|=\max\{\|x\|_{\infty},\psi(x^*x)^{\frac12},\psi(xx^*)^{\frac12}\}$.
Then there exists a von Neumann algebra $M$, a normal faithful
state $\phi$ with density $D$, a linear map $\al:A(X)\to L_2(M)$
and a family $\{u(x)\}_{x\in N_{sa,a}}$ of selfadjoint operators
affiliated to $M$ such that
 \begin{enumerate}
 \item[i)] $\al(x_1\ten \cdots \ten x_m)$ is in the domain of $u(x)$ and
 \[ u(x)\al(x_1\ten \cdots x_m)\lel \al(x\ten x_1\ten \cdots \ten x_m)
 \pl ,\]
 \item[ii)] $D^{\frac12}$ is in the domain of $u(x)$ and
  \[ u(x)D^{\frac12}\lel \al(x) \pl ,\]
 \item[iii)] $(D^{\frac12},\al(x_1\ten \cdots \ten x_m))\lel
 \summ_{\si\in P_2(m)} \beta(\si)\psi_{\si}[x_1,...,x_m]$,
 \item[iv)] $\si_t^{\phi}(u(x))\lel u(\si_t^{\psi}(x))$ holds for
 all $t\in \rz$,
 \end{enumerate}
holds for all $x_1,...,x_m\in N_{sa,a}$ and $x\in N_{sa,a}$.
\end{theorem}

\begin{proof} For fixed  $j$ and $\nen$
we define $T_j=\psi(e_j)$ and $\phi_{n,j}=\tau_n\ten \kla
\frac{\psi}{T_j}^{\ten_n}\mer$ on $M_{n,j}=\Mz_{2^n}\ten
(e_jNe_j)^{\ten n}$. For $x\in N_{sa,a}$ we define
 \[ u_{n,j}(x)\lel \sqrt{\frac{T_j}{n}} \summ_{k=1}^n v_k(n)\ten
 \pi_k(e_jxe_j) \pl .\]
Since $\tau_n$ is assumed to be a trace, we have
 \begin{equation}\label{modcondd} \si_t^{\phi_n}(u_{n,j}(x))\lel u_{n,j}(\si_t^{\psi}(x)) \pl
 . \end{equation}
In order to check condition i) we apply Corollary \ref{growth} for
$e_jxe_j$. Let us show that $|e_jxe_j|\kl |x|$. Indeed,
$\|e_jxe_j\|_{\infty}\le \|x\|_{\infty}$. Moreover, we deduce from
\eqref{modtr} that
 \begin{align*}
  \psi(e_jx^*e_jxe_j) &\le \psi(e_jx^*xe_j)\lel \noo
  xe_j\rrm_{L_2(N,\psi)}^2 \lel \noo Je_jJx\rrm_{L_2(N,\psi)} \kl
  \noo x\rrm_{L_2(N,\psi)}^2 \lel \psi(x^*x) \pl .
  \end{align*}
Let us denote by $D_{n,j}$ the density of $\phi_{n,j}$. We will
now apply Remark \ref{amom} 2) in order to verify the condition
ii) in Proposition \ref{limit1}. Indeed, we take  $p=4=q$ and
$x_j=e_jxe_j$. Then we consider
 \[ y_{kj} \lel \si_{\frac{i}{q}}(x_j^k)+ \si_{\frac{i}{q}}(x_j^k)^* \quad
 \mbox{and} \quad z_{kj} \lel
 \si_{\frac{i}{q}}(x_j^k)-\si_{\frac{i}{q}}(x_j^k)^* \pl .\]
We note that
\[ |y_{kj}| \kl 2|\si_{\frac{i}{q}}(x_j^k)|\kl 2
 \|\si_{\frac{i}{q}}(x_j)\|_{\infty}^{k-2} |\si_{\frac{i}{q}}(x_j)|
 \pl .\]
However $|\si_{\frac{i}{q}}(x_j)|=|e_j\si_{\frac{i}{q}}(x)e_j|\kl
|\si_{\frac{i}{q}}(x)|$ implies
\[ |y_{kj}| \kl 2\pl \|\si_{\frac{i}{q}}(x)\|_{\infty}^{k-2}
 \pl  |\si_{\frac{i}{q}}(x)|\pl .\]
For $k=1$ we simply have $|y_{kj}| \kl 2|\si_{\frac{i}{q}}(x)|$.
Similarly, we find the estimate
 \[ |z_{kj}| \le 2 \pl\|\si_{\frac{i}{q}}(x)\|_{\infty}^{k-2}\pl
  |\si_{\frac{i}{q}}(x)| \pl .\]
This implies with Lemma \ref{assoc} and Lemma \ref{growth} that
 \begin{align*}
   \noo \si_{\frac{i}{q}}(x_j^k)D_n^{\frac14}\rrm_4&\le
  \frac{\noo y_{kj}D^{\frac14}\rrm_4+  \noo
  z_{k,j}D^{\frac14}\rrm_4}{2}\kl
 \frac{\phi_n(y_{kj}^4)^{\frac14}+
  \phi_n(z_{kj}^4)^{\frac14}}{2} \\
  &\le
   4^2 \frac{|y_{kj}|+|z_{kj}|}{2} \kl 32 \pl
  \|\si_{\frac{i}{q}}(x)\|_{\infty}^{k-2} \pl  |\si_{\frac{i}{q}}(x)|
  \pl .
  \end{align*}
Therefore the condition ii) in Proposition \ref{limit1} is
satisfied uniformly in $n$ and $j$.  Let $\U'$ be an ultrafilter
on $I$ and $\tilde{\U}=\U'\times \U$ the ultrafilter on $I\times
\nz$ such that $A\in \tilde{\U}$ iff $\{j:\{\nen: (j,n)\in A\}\in
\U\}\}\in \U'$. In terms of limits, this means
$\lim_{(j,n),\tilde{\U}}a_{j,n}=\lim_{j,\U'}\lim_{n,\U}a_{j,n}$.
We apply Proposition \ref{limit1} to the double indexed family
$u_{n,j}(x)$. By linearity it suffices to define the linear map
$\al$ on $A(X)$ for tensors with the help of Raynaud's
isomorphism:
 \[ \al(x_1\ten \cdots \ten x_m)\lel I_2^{-1}\bigg((u_{n,j}(x_1)\cdots
 u_{n,j}(x_m)D_{n,j}^{\frac 12})^{\bullet}\bigg) \pl .\]
The state $\phi_{\tilde{U}}$ is the ultraproduct state
$(\phi_{n,j})^{\bullet}$ with density $D_{\tilde{\U}}$ in
$L_1(N_{\tilde{\U}})$. Note that by definition of $N_{\tilde{U}}$,
the state $\phi_{\tilde{U}}$ is faithful. We deduce from
Proposition \ref{limit1} that $u(x)$ is affiliated to
$N_{\tilde{\U}}$ and satisfies the domain properties i) and ii).
The modular condition iv) follows from \eqref{modcondd} and
\eqref{modular1}. The moment condition iii) follows from Raynaud's
isomorphism, Lemma \ref{ra} and Corollary \ref{weight}
 \begin{align*}
  &(D_{\U},\al(x_1\ten \cdots \ten x_m))
   \lel  \lim_{(j,n),\tilde{\U}} (D_{n,j}^{\frac12},u_{n,j}(x_1)\cdots
  u_{n,j}(x_m)D_{n,}^{\frac12}) \\
  &=
  \lim_{j,\U'}\lim_{n,\U} \phi_{n,j}(u_{n,j}(x_1)\cdots
  u_{n,j}(x_m))
   \lel  \lim_{j,\U'} \summ_{\si\in
  P_2(m)}\beta(\si)\psi_{\si}[e_jx_1e_j,...,e_jx_me_j] \\
  &=
     \summ_{\si\in
  P_2(m)}\beta(\si)\psi_{\si}[x_1,...,x_m] \pl .\qedhere
  \end{align*} \qd

\begin{cor}\label{all} Let $N$ be a von Neumann algebra and $\psi:N\to \cz$ be a strictly semifinite normal faithful
weight. Then there exists a von Neumann algebra $M[-1,1]$ with a
normal faithful state $\phi$ and a family $\{U_t(x,q): x\in
N_{sa,a}, -1\le q\le 1\}$ of unitary groups which generate $M$
with the following properties.
 \begin{enumerate}
 \item[i)] The generators  $u_q(x)=\frac{d}{idt} U_t(x,q)$ satisfy
 \begin{align*}
 \phi(u_{q_1}(x_1)\cdots u_{q_m}(x_m))
 \lel \summ_{\si\in  P_2(m), \si \le \rho} t(\si,q_1,...,q_m)
 \psi_{\si}[x_1,...,x_m] \pl .
 \end{align*}
 \item[ii)] For every subset  $I\subset [-1,1]$ there is
 $\phi$-preserving conditional expectation $E_I:M\to M_I$, $M_I$
generated by $\{U_t(x,q):x\in N_{sa,a}, q\in I\}$. For disjoint
sets $I$ and $J$ the algebras $M_I$ and $M_J$ are independent over
$\phi$ {\rm (}in the sense of \eqref{indepo}{\rm)}
\end{enumerate}
\end{cor}

\begin{proof} The first assertion is a direct consequence of
Corollary \ref{-one} and  Theorem \ref{limitth} but now applied
for the family $v_{j,q}(n)$ indexed by the additional parameter
$-1\le q\le 1$. For a subset $I\subset [-1,1]$, we deduce from
$\si_t^{\phi_n}(e^{isu_{n,q}(x)})=e^{isu_{n,q}(\si_t^{\psi}(x))}$
that $M_I$ is invariant under $\si_t^{\phi}$ and hence we can
apply Takesaki's theorem (see e.g. \cite[Theorem 10.1]{Strat}) and
find a $\phi$-preserving conditional expectation. Now, we assume
$I\cap J=\emptyset$, $q_1,...,q_k\in I$, $q_{k+1},...,q_{m}\in J$.
Then
 \[ \phi(u_{q_1}(x_1) \cdots u_{q_k}(x_k)u_{q_{k+1}}(x_{k+1}) \cdots
 u_{q_m}(x_m))= \phi(u_{q_1}(x_1) \cdots u_{q_k}(x_k))\phi(u_{q_{k+1}}(x_{k+1}) \cdots
 u_{q_m}(x_m))\]
follows from i). Using absolute convergence (uniform in $q$), we
may extend this relation to polynomials in $U_t(x,q_k)$. Since
$\phi$ is normal this implies $\phi(ab)=\phi(a)\phi(b)$ for $a\in
M_{J}$, $b\in M_I$.  \qd

We refer to \cite{Bla} for a short introduction to
$C(X)$-algebras. In our context we obtain the following continuity
result (without constructing an embedding of $C[-1,1]$ in the
center.)

\begin{cor} Let $T=\{(x,t):x\in N_{sa,a}, t\in \rz\}$ and $A(T)$
be the free algebra in $T$ noncommuting selfadjoint variables. Let
$\pi_q:A(T)\to M$, $\pi_q((x,t))=e^{itu_q(x)}$ the induced
representation. The function
 \[ f_x(q) \lel \|x+\ker \pi_q\| \]
is lower semi-continuous.
\end{cor}

\begin{proof} Let $p$ be a noncommutative polynomial in the
variables $(x_1,t_1),...,(x_m,t_m)$. For fixed $q_0$ we can find
noncommutative polynomials $p_{1},p_2$ in the variables
$y_1,...,y_k\in N_{sa,a}$ such that
 \[ (1-\eps)\|\pi_{q_0}(p)\|\kl
 |(p_1(u_{q_0}(y_1),...,u_{q_0}(y_k)),\pi_{q_0}(p)p_2(u_{q_0}(y_1),...,u_{q_0}(y_k))|
 \pl \]
and
$\|p_1(u_{q_0}(y_1),...,u_{q_0}(y_k))\|_2=1=\|p_1(u_{q_0}(y_1),...,u_{q_0}(y_k))\|_2$.
Using the combinatorial formula from Corollary \ref{comb00} and
Corollary \ref{prob}, we deduce
 \[ \lim_{q\to q_0} \|p_{j}(u_{q}(y_1),...,u_{q}(y_k))\|_2 \lel
   \|p_{j}(u_{q_0}(y_1),...,u_{q_0}(y_k))\|_2 \pl \]
for $j\in \{1,2\}$. By linearity it suffices to show
 \begin{align*}
 &\lim_{q\to q_0}
  \phi(u_q(y_1)\cdots u_{q}(y_k)e^{it_1u_q(x_1)}\cdots
  e^{it_mu_q(x_m)}u_q(y'_{1})\cdots u_q(y'_{k'})) \\
  &\lel  \phi(u_{q_0}(y_1)\cdots u_{q_0}(y_k)e^{it_1u_{q_0}(x_1)}\cdots
  e^{it_mu_{q_0}(x_m)}u_{q_0}(y'_{1})\cdots u_{q_0}(y'_{k'})) \pl
  .
  \end{align*}
According to Corollary \ref{growth}, we know that
 \begin{align*}
 & |\phi(u_q(y_1)\cdots u_q(y_k)u_q(x_1)^{l_1}\cdots
 u_q(x_m)^{l_m}u_q(y'_1)\cdots u_q(y'_{k'}))|\\
  &\le  C(y_1,...,y_k,y'_1,...,y'_{k'}) (k+k')^{\frac{d}{2}}d^{\frac{d}{2}}\prod
 |x_j|^{l_j}
 \end{align*}
where $d=\sum_i l_i$. Using the absolute convergence of  $\sum_l
\frac{a^l}{l!} l^{\frac{l}{2}}$, we may now approximate the
exponential functions by polynomials (uniformly in $q$). Hence the
convergence follows from continuity in $q$ of the combinatorial
formula for the joint moments (see Corollary \ref{comb00} and
Corollary \ref{prob}). For an arbitrary  element of $A(T)$ the
assertion follows by approximation. \qd

\section{A uniqueness result}

\setcounter{equation}{0}

The CAR and CCR  relations are usually defined using Fock-space
representations. We will prove a noncommutative version of the
Hamburger moment problem which allows us to show that the
ultraproduct construction and the Fock space models describe the
same von Neumann algebra.

We will say a {\it $^*$-probability space} $(\A,*,\phi)$ is given
by unital $^*$-algebra $\A$ over the complex numbers and a
positive linear functional $\phi:\A\to \cz$ with $\phi(1)=1$. Here
$^*$-algebra means that $^*$ is an antilinear involution and that
 \[ \A \lel  \A_{sa}+ i \A_{sa} \pl \quad
 \mbox{where} \quad
  \A_{sa} \lel \{a\in A \p|\p a^*=a\} \]
is called the  {\it selfadjoint part} of $\A$. The functional
$\phi$ is called positive, if $\phi(a^*a)\ge  0$ and
$\phi(a^*)=\overline{\phi(a)}$.  A \emph{representation} of $\A$
is given by a generating set  $S\subset \A_{sa}$,
 \begin{enumerate}
  \item[i)] a Hilbert space $H$, a unit vector $\xi$, and a
  linear map
  \begin{equation} \label{hilb}
   \al:\A\to H \quad \mbox{with}\quad  \al(1)\lel \xi \pl ;
   \end{equation}
  \item[ii)] a  map $\pi:S\cup\{1\}\to Sa(H)$, $Sa(H)$ the set of selfadjoint densely defined operators on $H$,
such that  $\pi(1)=1$, and $\al(\A)$ is a subset of the domain of
$\pi(a)^k$ for all $a\in S$, $k\in \nz$, and
   \begin{equation} \label{specpr}
    \pi(a)\xi \lel \al(a) \quad ,\quad
    \intt t^k (\al(b),dE_t^{\pi(a)}\al(b)) \lel
    (\al(b)
    ,\al(a^kb)) \end{equation}
  holds  for all $a\in S$ with spectral resolution
  $E_t^{\pi(a)}$ and for all  $b\in \A$.
   \item[iii)] Moreover, we require
    \begin{equation} \label{scprp} \phi(a^*b)\lel
   (\al(a),\al(b))_H \pl
   \end{equation}
for all  $a,b\in   \A$.
\end{enumerate}

Let us say that $\A$ satisfies the \emph{growth condition} if
there exists a generating subset $S\subset \A_{sa}$ such that
$S\cup \{1\}$ generates $\A$ as an algebra and there exists a
length function $|\pll|:S\to (0,\8)$ such that
 \[ |\phi(a_1\cdots a_n)| \kl n^n \prod_{i=1}^n |a_i| \pl
 \]
holds for all $a_1,..,a_n\in S$.  We will show that under these
assumptions the von Neumann algebra generated  by the spectral
projections of elements $\pi(a)$, $a\in S$ is uniquely
determined. We will need the following formulation of the
Hamburger moment problem (see e.g. \cite{Wid}).

\begin{theorem}[Hamburger moment problem]\label{hammom} Let $\mu_1$ and $\mu_2$ be regular Borel measures
on $\rz$ such that the moments
 \[ \intt t^kd\mu_1(t)\lel m_k\lel \intt t^kd\mu_2(t) \]
coincide.  If there exists a constant $c$ such that $m_k\le
c^{k+1}k^k$, then $\mu_1\lel \mu_2$. Moreover, under these
assumptions on $m_k$ the polynomials are dense in $L_p(\mu_1)$ for
every $p<\infty$.
\end{theorem}

\begin{proof}[Sketch of Proof.] It is well-known that the Hamburger moment
problem has a positive solution under these conditions. We will
only indicate the proof of the  density assertion needed in this
paper. We consider the function $h(x)\lel \frac{\sin(x)}{x}$ and
the translates $h_{\la,s}(x) \lel h(s(x-\la))$. One first shows
that for $2\le p<\infty$ and $s(p)=(6cep)^{-1}$ the functions
$(h_{\la,s})_{s<s(p)}$ belong to the closure of the polynomials in
$L_p(\rz,\mu_1)$. Indeed, let $q_m(x)= \sum_{k=0}^m (-1)^k s^{2k}
\frac{(x-\la)^{2k}}{(2k+1)!}$ and $p\in 2\nz$. Assume $m>|\la|+1$.
Then
 \begin{align*}
 \noo  h_{\la,s}-q_m\rrm_p
  &\le \summ_{k>m} \frac{s^{2k}}{(2k+1)!}
  \noo (x-\la)^{2k}\rrm_p \kl
   \summ_{k>m} \frac{s^{2k}}{(2k+1)!}
   c^{2k+\frac{1}{p}}
   (2kp+|\la|)^{2k} \\
  &\le   c^{\frac1p}   \pl
   \summ_{k>m} \frac{(cs)^{2k}}{(2k+1)!}  (2kp+k)^{2k}\kl  c^{\frac1p}  \pl
   \summ_{k>m} (3cpes)^{2k} \frac{1}{2^{k}(2k+1)} \pl .
   \end{align*}
Here, we used the the  estimate $\noo (x-\la)^{2k}\rrm_p\le
c^{2k+\frac{1}{p}}(2kp+|\la|)^{2k}$ for $k>|\la|+1$. Let
$s<s(2p)$. Now we show by induction on $m$ that for every
polynomial $q$ the functions $h_{\la_1,s}\cdots h_{\la_m,s}q$
belong to the closure of the polynomials in the $L_p(\rz,\mu_1)$.
In particular, the algebra $A$ generated by
$(h_{\la,s})_{s<s(2p),\la \in \rz}$ is in the closure of the
polynomials in $L_p(\mu_1,\rz)$. Finally, we deduce from the
Stone-Weierstrass theorem that $A$ is dense in $C(\rz\cup
\{\infty\})$ and hence in $L_p(\mu_1)$ because $\mu_1$ is finite.
We deduce that $1_{(a,b)}$ is  in the $L_p$-closure of polynomials
for all $-\infty\le a<b\le \infty$ and $p<\infty$. By regularity
of the Borel measure $\mu_1$ the step functions are dense in $L_p$
and this completes the proof.\qd

Let us now verify that the growth condition allows us to apply the
Hamburger moment problem. We consider a representation
$(\pi,\al,\xi)$ of $\A$, and $a\in S$ and $b\in \A$. Then  we may
define the regular  Borel measure  (see e.g. \cite[I Theorem
5.2.6]{Kad})
 $\mu$ on the real line given by
 \[ \mu(B) \lel \int_B (\al(b),dE_t^{\pi(a)} \al(b))  \pl .\]
Clearly, $\mu$ is a finite measure satisfying
 \[ \mu(\rz)\lel \phi(b^*b) \pl .\]
The following estimate follows similarly as in \eqref{momentum}
and we leave it to the reader.

\begin{lemma} \label{m1} Let $b\in \A$ and $a\in S$.
Then there exists a constant $C$  such that for all $k\in \nz_0$,
$\la \in \rz$
\[ |\int_\rz (x+\la)^k d\mu(x)|  \kl  C^{k+1}
 (k+|\la|)^k\pl .\]
\end{lemma}
\hz

\begin{theorem} \label{uniq} Let $(\A,\phi)$ be a $^*$-probability space
satisfying the growth condition with respect to $S$. Let
$(\pi_1,\al_1,\xi_1)$ and  $(\pi_2,\al_2,\xi_2)$ be
representations. Let $\N_1$, $\N_2$ be the von Neumann algebras
generated by the spectral projections of $\{\pi_1(x)\}_{x\in S}$,
$\{\pi_2(x)\}_{x\in S}$, respectively. Assume in addition that
the restrictions $\phi_1(x)=(\xi_{1},x\xi_{1})$,
$\phi_2=(\xi_{2},x\xi_{2})$, are faithful on $\N_1$, $\N_2$,
respectively. Then there is a normal homomorphism $\pi:\N_1\to
\N_2$ such that $\phi_2\circ \pi=\phi_1$.
\end{theorem}

\begin{proof} Let $\A_{\phi}$ be the pre-Hilbert space $\A$
equipped with the scalar product $(b,c)=\phi(b^*c)$ and $H_{\phi}$
its completion. Let us first consider a single representation
$\al_{1}:\A_\phi\to H_1$. Then $\al_1$ extends to an isometric
isomorphism between $H_{\phi}$ and the closure of $\al_1(\A)$.
Given $b\in \A$, $a\in S$, we may consider the measure
 \[   \mu_1(B) \lel  \int_B
 \big( \al_1(b),dE_t^{\pi_1(a)}(\al_1(b))\big) \pl .\]
The measure  is a  regular Borel measure (see e.g. \cite[I
Theorem 5.2.6]{Kad}) and satisfies the growth condition.
Moreover, given a bounded function $f$, there exists a polynomial
$q_n$ such that $\noo f-q_n\rrm_{L_2(\mu_1)} \kl \frac{1}{n}$.
Therefore, we deduce
 \for
  \noo f(\pi_1(a))\al_1(b)-
  q_n(\pi_1(a))\al_1(b)
  \rrm_2 &=&
  \noo f-q_n\rrm_{L_2(\mu_1)} \kl n^{-1} \pl .
  \mel
However,
 \[ q_n(\pi_1(a))(\al_1(b))  \lel \al_1(q_n(a)b)  \]
and thus  $q_n(a)(b)$ forms a Cauchy sequence in $H_{\phi}$
converging to some element $h\in H_{\phi}$. Thus we get
 \for
  f(\pi_1(a))(\al_1(b)) &=&
  \lim_n q_n(\pi_1(a))(\al_1(b)) \lel
  \lim_n \al_1(q_n(a)b)
  \lel \al_1(h)  .
 \mel
We observe that $K_1=cl(\al_1(\A))$ is an invariant subspace for
$\N_1$ and hence there exists a projection $e_1\in \N_1'$ such
that $K_1=e_1H_1$. Then $\al_1$ provides a unitary between
$H_{\phi}$ and $K_1$. We apply the same argument to $\al_2$ and
obtain $e_2\in \N_2'$, $K_2=e_2H_2$ and an isometry
$\al_2:H_{\phi}\to K_2$. We define the unitary
$u=\al_2\al_1^{-1}:e_1H_1 \to e_2H_2$  and want to show
 \begin{equation}\label{rel}  u^{-1}f(\pi_2(a))u \lel
 f(\pi_1(a)) \end{equation}
holds for all bounded measurable functions $f$ and $a\in S$. Let
$a\in S$ and $b\in \A$. The argument above shows that
 \[  \al_1(h) \lel  f(\pi_1(a))\al_1(b)  \lel \lim_{n} \al_1(q_n(a)b)  \pll
 \mbox{and}\pll
 \al_2(h) \lel  f(\pi_2(a)) \al_2(b) \lel \lim_n \al_2(q_n(a)b) \pl .\]
Therefore, it suffices to observe that the sequence $(q_n(a))$
from above can be chosen to work for $\pi_1(a)$ and $\pi_2(a)$
simultaneously. This follows immediately from the fact that the
measure $\mu_2$ given by
 \[   \mu_2(B) \lel  \int\nolimits_B
 \big( \al_2(b),dE_t^{\pi_2(a)}(\al_2(b))\big) \pl \]
has the same moments as $\mu_1$ and thus $\mu_1=\mu_2$ by Theorem
\ref{hammom}. This completes the proof of  \eqref{rel}. We can
now define the  state preserving normal $^*$-homomorphism
$\tilde{\pi}: e_1\N_1\to e_2\N_2$ given by
$\tilde{\pi}(x)=u^{-1}xu$. The conclusion follows from the fact
$\N_1$ and $e_1\N_1$ are isomorphic and $\N_2$ and $e_2\N_2$ are
isomorphic. Indeed, the induction $\rho_1:\N_1\to e_1\N_1$ defined
by $\rho_1(x)=e_1x$ is a normal $^*$-homomorphism. For a positive
element $\rho_1(x)=0$ implies
 \[ \phi_1(x)\lel (\xi_{1},\rho_1(x)\xi_{1}) \lel
 0 \pl .\]
Since $\phi_1$ is faithful, $\rho_1$ is injective and therefore a
normal isomorphism (see \cite[I.4.3 Corollary 1]{Dix}). The same
applies for $\rho_2$. Hence $\pi=\rho_2^{-1}\tilde{\pi}\rho_1$
yields the state preserving homomorphism. \qd

\begin{rem} {\rm If we assume only that $\phi_1$ is faithful we still
have an isomorphism $\pi:\N_1\to e_2\N_2$.}
\end{rem}

In the following we will often have some control on the modular
group:

\begin{prop} \label{mod} Let $(\A,\phi)$ be a $^*$-probability
space with generating system $S$ satisfying the growth condition.
Let $(\si_t)_{t\in \rz}$ be a family of maps such that
$\si_t(S)\subset S$. Let $(\pi,\al,\xi)$ be a representation of
$(\A,\phi)$ and $\M$ be a von Neumann algebra such that $\pi(x)$
is affiliated to $\M$ for all $x\in S$ and
$\phi_{\xi}(y)=(\xi,y\xi)$ is faithful on $\M$. Assume that
 \[ \pi(\si_t(x)) \lel \si_t^{\phi_{\xi}}(\pi(x)) \pl \]
holds for all $t\in \rz$ and $x\in S$. Let $\N$ be the subalgebra
of $\M$ generated by the spectral projections of $\pi(S)$. Then
there is a conditional expectation $\E:\M\to \N$ onto $\N$ such
that $\phi\circ \E\lel \phi$.
\end{prop}

\begin{proof} Let $x\in S$ and $B\subset \rz$ be measurable.
Since $\si_t^{\phi_{\xi}}$ is an automorphism on $\M$, we deduce
 \[ \si_t^{\phi_{\xi}}(1_B(\pi(x))) \lel 1_B(\si_t^{\phi_{\xi}}(\pi(x)))
 \lel 1_B(\pi(\si_t(x)))
  \pl .\]
Therefore $\si_t^{\phi_{\xi}}$ leaves $\N$ invariant.  The result
follows by an application of Takesaki's theorem (see \cite[Theorem
10.1]{Strat}).\qd

As an application we will show that the von Neumann algebras
$\Gamma_q(K,U_t)$ can be obtained from the central limit
procedure. This is an important link for our norm estimates.

\begin{cor}\label{identify} Let $(K,H,U_t)$ be as in section 2. Then
$\Gamma_q(K,U_t)$ is isomorphic to a complemented subalgebra $N
\subset M$ obtained in Theorem \ref{limitth}. More precisely,
there exists a homomorphism $\pi: \Gamma_q(K,U_t)\to M$ such that
$\phi\circ \pi$ is the vacuum vector state $\phi_{vac}$,
$\si_t^{\phi}\circ \pi=\pi \circ \si_t^{\phi_{vac}}$ and
 \begin{equation} \label{p-rel}  \pi(s_q(h)D_{vac}^{\frac1p})\lel u_q(h)D^{\frac1p}\
 \end{equation}
holds for all $h\in K$. Here $u_q$ is the map $u$ constructed for
$\beta(\si)=q^{I(\si)}$.
\end{cor}

\begin{proof} According to Lemma \ref{arak-form}, we may assume
$K=K_0^{\rz}\oplus L_2(\Om_1,\mu;\rz)$. Using a conditional
expectation at the end of our proof, we may assume that the
dimension of $K_0$ is even  and hence given by
$K_0=L_2(\Om_0,\mu_0;\cz)$. Then we use the disjoint union
$\Om=\Om_0\cup \Om_1$ with the measure $\mu=\mu_0+\mu_1$. We
extend the functions $f_1$ and $f_2$ defined on $\Om_1$ via Lemma
\ref{arak-form} by $f_1(\om)=f_2(\om)=\frac12$ on $\Om_0$. We
work with our standard model $N=L_{\infty}(\Om,\mu;\Mz_2)$ from
section 3.  We apply the central limit procedure Theorem
\ref{limitth} for the partition function $\beta(\si)=q^{I(\si)}$
using  Speicher's random model (see Corollary \ref{prob}). Our
generating set is given by
 \[ S \lel \{ j(f)\pl:\pl  f\in L_{\infty}(\Om,\mu;\cz)\}
 \pl .\]
Let $\N(S)$ be the von Neumann algebra generated by the spectral
projections of elements $u_q(s)$, $s\in S$. Since $U_t$ leaves $S$
invariant, we can apply Proposition \ref{mod} and find a normal
faithful state preserving conditional expectation $E:M \to \N(S)$.
Using the combinatorial formulae \eqref{q-com} and \eqref{q-rel},
we see  $u_q(S)$ and $s_q(S)$ have the same moment formulae and
satisfy the growth condition. By Corollary \ref{growth} and
Theorem \ref{uniq}, we conclude that $\N(S)$ and the subalgebra
$\Gamma_q(S)\subset \Gamma_q(K,U_t)$ generated by $s_q(S)$ in the
Fock space construction are isomorphic via an isomorphism
respecting the state and the modular group. Moreover
$L_{\infty}(\Om,\mu)\cap L_2(\Om,\mu)$ is dense in $L_2(\Om,\mu)$
and the map $u_{q,2}:L_2(\Om,\mu)\to L_2(\Gamma_q(K,U_t))$,
$u_{q,2}(f)=s_q(f)D_{vac}^{1/2}$ is an isometry. Here $D_{vac}$ is
the density of the state given by the vacuum. For $-1\le q<1$ we
know that $s_q(f)$ is bounded and hence $s_q(f)=SOT-\lim s_q(f_k)$
also belongs to the von Neumann algebra $\Gamma_q(S)$ (see
\cite{S} and \cite{Hiai}). Thus $\Gamma_q(S)=\Gamma_q(K,U_t)$ and
\eqref{p-rel} holds for arbitrary elements by approximation (see
e.g. \cite[Lemma 2.3]{JD}).

For $q=1$ the argument is slightly different because we can no
longer define $\Gamma_{1}(K,U_t)$ as generated by the unbounded
operators $\{s_1(h):h\in K\}$. In that case we recall Segal's
formulation of the commutation relations given by a complex
Hilbert space $H_U$ and a family of unitaries $W(h)$:
 \begin{equation}\label{weyl}
  W(h_1+h_2)\lel e^{-\frac{i}{2}{\rm Im}(h_1,h_2)}\quad \mbox{and} \quad
 (\Om,W(h)\Om)\lel e^{-\frac{1}{4}\|h\|^2} \pl .
 \end{equation}
Then we have
 \[
  \Gamma_1(K,U_t) \lel R_{Segal}(K/H_U) \lel \{W(h):h\in K\}''  \pl .\]
In Example \ref{22-ext} we have calculated the form $B(f,g)={\rm
Im}\psi(j(f)j(g))=(A(f),g)$ which determines $\Gamma_1(K,U_t)$.
Here we have $A(f)=(f_2-f_1)f$. For $f\in L_{\infty}(\mu;\cz)\cap
L_2(\mu;\cz)S$ we know by Theorem \ref{limitth} that the unitary
$\tilde{W}(f)= e^{iu_1(j(2^{-1/2}f))}\in \N(S)$. Using Remark
\ref{contccr} and the fact that $\al(A(S))$ is dense in
$L_2(\N(S))$ (see also the proof of Theorem \ref{uniq}), we find
the same commutation relations as in \eqref{weyl}. Since the
commutation relations uniquely determine $R_{Segal}(K/H_U)$, we
see that the $\Gamma_1(K,U_t)$ and the algebra
$\tilde{\N}(S)\subset \N(S)$ generated by the family
$\{\tilde{W}(f):f\in S\}\subset \N(S)$ are isomorphic with respect
to a state preserving homomorphism. (Alternatively, we may apply
Theorem \ref{uniq} to the real and imaginary parts of these
unitaries). Note
 that $\tilde{\N}(S)$ is invariant under the modular group
and thus Takesaki's theorem (see \cite[Theorem 10.1]{Strat})
yields a normal state preserving conditional expectation onto
$\tilde{\N}(S)$ (in fact we have $\tilde{\N}(S)=\N(S)$). Formula
\eqref{p-rel} is obtained by differentiation for elements $s\in
S$. In order to show that \eqref{p-rel} holds for arbitrary
elements in $K$,  we recall that the distribution of $u_1(j(f))$
with respect to $\phi$ is the normal distribution
$N(0,\|f\|_2^2)$. Thus we have $\|u_1(j(f))\|_p\kl
c\sqrt{p+1}\|f\|_2$. According to Lemma \ref{assoc}, we deduce
$\|u_1(j(f))D^{\frac1p}\| \kl c \sqrt{p+1}\|f\|_2$. By density we
obtain \eqref{p-rel}. \qd

\begin{rem}\label{typp} {\rm In order to determine the type of
$\Gamma_1(K,U_t)$, we may transform Segal's representation into
Weyl's representation $W(h)=e^{i/2({\rm Re}(h),{\rm Im}(h)}U({\rm
Re}(h))V({\rm Im}(h))$. According to \cite{Ar} we get
 \begin{align}
  \Gamma_1(K,U_t) &=  \{W(h):h\in K\}'' \nonumber \lel
  \{U({\rm Re}(h))V({\rm Im}(h)): h\in K \}'' \lel
 R(K_1,\rho(K_1)/K)  \label{weyl-segal} \pl .
 \end{align}
Here $K_1,K_2$ and $\rho$ are chosen according to the definitions
before Corollary \ref{typeIII1} such that
$(f,\rho(g))=(f,(f_2-f_1)g)$. Thus the type depends only on the
spectrum of $1+|f_2-f_1|^{-1}$.}
\end{rem}

\begin{rem}\label{qwep} {\rm
1) As an application we obtain a generalization of a very recent
result of Nou \cite{Nou}: The von Neumann algebras
$\Gamma_{q}(K,U_t)$ are QWEP or equivalently allow `matrix
models'.  Although it is open whether every von Neumann algebra is
QWEP, it is sometimes difficult to verify this property
explicitly. Let us recall that a $C^*$-algebra has WEP if
$A\subset A^{**}\subset B(H)$ and there is a contraction
$P:B(H)\to A^{**}$ such that $P|_A=id|_A$. A $C^*$-algebra  $B$ is
QWEP if there exists a $C^*$-algebra A with WEP and a two sided
ideal $I\subset A$ such that $B=A/I$. Since the QWEP property is
stable under ultraproducts and conditional expectations (see e.g.
\cite{Joh}), we deduce that the algebra $\N(S)$ is QWEP because in
this case the von Neumann algebra $M$ in Theorem \ref{limitth} is
constructed from an ultraproduct of type I von Neumann algebras.
The same argument applies for the algebra $M[-1,1]$ from Corollary
\ref{all}.

2) More generally we may define a von Neumann algebra
$\Gamma_{\beta}(K,U_t)$ given by a partition function $\beta$
constructed as in \ref{gen-part} and an inclusion $K\subset H_U$
as above. We define a positive functional on $A(K)$ by
 \begin{equation}\label{beta}  \phi_{\beta}(x_1\ten \cdots \ten x_n) \lel
 \summ_{\si=\{\{i_1,j_1\},...,\{i_{\frac{m}{2}},j_{\frac{m}{2}}\}\}
 \in P_2(m)}
 \beta(\si) \prod_{l=1}^{\frac{m}{2}} (x_{i_l},x_{j_l})_U \pl .
 \end{equation}
This induces a Hilbert space $(A(K),(\pl ,\pl)_{\beta})$ obtained
by completion of $(a,b)_{\beta}=\phi_{\beta}(a^*b)$. According to
Lemma \ref{arak-form} it suffices to consider our standard model
$N=L_{\infty}(\Om,\mu;\Mz_2)$ with functions $f_1+f_2=1$. The
growth condition is satisfied for
$|f|=\max\{\|f\|_2,\|f\|_{\infty}\}$ and $S$ the set of
selfadjoint elements of the form $x=e_{12}\ten f+e_{21}\ten
\bar{f}$, $f\in L_{\infty}$.  Hence the ultraproduct construction
Theorem \ref{limitth} provides unbounded operators $s_{\beta}(x)$
affiliated to $N_\U$. Using Proposition \ref{mod}, we see that the
von Neumann algebra $\N_{\beta}(S)$ generated by the spectral
projections admits a state preserving conditional expectation. In
order to extend this representation to all elements in $K$, we
introduce the new length function $|x|_U=\|x\|_U$ and note that
 \[ \lim_{j} \phi_{\beta}(a_j\ten e^{zx_j}\ten b_j) \lel
 \phi_{\beta}(a\ten e^{zx}\ten b) \]
holds whenever we have convergence $\lim_j x_j=x$, $\lim_j a_j=a$,
$\lim_j b_j=b$ in $H_U$. By density we find an isomorphism $u$
between the completion of $A(K)$ with respect to the scalar
product induced by $\phi_{\beta}$ and $L_2(\N(S),\phi_\U)$. For a
selfadjoint element $x=e_{12}\ten f+e_{21}\ten \bar{f}$ we choose
an approximating sequence with $f_j\in L_{\infty}$ and define
$U_t^x=w^*-\lim_j e^{its_{\beta}(x_j)}$ using an ultra-limit in
the weak operator topology. For $b\in A(K)$, we may also define
the measure $\mu=w^*-\lim_j \mu^j_b$ of the measures
$\mu^j_b(A)=\int_A (u(b),dE_t^{\pi(x_j)}u(b))$. Using the length
function $|x_j|_U$ we see that the growth condition is uniformly
satisfied for all $j$, hence also for $\mu$. From this we deduce
strong convergence $U_t^x=\lim_j e^{itx_j}$ and multiplicativity
$U_s^xU_t^x=U_{s+t}^x$. We find a von Neumann algebra
$\Gamma_{\beta}(K,U_t)$ generated by $(e^{its_{\beta}(x)})_{x\in
K}$, complemented in $\N_\U$. Moreover, by construction
$\Gamma_{\beta}(K,U_t)$ is QWEP and \eqref{beta} holds for the
generators. It would be interesting to know what kind of functions
$\beta$ defined on the pair partitions have similar properties.}
\end{rem}

\section{An inequality for sums of independent copies } \setcounter{equation}{0}

In this section we will prove the main inequality of the paper
based on an appropriate concept of independent copies. Let us fix
von Neumann algebras $\Mz\subset \N$, and a normal, faithful,
conditional expectation $E:\N\to \Mz$. Let $\Mz\subset A,B\subset
\N$ be subalgebras. $A$ and $B$ are said  to be {\it independent
over $E$} if
 \begin{equation}\label{indepo}  E(ab)\lel E(a)E(b) \end{equation}
holds for all $a\in A$ and $b\in B$. We are indebted to C.
K\"ostler (see \cite{ko2}) for pointing out the condition ii) in
the following Lemma:

\begin{lemma}\label{pyind} Let $A$ and $B$ independent over $E$.
Let $\phi_{\Mz}$ be  a normal faithful state on $\Mz$ such that
$\si_t^{\phi_{\Mz}\circ E}(A)\subset A$.
  \begin{enumerate}
  \item[i)] There   exists a conditional expectation $E_A:\N\to A$ such that
   \[ E_A(b)\lel E(b) \]
 for all $b\in B$.
  \item[ii)] Let $a_1,a_2\in A$ and $b\in B$. Then
   \[ E(a_1ba_2)\lel E(a_1E(b)a_2) \pl .\]
  \item[iii)] Let $a,b\in \Mz_m(A)$. Then
   \[ \noo (id\ten E)(a^*b^*ba)\rrm \kl \noo (id\ten  E)(a^*a)\rrm \noo
 (id \ten    E)(b^*b)\rrm \pl .\]
\end{enumerate}
\end{lemma}
\begin{proof} If $\phi_{\Mz}$ is a normal faithful state on $M$, then
$\phi=\phi_{\Mz}\circ E$ is a normal faithful  state on $\N$ such
that $\si_t^{\phi_{\Mz}}\circ E\lel E\circ \si_t^{\phi}$ (see
\cite{Co}). Then Takesaki's theorem implies that $E$ is the unique
conditional expectation such that
 \[ \phi(ab)\lel \phi(aE(b)) \]
holds for all $a\in M$ and $b\in \N$. By assumption on the modular
group, we may apply Takesaki's theorem (see \cite[Theorem
10.1]{Strat}) and find a normal conditional expectation $E_A:\N\to
A$ which is characterized by
 \[ \phi(ab) \lel \phi(aE_A(b)) \]
for all $a\in A$, $b\in \N$. However, by independence and the
module property of $E$  we deduce
 \begin{align*}
 \phi(aE_A(b))&=  \phi(ab) \lel \phi_{\Mz}(E(ab)) \lel
 \phi_{\Mz}(E(a)E(b))\lel  \phi_{\Mz}(E(aE(b))) \lel \phi(aE(b))
  \end{align*}
for   all $a\in A$ and $b\in B$. Now, we prove ii). Let us
consider $a_1\in A$, a $\phi$-analytic element $a_2\in A$ and
$m\in M$. Then, we deduce from i)
 \begin{align*}
 &\phi\big(mE(a_1ba_2)\big) \lel  \phi(ma_1ba_2) \lel
 \phi(\si_{i}^{\phi}(a_2)ma_1b) \lel
  \phi(\si_{i}^{\phi}(a_2)ma_1E_A(b)) \\
  &=
 \phi(\si_{i}^{\phi}(a_2)ma_1E(b)) \lel
  \phi(ma_1E(b)a_2)\lel  \phi\big(mE(a_1E(b)a_2)\big) \pl
 .
 \end{align*}
Since this is true for all $m\in \Ma$, we deduce ii) for analytic
elements $a_2$. By approximation with bounded nets of analytic
elements in the strong topology, the assertion follows for general
$a_2$. For the proof of iii) we first note that $\Mz_m(A)$ and
$\Mz_m(B)$ are independent over $id_{\Mz_m}\ten E$. Indeed, we
deduce from linearity that
 \begin{align*}
 (id_{\Mz_m}\ten E)(ab) &=  [E(\summ_k a_{ik}b_{kj})]_{ij} \lel
 [\summ_k E(a_{ik})E(b_{kj})]_{ij} \lel
  (id_{\Mz_m}\ten E)(a) (id_{\Mz_m}\ten
 E)(b) \pl .
 \end{align*}
If $\phi$ is a normal faithful state, then $\phi_m=
\frac{tr}{m}\ten \phi$ is a normal faithful state on $\Ma_m(\N)$.
It is easily checked that the assumption on the modular group on
$A$ implies the same assumption on the modular group for
$\Mz_m(A)$. Hence condition ii) holds on the  matrix level. By
positivity we deduce
 \begin{align*}
  (id_{\Mz_m}\ten E)(a^*b^*ba) &\lel  (id_{\Mz_m}\ten E)(a^*(id_{\Mz_m}\ten
  E)(b^*b)a)\\
 &\le
   \noo (id_{\Mz_m}\ten E)(b^*b)\rrm \pl (id_{\Mz_m}\ten E)(a^*a) \pl .
  \end{align*}
Taking norms implies the assertion.\qd

We will now discuss  the notion of subsymmetric copies needed for
our key inequality. This notion is closely related to the order
invariance discussed in \cite{KS}. We refer to \cite{Ku} for more
information on the notion of white noise in the continuous
setting which seems to be closely related. We will consider an
inclusion of three von Neumann algebras
 \[ \Ma\subset \M \subset \N \pl. \]
We will always assume that we have  normal faithful conditional
expectations $E:\N\to \Ma$ and $\E:\N\to \M$ such that
$E|_{\M}\circ \E=E$. A system of \emph{subsymmetric copies}
$(\Ma,\M,\N,\al_1,...,\al_n,E)$ over $\Ma$ (strictly speaking over
$E$) is given by faithful normal isomorphisms $\al_i:\M\to
\M_i\subset \N$ such that $\al_i\circ E\lel E\circ \al_i$ holds
for all $i=1,...,n$ and
 \begin{equation}\label{subsym} E(\al_{i_1}(a_1)\cdots \al_{i_m}(a_m))
 \lel E(\al_{j_1}(a_1)\cdots \al_{j_m}(a_m))
 \end{equation}
holds for all $a_1,...,a_m\in \M$ and order-equivalent functions
$\i,\j :\{1,...,m\}\to \{1,...,n\}$. We recall that two functions
$\i,\j$ are order equivalent if
 \[ i_k\le i_l \quad \Longleftrightarrow \quad  j_k\le j_l \]
holds for all $1\le k,l\le m$. Subsymmetric, not necessarily
symmetric, copies appear naturally in the context of iterated
crossed products (see \cite{JRo}). For our proof we need a
slightly weaker assumption. We say that
$(\Ma,\al_1(\M),...,\al_n(\M),\N,E)$ are top-subsymmetric copies
if \eqref{subsym} holds for $\i$ and $\j$ of the special form
 \[ \i|_{\{1,...,m\}\setminus A}=\j|_{\{1,...,m\}\setminus A} \]
$|A|\le 2$ such that $k\in A$ implies
 \[  i_k=\max\{i_l:1\le l\le m\} \quad \mbox{and} \quad  j_k=\max\{j_l:1\le l\le m\} \pl .\]
Intuitively speaking we are only exchanging at most two top
values. $(\Ma,\M,\N,\al_1,...,\al_n,E)$ is called a system of
\emph{independent subsymmetric copies} if in addition the von
Neumann algebra $\M_i$ generated by $\bigcup_{j<i} \al_i(\M)$ is
independent of $\al_i(\M)$ over $\Ma$ for  all $i=1,...,n$ (see
e.g. \cite{Sch} for more information).  We will use the term
\emph{top-subsymmetric independent copies} if the system is
top-subsymmetric and the dependence condition is satisfied.

In our applications, we will often find a stronger assumption of
a system of \emph{symmetric independent copies}
$(\Ma,\al_1(\M),...,\al_n(\M),\N,E)$. This means that the
independence condition is satisfied and
 \begin{equation}\label{sym} E(\al_{i_1}(a_1)\cdots \al_{i_m}(a_m))
 \lel E(\al_{\si(i_1)}(a_1)\cdots \al_{\si(j_m)}(a_m))
 \end{equation}
holds for every function $\i:\{1,...,m\}\to \{1,...,n\}$ and every
permutation $\si$ of $\{1,...,n\}$. It is clear that this
condition implies top-subsymmetry by considering an inversion
$(i_kj_k)$ of the largest elements. We observe that \eqref{sym}
induces an automorphisms $\al_{(k,l)}$ of the von Neumann algebra
$\M_n$ generated by $\al_1(\M),...,\al_n(\M)$  which interchanges
elements of $\al_k(M)$ and $\al_l(\M)$ while leaving all other
elements fixed. Here are some examples.

\begin{exam} {\rm Let $\N=\Mz\ten N^{\ten_n}$ and $\M=\Mz\ten
\pi_n(N)$. We define an automorphism $\al_i$ by
 \[ \al_i(m\ten x_1\ten \cdots \ten x_{i-1}\ten x_i\ten
 x_{i+1}\ten \cdots \ten x_n)\lel
 m\ten x_1\ten \cdots \ten x_{i-1}\ten x_n\ten
 x_{i+1}\ten \cdots \ten x_i \pl.\]
If $\phi_{\Mz}$, $\phi_N$ are  faithful normal states then
$\phi=\phi_M\ten \phi_N^{\ten_n}$ is a faithful normal state on
$\N$. The conditional expectation $E$ is given by $id\ten
\phi_N^{\ten_n}$. Then it is easily checked that $(\Mz\ten
1,\M,\N,(\al_i)_{i=1,..,n},\phi)$ is a system of independent
copies. The same argument applies in the free product situation
$\N=\Mz\ten \ast_{i=1,..,n}(N,\phi)$ where $\M=\Mz\ten\pi_n(N)$ is
given by the n-th copy. In the next example we extend this to
$q$-independent copies.}
\end{exam}

\begin{exam} {\rm Let $N=\ell_\infty^n(\ell_{\infty}(\Mz_2))$ and
 \[ \psi(x) \lel \summ_{i=1}^n \summ_{j\in \nz}
 [(1-\mu_j)x_{11}(i,j)+ \mu_jx_{22}(i,j)] \pl .\]
We assume that $(v_k(n))$ satisfies the assumptions of Theorem
\ref{limitth} and denote by $\beta$ the resulting partition
function. We consider the algebra $\N_{\beta} \subset N_\U$
generated by the spectral projections of the selfadjoint operators
$\{u(x)\}_{x\in N_{sa,a}}$ from Proposition \ref{limit1}. Recall
that by Corollary \ref{growth} the growth condition is satisfied
and hence $N_{\beta}$ is uniquely determined by the partition
function $\beta$. Now, we consider the $\psi$-preserving
automorphisms
\[ \gamma_i(x_1,..,x_{i-1},x_i,x_{i+1},..,x_n)\lel (x_1,..,x_{i-1},x_n,x_{i+1},..,x_i) \pl .\]
According to Theorem \ref{uniq} we find $\phi_\U$-preserving
automorphisms $\al_i:\N\to \N$ such that
\[ \al_i(u(x))\lel u(\gamma_i(x)) \pl .\]
Let $f_n$ be the projection onto the last coordinate in
$\ell_{\infty}^n$. Let $\M\subset \N_{\beta}$ be the subalgebra
generated by $\{u(x)\}_{f_nxf_n\in (f_nNf_n)_{sa,a}}$. Using
$\si_t^{\phi_\U}(u(x))=u(\si_t^{\psi}(x))$ we may apply
Proposition \ref{mod} and find a conditional expectation
$E:\N_{\beta}\to \M$. Then $(\cz,\N_{\beta},
\al_1(\M),...,\al_n(\M),\phi_{\beta})$ is a system of independent
copies over $\phi_\U$. Tensoring with matrices $\Mz_m$, we find
examples of independent copies over $id_{\Mz_m}\ten \phi_\U$.}
\end{exam}\lz

For our proof we need some facts from modular theory. Unless
stated otherwise we will assume that $\Ma$ is $\si$-finite. Let
$\phi_{\Ma}$ be a normal faithful state. Then we find normal
faithful states  $\phi_\M=\phi_{\Ma}\circ E|_{M}$ and
$\phi=\phi_{\Ma}\circ E$ on $\M$ and $\N$ respectively. Due to the
work of Connes \cite{Co} it is well-known that then
 \begin{equation} \label{inv00}
  \si_t^{\phi_\M}(\Ma)\subset \Ma \quad , \quad
  \si_t^{\phi}(\Ma)\subset \Ma
  \end{equation}
holds for the modular group of the corresponding states. Moreover,
from $E=E\E$, we deduce  $\phi\lel \phi_{\M}\circ \E$ and hence we
also have
 \begin{equation}\label{inv0} \si_t^{\phi}(\M)\subset \M \pl .
 \end{equation}

\begin{lemma}\label{condit} Let $(\Ma,\al_1(\M),...,\al_n(\M),\N,E)$ be a system
of symmetric independent copies and $\N=\M_n$ the von Neumann
algebra generated by  $\al_1(\M),...,\al_n(\M)$. Then there exist
conditional expectations $\E_i:\N\to \M_i$ such that
 \begin{equation}\label{condi} \E_{i}(x)\lel E(x)
 \end{equation}
holds for all $x\in \al_{i+1}(\M)$. Moreover, we have
 \begin{equation}\label{modul}
  \al_i\circ \si_t^{\phi_\M}\lel  \si_t^{\phi} \circ
 \al_i \pl. \end{equation}
and
 \begin{equation}\label{compat} \phi\circ \al_i\lel \phi \quad  \mbox{and} \quad \E_i\lel \al_i \circ \E \circ
 \al_{1,i}^{-1}
 \end{equation}
for all $i=1,...,n$.
\end{lemma}

\begin{proof} We observe that $\al_{1,i}:\M_n\to \M_n$ is an
automorphism of $\M_n$ such that $\al_{1,i}|_M=\al_i$. Moreover,
our assumption implies  that we have $E\circ \al_{1,i}\lel E$.
This allows us to define $E_i\lel \al_i \E\al_{1,i}$. We observe
that $\phi_{|\al_i(M)}\circ \E_i\lel \phi$. By Connes's result
\cite{Co} we find \eqref{modul}. The equations \eqref{compat}
follow from $\al_i|_{\Ma}=id$ and the uniqueness of the
$\si_t^{\phi}$ invariant conditional expectation. Since $\M_i$ is
generated by $\al_1(\M),...,\al_i(\M)$ we deduce that $\M_i$ is
also $\si_t^{\phi}$ invariant. This allows us to apply Lemma to
apply Lemma \ref{pyind} i) and deduce \eqref{condi}.\qd

We will say that $(\Ma,\al_1(\M),...,\al_n(\M),\N,E)$ is a
\emph{conditioned system of top-subsymmetric copies} of in
addition
 \[ \si_t^{\phi_{\Ma}\circ E}(\al_i(M))\subset \al_i(M) \]
holds for a normal faithful state  $\phi_{\Ma}\circ E$ on $\Ma$.
We will keep  this (minimal) assumption for the rest of this
section and the notation $\phi=\phi_{\Ma}\circ E$. Let us note
that in particular, the von Neumann algebra $\M_i$ generated by
$\al_1(\M),...,\al_i(\M)))$ satisfies
 \[ \si_t^{\phi}(\M_i)\subset \M_i \]
In particular, we find $\si_t^{\phi}$-invariant  conditional
expectations $E_i:\N\to \al_i(\M)$ and $\E_i:\N\to \M_i$. Let us
fix some further notations which will be used in the proof. We
denote by $D$ the density of the faithful normal state $\phi$.
For an element $x\in \N$, we may define
 \[ \al_i(xD)\lel \al_i(x)D \pl .\]
Then
 \begin{align*}
  \noo \al_i(xD)\rrm_1 &= \sup_{\noo y\rrm\le 1}
 |\phi(y\al_i(x))| \lel
 \sup_{\noo y\rrm\le 1}
 |\phi(E_i(y)\al_i(x))| \lel
 \sup_{\noo y\rrm\le 1}
 |\phi(\E(y)x)| \lel  \noo xD\rrm_{L_1(\M)} \pl .
 \end{align*}
We refer to \cite{JX} for the natural extension of $E$ and $\E_i$,
$i=1,...,n$ on $L_p(\N)$. We refer to \cite{JD} for the space
$L_1^c(\M,E)$ to be the completion of $\M D$ under then norm

 \[ \noo xD\rrm_{L_1^c(\M,E)}\lel \noo
 DE(x^*x)D\rrm_{L_{\frac{1}{2}}(\Ma)}^{\frac12}  \lel \noo
 (DE(x^*x)D)^{\frac12}\rrm_{L_1(\Ma)} \pl
 .\]
It is easily checked that $L_1^c(\M,E)=L_2(\M)L_2(\Ma)$. Moreover,
we have a natural duality bracket between $L_1^c(\M,E)$ and the
Hilbert $C^*$-module $L_\infty(\M,E)$ defined as the completion of
$\M$ with respect to the norm $\noo x\rrm_{L_\infty^c(\M,E)}\lel
\noo E(x^*x)\rrm^{1/2}$. Indeed, we have
 \[ |tr(x^*y)|\kl \noo E(x^*x)\rrm_{\Ma}^{\frac12} \noo
 E(y^*y)^{\frac12}\rrm_{L_1(\Ma)} \pl .\]
For more information on this Cauchy-Schwarz type inequality and
the fact that $L_\infty(\M,E)$ embeds isometrically into the
antilinear dual of $L_1^c(\M,E)$, we refer to \cite{JD}. We will
also  use the notation $L_1^r(\M,E)=L_2(\Ma)L_2(\M)$ for the
completion of $D\M$ with respect to the norm $\noo
x\rrm_{L_1^r(\M,E)}=\noo x^*\rrm_{L_1^c(\M,E)}$. The following
3-term quotient norm is the central object in this section
 \begin{align*}
 \noo x\rrm_{\kz_{n,\eps}} \lel
  \inf_{x=x_1+x_2+x_3}  \kla n\noo x_1\rrm_{L_1(\M)}
     +\sqrt{n \eps } \noo E(x_2^*x_2)^{\frac12}\rrm_{L_1(\Ma)}
  +\sqrt{n \eps } \noo E(x_3x_3^*)^{\frac12}\rrm_{L_1(\Ma)}
  \mer
  \pl .
 \end{align*}
Here, we allow  $x_2\in L_1^c(\M,E)$ and $x_3\in L_1^r(\M,E)$. The
parameter $\eps>0$ will be chosen in a convenient way later. We
will also use the symbol $\kz_{n,\eps}(\M,E,\phi)$ for the space
$L_1(\M)$ equipped with this norm. In our context it is convenient
to use the following  antilinear duality
 \begin{equation} \label{brn}  \lb y,x\rb_n \lel n tr(y^*x) \pl.
 \end{equation}

The following Lemma is a minor modification of \cite[Lemma
6.9]{Joh}.\lz

\begin{lemma}\label{dualk} The dual space of $\kz_{n,\eps}(\M,E,\phi)$ with
respect to the duality bracket $\lb \pl \pl ,\pl \rb_n$ is
$\bar{\M}$ and the norm is given by
 \for
  \noo y\rrm_{\kz_{n,\eps}^*} &=&  \max\{ \noo y\rrm_{\M}, \eps^{-\frac12} \sqrt{n} \noo
  E(y^*y)\rrm_\Ma^{\frac12}, \eps^{-\frac12} \sqrt{n} \noo   E(yy^*)\rrm_\Ma^{\frac12} \}
  \pl .\
 \mel
\end{lemma}
The main inequality of this section is to show that for
$\eps=0.01$ we have
 \begin{equation}\label{main}
  \noo x\rrm_{\kz_{n,\eps}} \kl c \pl \ez \noo \summ_{i=1}^n \eps_i \al_i(x)
 \rrm_{L_1(\N)} \pl . \end{equation}
Here $(\eps_i)_{i=1}^n$ are independent Bernoulli variables, i.e.
$Prob(\eps_i=\pm 1)=\frac12$. Of course, we will first apply the
noncommutative Khintchine inequality (see \cite{LPP}). In
\cite{LPP} the  Khintchine inequality was formulated for the
characters  $(e^{2\pi \sqrt{-1}\pl 2^i})_i$.  Together with the
well-known contraction principle (see \cite{LT-I, LT-II}) we may
replace the characters $(e^{2\pi \sqrt{-1}\pl 2^i})_i$ by
Rademacher variables $(\eps_i)_i$. This implies
 \begin{eqnarray}\label{kh}
 \begin{minipage}{13cm} \vspace{-0.1cm}
 \for
 \ez \noo \summ_{i=1}^n \eps_i \al_i(x)
 \rrm_{1}
  &\le&   \inf_{\al_i(x)=c_i+d_i} \noo (\summ_{i=1}^n
  c_i^*c_i)^{\frac12}\rrm_{1} +
  \noo (\summ_{i=1}^n
   d_id_i^*)^{\frac12}\rrm_{1}  \\
  &\le& 2(1+\sqrt{2})   \pl   \ez \noo \summ_{i=1}^n \eps_i \al_i(x)
  \rrm_{1} \pl .
 \mel\end{minipage}
 \end{eqnarray}
Therefore a lower estimate for $\sum_{i=1}^n \eps_i \al_i(x)$ may
be obtained  by finding elements in $y_1,...,y_n\in \N$ whose row
and column norms are  controlled simultaneously. The first impulse
is to use the elements $y_1=\al_1(y)$,.....,$y_n=\al_n(y)$.
However, there we have no  good norm estimate for the square
function of these elements. The key idea in our proof is to define
the following elements starting from a contraction $y\in \M$:
 \begin{align*}
  a &= \sqrt{1-yy^*} \quad \mbox{and}\quad  b\lel    \sqrt{1-y^*y} \pl ,\\
  Y_i &=  \al_i(y)   \pl ,\\
   A_i&= \al_1(a)\al_2(a)\cdots \al_{i-1}(a) \quad \mbox{and}
   \quad
   B_i \lel     \al_{i-1}(b)\cdots \al_2(b) \al_1(b)  \pl .
 \end{align*}
Note  $a$, $b$ are well-defined and $A_i$ and $B_i$ are
contractions.

\begin{lemma} \label{updown}
With the definitions above, we have
 \[ \noo \summ_{i=1}^n (A_iY_iB_i)(A_iY_iB_i)^*\rrm_{\N} \kl 1
 \quad ,\quad
  \noo \summ_{i=1}^n (A_iY_iB_i)^*(A_iY_iB_i)\rrm_{\N} \kl 1  \pl .
 \]
\end{lemma}

\begin{proof} By symmetry it suffices to consider the first
term. Since the $B_i$'s are contractions, we have
 \for
 \summ_{i=1}^n (A_iY_iB_i)(A_iY_iB_i)^* &=&
  \summ_{i=1}^n A_iY_i B_iB_i^*Y_i^*A_i^* \kl
  \summ_{i=1}^n A_iY_iY_i^*A_i^*  \pl .
 \mel
We will inductively rewrite this sum as follows:
 \begin{align*}
 &1- \summ_{i=1}^n A_iY_iY_i^*A_i^*  \lel 1-\al_1(yy^*) -\summ_{i=2}^n
 A_iY_iY_i^*A_i^* \\
 &\pll \lel  1- \al_1(yy^*)
 -\summ_{i=2}^n \al_1(a) \al_2(a) \cdots \al_{i-1}(a) \al_i(yy^*)
  \al_{i-1}(a)\cdots \al_2(a) \al_1(a) \\
 &  \lel  \al_1(a)\al_1(a)
 - \al_1(a)\kla
 \summ_{i=2}^n  \al_2(a) \cdots \al_{i-1}(a) \al_i(yy^*)
 \al_{i-1}(a)\cdots \al_2(a) \mer \al_1(a) \\
 &\pll =  \al_1(a)\kla 1-\summ_{i=2}^n  \al_2(a) \cdots \al_{i-1}(a) \al_i(yy^*)
 \al_{i-1}(a)\cdots \al_2(a) \mer \al_1(a) \\
 &\pll =   ... \\
 &\pll =  \al_1(a) \al_2(a) \cdots \al_n(a)\al_n(a)\cdots \al_1(a)
 \gl 0 \pl .
 \end{align*}
This implies $\sum_{i=1}^n A_iY_iY_i^*A_i^* \kl 1$ and the
assertion is proved.\qd

The next (slightly technical) lemma explains why we need
independent copies.

\begin{samepage}\begin{lemma} \label{rechnen} Let  $\eps<e^{-1}$
and $\|y\|_{K_{n,\eps}^*}\kl 1$. Let $i\in \{1,..,n\}$. Then we
have
\begin{enumerate}
 \item[i)] $E(\al_1(a)\al_2(a)\cdots \al_i(a)) \lel E(a)^i$,
 $E(\al_i(a)\al_{i-1}(a)\cdots \al_1(a)) \lel E(a)^i$,  \newline
 $E(\al_i(b)\cdots \al_2(b)\al_1(b)) \lel E(b)^i$,
 $E(\al_1(b)\al_2(b) \cdots \al_i(b)) \lel E(b)^i$.
  \item[ii)] $\noo 1-E(a)\rrm \kl
  \frac{\eps}{n}$,  $\noo 1-E(b)\rrm \kl \frac{\eps}{n}$.
  \item[ iii)] $\|\sum_{i=1}^n 1-E(a)^{i-1}\|\kl e\eps n$,
  $\| \sum_{i=1}^n 1-E(b)^{i-1}\|\kl e\eps n$.
   \item[ iv)]
 \[  \noo \summ_{i=1}^n E\bigg( (1-A_i)
  \big( 1-A_i)^* \bigg) \rrm \kl 2e\eps n \quad, \quad
   \noo \summ_{i=1}^n E\bigg( \big(1-B_i\big)^* \big(1-B_i)\big) \bigg) \rrm \kl 2e\eps n \pl .\]
\end{enumerate}
\end{lemma}\end{samepage}

\begin{proof} It suffices to prove the inequalities involving
$a$'s. In order to prove i), we note that by
$\si_t^{\phi}$-invariance, the conditional expectation $E$ is
unique. Therefore $E\lel E\circ \E_i$. From \eqref{condi}  and
independence we deduce that
 \begin{align*}
  E(\al_1(a)\cdots \al_{i}(a)) &= E\bigg(\E_{i-1}\big(\al_1(a)\cdots
  \al_{i-1}(a)\big)\al_i(a)\bigg)
  \lel E\bigg(E\big(\al_1(a)\cdots \al_{i-1}(a)\big) \al_i(a)\bigg) \\
  &= E(\al_1(a)\cdots \al_{i-1}(a))E(a) \pl .
 \end{align*}
By induction, we get
 \[ E(\al_1(a)\cdots \al_{i}(a))  \lel E(a)^{i} \pl .\]
For the  proof of  ii), we note that for $0\le t\le 1$, we have
$1-\sqrt{1-t}\le t$. By functional calculus, we obtain from the
assumption on $y$
 \for
 1-E(a) \lel   1-E(\sqrt{1-yy^*}) &=& E(1-\sqrt{1-yy^*}) \kl   E(yy^*) \kl \frac{\eps}{n} \pl .
 \mel
In order to prove iii), we first recall  that for $t\le e^{-1}$,
we have $1-t\ge \exp(-et)$. By functional calculus, this implies
 \[ 1\ge E(a) \ge 1- \frac{\eps}{n} \ge \exp(-\frac{e\eps}{n})
 \pl . \]
Hence, for all $i=1,..,n$, we have
 \[ E(a)^{i}\ge \exp(-e\eps) \pl .\]
Using $1-\exp(-t)\le t$, this yields
 \[ \summ_{i=1}^n [1-E(a)^{i-1}] \kl n(1-\exp(-e\eps)) \kl n e\eps
 \pl . \]
For the proof of the last assertion iv), we deduce with   $\noo
a\rrm\le 1$ and i)

 \begin{samepage}\begin{align*}
 &E\big( (1-\al_1(a)\cdots \al_{i-1}(a)) (1-\al_1(a)\cdots
  \al_{i-1}(a))^*\big) \\
  &= 1- E\big(\al_1(a)\cdots
  \al_{i-1}(a)\big)  -   E\big(\al_{i-1}(a)\cdots \al_1(a)\big)
 +
   E\big (\al_1(a)\cdots \al_{i-1}(a)\al_{i-1}(a)\cdots \al_1(a)\big)
   \\
  &\le  2- E\big(\al_1(a)\cdots \al_{i-1}(a)\big) -
  E\big(\al_{i-1}(a)\cdots \al_1(a)\big) \\
  &  =  2(1-E(a)^{i-1}) \pl.
  \end{align*}\end{samepage}
Therefore, iii) implies
 \begin{align*}
  \noo \summ_{i=1}^n  E\bigg( \big( 1-\al_1(a)\cdots \al_{i-1}(a)\big) \big(1-\al_1(a)\cdots
  \al_{i-1}(a)\big)^*\bigg) \rrm &\le  2ne \eps \pl . \qedhere
   \end{align*}
\qd

\begin{lemma}\label{expest} Let $E_n:\N\to \al_n(M)$ be the
$\si_t^{\phi}$ invariant conditional expectation and $A_i$, $B_i$
as above.  Let  $e\eps \le 2$ and $z\in \M$. Then
 \[   \noo \summ_{i=1}^n \al_n^{-1}E_n\bigg((1-A_i)\al_n(z)(1-B_i)\bigg)\rrm_{\kz_{n,\eps}^*} \kl
 2ne\sqrt{\eps}
  \noo z\rrm_{\kz_{n,\eps}^*} \pl .\]
\end{lemma} \lz

\begin{proof}  We first consider the norm estimate in $\M$.
By   the Cauchy-Schwarz inequality \cite{JD} and Lemma \ref{pyind}
i), we deduce
 \begin{align*}
 &\noo \summ_{i=1}^n \al_n^{-1}E_n\big((1-A_i)\al_n(z)(1-B_i)\big)\rrm_{\N}  \\
 &\pl \kl
 \noo \summ_{i=1}^n E_n\big((1-A_i)(1-A_i)^*\big)\rrm_{\N}^{\frac12}
 \noo \summ_{i=1}^n E_n\big((1-B_i)^*\al_n(z^*z)(1-B_i)\big)\rrm_{\N}^{\frac12} \\
 &\pl   \kl  \noo z^*z\rrm^{\frac12} \pl   \noo \summ_{i=1}^n
 E_n\big((1-A_i)(1-A_i)^*\big)\rrm_{\N}^{\frac12}
 \noo \summ_{i=1}^n
 E_n\big((1-B_i)^*(1-B_i)\big)\rrm_{\N}^{\frac12}\\
 & \pl \lel \noo z^*z\rrm^{\frac12} \pl   \noo \summ_{i=1}^n
 E\big((1-A_i)(1-A_i)^*\big)\rrm_{\N}^{\frac12}
 \noo \summ_{i=1}^n
 E\big((1-B_i)^*(1-B_i)\big)\rrm_{\N}^{\frac12}
 \pl .
 \end{align*}
Here we used $1-A_i,1-B_i\in \M_{n-1}$ for all $i=1,...,n$. Then
Lemma \ref{rechnen} implies that
 \begin{equation}\label{above}
  \noo \summ_{i=1}^n E((1-A_i)(1-A_i)^*)\rrm
  \kl 2en \eps \quad \mbox{and}\quad
   \noo \summ_{i=1}^n E((1-B_i)^*(1-B_i))\rrm  \kl 2en \eps \p.
  \end{equation}
Thus we get
 \[ \noo \summ_{i=1}^n \al_{n}^{-1}E_n((1-A_i)\al_n(z)(1-B_i))\rrm_{\M}
  \kl 2e\eps n \noo z\rrm_{\M} \pl .\]
Now, we prove the estimates for the conditioned norm. We observe
that
 \[ B \lel   \sum_{i=1}^n e_{i,1}\ten (1-B_i)\in
 \Ma_n(\M_{n-1}) \pll \mbox{and} \pll  A\lel \sum_{i=1}^n
 e_{1,i}\ten (1-A_i)\in \Ma_n(\M_{n-1}) \pl .\]
Since the $A_i$'s are  contractions, we deduce
 \[ \noo A\rrm\kl \sqrt{n} \sup_{i=1,..,n}\noo 1-A_i\rrm \kl
 2\sqrt{n} \pl .\]
Similarly, $\noo B\rrm\le 2\sqrt{n}$. Now we apply  Kadison's
inequality $E_n(x)^*E_n(x)\le  E_n(x^*x)$ (see \cite{Pau}), Lemma
\ref{pyind} iii)  and \eqref{above}:
 \begin{align*}
  &\noo E\Bigg (
  \bigg(\summ_{i=1}^n \al_n^{-1}E_n\big((1-A_i)\al_n(z)(1-B_i)\big)\bigg)^* \bigg(\summ_{i=1}^n
  \al_n^{-1}E_n\big((1-A_i)\al_n(z)(1-B_i)\big)\bigg)\Bigg)\rrm \hspace{1cm} {\atop }\\
  &\pll =    \noo E \bigg(E_n\big(B^*(1\ten \al_n(z)^*)A^*\big)E_n\big(A(1\ten \al_n(z))B\big)\bigg)
  \rrm\\
  &\pll
  \kl    \noo E \circ  E_n\bigg(B^*(1\ten \al_n(z)^*)A^*A(1\ten \al_n(z))B\bigg)
  \rrm \kl
   \noo A^*A\rrm  \pl  \noo E(B^*(1\ten \al_n(z)^*\al_n(z))B)  \rrm \\
   &\pll \kl
       4n \pl  \noo \summ_{i=1}^n E((1-B_i)^*(1-B_i))  \rrm \noo
       E(\al_n(z^*z))
       \rrm \kl
   8n^2 e \eps  \noo E(z^*z) \rrm   \pl .
 \end{align*}
The argument in the row case is identical. Since $\sqrt{2e}\le e$
and $\eps\le \sqrt{\eps}$,   we get
 \[ \noo \summ_{i=1}^n \al_n^{-1}E_n((1-A_i)\al_n(z)(1-B_i))\rrm_{\kz_{n,\delta}^*} \kl
  2ne\sqrt{\eps}  \noo z\rrm_{\kz_{n,\delta}^*} \pl \]
for all $\delta>0$. Using $\delta=\eps$ yields the assertion. \qd

The next argument provides  the lower estimate.

\begin{prop} \label{low1} Let $\eps=\frac{1}{100}$ and $N$ be a von Neumann algebra with normal faithful tracial state
$\tau$. Let $v_1,...,v_n$ be unitaries in $N$. Then for every
$x\in L_1(\N)$
 \[ \noo x\rrm_{\kz_{n,\eps}}  \kl 4 \pl \inf_{v_i.\tau \ten \al_i(x)=c_i+d_i}
  \noo \kla \summ_{i=1}^n c_i^*c_i\mer^{\frac12} \rrm_{L_1(N\ten \N)} +
  \noo \kla \summ_{i=1}^n d_id_i^*\mer^{\frac12} \rrm_{L_1(N\ten \N)} \pl
  .
  \]
\end{prop}

\begin{proof} We may assume $\noo x\rrm_{\kz_{n,\eps}}=1$.  We
denote the right hand side without the factor $4$ by $INF$. By the
Hahn-Banach theorem, we can find $y\in \M$ such that
 \[ ntr(y^*x) \lel 1 \quad\mbox{and} \quad \noo y\rrm_{\kz_{n,\eps}^*} \lel
 1 \pl .\]
For any decomposition $v_i.\tau\ten \al_i(x)=c_i+d_i$, we deduce
from Lemma \ref{updown}
  \begin{align*}
 &\bet \summ_{i=1}^n tr\bigg((v_i\ten A_iY_iB_i)^*(v_i\tau \ten
\al_i(x))\bigg) \rag
 \lel  \bet \summ_{i=1}^n tr\bigg((v_i\ten A_iY_iB_i)^*(c_i+d_i)\bigg) \rag
 \end{align*}

 \begin{samepage}\begin{align*}
  & \le   \max\left \{ \noo \summ_{i=1}^n (A_iY_iB_i)^*(A_iY_iB_i)\rrm
  ^{\frac12}, \noo \summ_{i=1}^n
  (A_iY_iB_i)(A_iY_iB_i)^*\rrm^{\frac12} \right \}
  \\
  &  \pll \pll \pll \pll \pll \pll \pl
  \kla \| \big(\summ_{i=1}^n c_i^*c_i\big)^{\frac12} \|_{L_1(N\ten \N)} +
  \| \big(\summ_{i=1}^n d_id_i^*\big)^{\frac12} \|_{L_1(N\ten \N)}
   \mer \kl   INF \pl .
  \end{align*}\end{samepage}
We will now complete the argument using an error estimate. For
general $y$, $a,b$ we will use  the following algebraic identity.
 \begin{align*}
  y- ayb &=
  y-yb +(1-a)yb \lel y(1-b) + (1-a)y - (1-a)y(1-b) \pl.
 \end{align*}
Let us introduce $\beta_i(x)\lel v_i.\tau\ten \al_i(x)$ and
$y_i=v_i\ten Y_i$. Then, we deduce from the above that
 \begin{align*}
 \bet n tr(y^*x) \rag  &= \bet \summ_{i=1}^n tr(Y_i^*\al_i(x))
 \rag
 \lel \bet \summ_{i=1}^n tr(y_i^*\beta_i(x)) \rag
 \\
 &\le \bet \summ_{i=1}^n tr\big((A_iy_iB_i)^*\beta_i(x)\big) \rag + \bet \summ_{i=1}^n
  tr\big(y_i^*\beta_i(x)\big)-tr\big((A_iy_iB_i)^*\beta_i(x)\big)\rag
  \\
  &\le INF + \bet \summ_{i=1}^n tr\big(((A_i-1)y_i(1-B_i))^*\beta_i(x)\big) \rag \\
  &  \pll \pll + \bet \summ_{i=1}^n tr\big(((1-A_i)y_i)^*\beta_i(x)\big) \rag
  +   \bet \summ_{i=1}^n tr\big((y_i(1-B_i))^*\beta_i(x)\big)
 \rag \pl .
 \end{align*}
The three error terms will be treated separately. In order to
estimate the first term on the last line,  we remind the reader of
\eqref{brn} and that the trace is invariant under the conditional
expectations. Using  Lemma \ref{rechnen} i) together with the fact
that $(\M,\noo \pl \rrm_{\kz_{n,\eps}^*})$ is a left  $\Ma$
module, we obtain  that
 \begin{align*}
 &|\summ_{i=1}^n tr\big(y_i^*(1-A_i^*)\beta_i(x)\big)|
 \lel
  |\summ_{i=1}^n tr\big(\al_i(y)^*(1-A_i^*)\al_i(x)\big)| \lel
   |\summ_{i=1}^n tr\big(\al_i(xy^*)\E_i(1-A_i^*)\big)| \\
  &=
    |\summ_{i=1}^n tr(\al_i(xy^*)E(1-A_i^*))| \lel
 |\summ_{i=1}^n tr(y^*E(1-A_i^*)x)|  \lel
   |tr(y^*\big(\summ_{i=1}^n [1-E(a)^{i-1}]\big)x)| \\
   &\le
      \frac{1}{n}  \noo \big(\summ_{i=1}^n [1-E(a)^{i-1}] \big)y\rrm_{\kz_{n,\eps}^*}
  \noo x\rrm_{\kz_{n,\eps}} \kl  e\eps \pl .
 \end{align*}
Similarly, we obtain
 \for
  \bet  \summ_{i=1}^n tr\big((y_i(1-B_i))^*\beta_i(x)) \rag
  &\le&  \eps e \pl .
  \mel
For the remaining term we first top-subsymmetry  and then Lemma
\ref{expest}. This yields that
 \begin{align*}
 & \bet  \summ_{i=1}^n tr\bigg(((1-A_i)y_i(1-B_i))^*\beta_i(x)\bigg)
  \rag
   \lel
 \bet  \summ_{i=1}^n tr\bigg(((1-A_i)\al_i(y)(1-B_i))^*\al_i(x)\bigg)  \rag
 \\
 & =
   \bet  \summ_{i=1}^n tr\bigg(
   ((1-A_i)\al_n(y)(1-B_i))^*\al_n(x)\bigg)
  \rag \lel
   \bet  \summ_{i=1}^n tr\bigg(E_n\big((1-A_i)\al_n(y)(1-B_i))^*\big)\al_n(x)\bigg)
  \rag\\
  &\le
    \frac{1}{n} \noo \summ_{i=1}^n
  \al_n^{-1}E_n\big((1-A_i)\al_n(y)(1-B_i)\big)\rrm_{\kz_{n,\eps}^*} \noo
  x\rrm_{\kz_{n,\eps}}
   \kl 2e\sqrt{\eps} \pl .
 \end{align*}
Hence, we deduce that
 \for
 1 &=& |ntr(xy^*)| \kl INF + 2e\eps + 2e\sqrt{\eps} \pl .
 \mel
It turns out that for $\eps=0.01$, we have
$2e(\eps+\sqrt{\eps})\le \frac{3}{4}$ and thus $\frac14 \kl INF$.
This completes  the proof.\qd

For the sake of completeness, we also prove the converse
implication.

\begin{lemma} \label{up1}
Let $x\in L_1(\M)$. Then
 \for
 \lefteqn{ \inf_{v_i.\tau \ten \al_i(x) \lel c_i+d_i}
  \noo \kla \summ_{i=1}^n c_i^*c_i\mer^{\frac12} \rrm +
  \noo \kla \summ_{i=1}^n d_id_i^*\mer^{\frac12} \rrm }\\
  &\le& \inf_{x=x_1+x_2+x_3} n \noo x_1 \rrm_{L_1(\M)}
  + \sqrt{n} \noo E(x_2^*x_2)^{\frac12} \rrm_{L_1(\Ma)} +
  \sqrt{n} \noo E(x_3x_3^*)^{\frac12} \rrm_{L_1(\Ma)}
   \pl .
 \mel
\end{lemma}\lz

\begin{proof} We define $c_i= v_i.\tau \ten \al_i(x)$. Then
we have
 \for
  \noo \kla \summ_{i=1}^n c_i^* c_i\mer^{\frac12} \rrm_1
  &\le& \summ_{i=1}^n \noo v_i\tau \ten \al_i(x) \rrm \lel n \noo
  x\rrm_{L_1(\M)}  \pl.
 \mel
Moreover, let us recall that $\noo c^*c \rrm_{\frac12}\le \noo
E(c^*c)\rrm_{\frac12}$ (see \cite{JD}). Since the $v_i$'s are
unitaries, we have
  \begin{align*}
   \noo  \summ_{i=1}^n c_i^* c_i\rrm_{\frac12}
   &=  \noo  \summ_{i=1}^n \al_i(x)^*\al_i(x)\rrm_{\frac12}
   \kl   \noo \summ_{i=1}^n E(\al_i(x)^* \al_i(x))  \rrm_{\frac12}
   \lel    n \noo E(x^*x) \rrm_{\frac12}  \pl.
   \end{align*}
Similarly, we deduce $\noo \sum_{i=1}^n c_i c_i^*
\rrm_{\frac12}\le
 n \noo E(xx^*) \rrm_{\frac12}$. Hence for any decomposition $x=x_1+x_2+x_3$ we deduce from the
triangle inequality for $c_i\lel v_i.\tau\ten \al_i(x_1+x_2)$ and
$d_i=v_i\tau\ten \al_i(x_3)$ that
 \begin{align*}
 \noo (\summ_i c_i^*c_i)^{\frac12}\rrm +
 \noo (\summ_i d_id_i^*)^{\frac12}\rrm &\le
   n\noo x_1\rrm+ \sqrt{n}\noo E(x_2^*x_2)^{\frac12}\rrm +
 \sqrt{n}\noo E(x_3x_3^*)^{\frac12}\rrm \pl . \qedhere
 \end{align*}
\qd

A combination of  Proposition \ref{low1} and Lemma \ref{up1}
yields the main results of this section.\lz

\begin{theorem} \label{main2} Let $(\Ma,\M,\N,(\al_{i})_{i=1}^n,E,\phi)$ be   a system of conditioned top-subsymmetric
independent copies. Let $N$ be a von Neumann algebra with
normalized faithful trace $\tau$ and $v_i$ be unitaries in $N$.
Let $x\in L_1(\M)$. Then
%\begin{samepage}
\for
 \lefteqn{ \frac{1}{40}  \inf_{x=x_1+x_2+x_3}  \kla n\noo x_1\rrm_{L_1(\N)}
   +\sqrt{n} \noo E(x_2^*x_2)^{\frac12}\rrm_{L_1(\Ma)}
   +\sqrt{n} \noo E(x_3x_3^*)^{\frac12}\rrm_{L_1(\Ma)} \mer } \\
 & & \kl
 \inf_{v_i\tau \ten \al_i(x) \lel c_i+d_i}
  \noo \kla \summ_{i=1}^n c_i^*c_i\mer^{\frac12} \rrm +
  \noo \kla \summ_{i=1}^n d_id_i^*\mer^{\frac12} \rrm \\
  & & \kl  \inf_{x=x_1+x_2+x_3} n \noo x_1\rrm_{L_1(\N)}
  + \sqrt{n} \noo E(x_2^*x_2)^{\frac12} \rrm_{L_1(\Ma)} +
  \sqrt{n} \noo E(x_3x_3^*)^{\frac12} \rrm_{L_1(\Ma)}
   \pl .
 \mel
% \end{samepage}
\end{theorem}

\begin{proof}  For $\eps=0.01$,  we have  $\noo x
\rrm_{\kz_{n,1}}\le 10 \noo x \rrm_{\kz_{n,\eps}}$. Hence
Proposition \ref{low1} implies the first inequality and the
converse inequality is Lemma \ref{up1}. \qd

\begin{cor} \label{main1} Let $(\Ma,\M,\N,(\al_{i})_{i=1}^n,E,\phi)$ be   a system of conditioned top-subsymmetric
independent copies.  Let $N$ be a von Neumann algebra with
normalized faithful trace $\tau$ and $v_i$ unitaries in $N$. Let
$x\in L_1(\M)$ and $(\eps_i)$ be independent Bernoulli random
variables. Then
 \for
  & & \vspace{-0.5cm}  \inf_{x=x_1+x_2+x_3}  \kla n\noo x_1\rrm_{L_1(\N)}
   +\sqrt{n} \noo E(x_2^*x_2)^{\frac12}\rrm_{L_1(\Ma)}
   +\sqrt{n} \noo E(x_3x_3^*)^{\frac12}\rrm_{L_1(\Ma)} \mer  \\
 & & \pll\pll   \sim_{200} \ez \noo \summ_{i=1}^n \eps_i v_i\tau \ten \al_i(x)
 \rrm_{L_1(\N)} \pl .
 \mel
\end{cor}

\begin{proof} This follows immediately from the noncommutative
Khintchine   inequality, see \eqref{kh}.\qd

Our main application is obtained for tensor products.

\begin{samepage}\begin{cor}\label{mainten} Let $\Ma$ and $M$ be von Neumann
algebras with  normal, faithful states $\phi_{\Ma}$ and $\psi$,
respectively. Let $N$ be a finite von Neumann algebra with
faithful normal trace $\tau$ and $v_1,...,v_n\in N$ be unitaries.
Let $D_n$ be the density of $\phi_{\Ma}\ten \psi^{\ten_n}$, $D$ be
the density of $\phi_{\Ma}\ten \psi$. Then
 \begin{align*}
  &\ez \noo \summ_{k=1}^n \eps_k  \pl v_k.\tau \ten (1\ten
  \pi_k)(y)D_n\rrm_{L_1(N\ten \Ma\ten M^{\ten_n})}
  \pl \pl \sim_{200}\pl
 \inf_{yD=y_1D +y_2D +Dy_3}
     n\noo y_1D\rrm_1 \\
  &\quad \quad  \pll  +
  \sqrt{n} \noo (DE(y_2^*y_2)D)^{\frac12}
  \rrm_{1} +
   \sqrt{n} \noo (DE(y_3y_3^*)D)^{\frac12}\rrm_{1}
 \end{align*}
holds for every $y\in (\Ma\ten M)_{a}$. Here the infimum is taken
over analytic elements.
\end{cor}\end{samepage}

\begin{proof} Let us assume $200 \pl  \ez \noo \sum_{k=1}^n \eps_k \pl  v_k.\tau \ten (1\ten
  \pi_k)(y)D\rrm_{L_1(N\ten \Ma\ten M^{\ten_n})}\le 1$.
For $x=yD\in L_1(\Ma\ten M)$ we  consider a decomposition
 \[ yD\lel x_1+x_2+x_3 \]
where $x_1$, $x_2$ and $x_3$ satisfy the corresponding norm
estimates given by Theorem \ref{main}. Then we may approximate
$x_2$ and $x_3$ by elements of the form
 \[ x_2\lel y_2D+x_2' \quad \mbox{and} \quad  x_3\lel Dy_3+ x_3'
 \]
such that $y_2$ and $y_3$ are analytic and
 \[ \noo x_2'\rrm_{L_1^c(\Ma\ten M,E)}\kl \frac{\eps}{n}
 \quad \mbox{and} \quad  \noo x_3'\rrm_{L_1^r(\Ma\ten M,E)}\kl \frac{\eps}{n} \pl
 .\]
This yields
   \[   yD\lel x_1+x_2'+x_3'+ y_2D+Dy_3  \quad \mbox{and} \quad
 x_1+x_2'+x_3' \lel (y-y_2-\si_{-i}(y_3))D \pl .\]
Since the inclusions $L_1^c(\Ma\ten M,E)\subset L_1(\Ma\ten M)$
and $L_1^r(\Ma\ten M,E)\subset L_1(\Ma\ten M)$ are contractive, we
deduce
 \begin{align*}
 & n \noo x_1+x_2'+x_3'\rrm_1 \kl
   n\noo x_1\rrm_1+ n \noo
 x_2'\rrm_{L_1^c(\Ma\ten M,E)} + n \noo x_3'\rrm_{L_1^r(\Ma\ten
 M,E)} \kl   (1+2\eps) \pl .
 \end{align*}
Thus the assertion follows with $y_1=y-y_2-\si_{-i}(y_3)$.\qd

\section{Khintchine type inequalities}

In this section we combine the central limit procedure with the
norm estimates from the previous section. Let us  recall that for
a sequence $\mu=(\mu_k)\subset (0,1)$ the von Neumann algebra
$\N(\mu)$ is the completion of $\ten_{k} \Mz_2$ with respect to
the GNS construction of the tensor product state
$\phi_{\mu}=\ten_k \phi_{\mu_k}$.

\begin{theorem}\label{carkh} Let $(a_k)$ be the generators of the CAR
algebra, $\phi_{\mu}$ the quasi-free state satisfying
\eqref{carrk} and $D_{\mu}\in L_1(\N(\mu))$ the density of
$\phi_{\mu}$. Let $\Mz$ be a von Neumann algebra and $x_k\in
L_1(\Mz)$. Then
 \begin{align*}
  &\noo \summ_k
  D_{\mu}^{\frac12}a_kD_{\mu}^{\frac12}\ten x_k\rrm_{L_1(\N(\mu))\wet
  L_1(\Mz)} \!\!\!\!
  \sim \!\!\!\!
   \inf_{x_k\lel c_k+d_k}
 \noo (\summ_k (1-\mu_k) c_k^*c_k)^{\frac12}\rrm_{L_1(\Mz)}+
   \noo (\summ_k \mu_k d_kd_k^*)^{\frac12}\rrm_{L_1(\Mz)} \!\!\!\! .
   \end{align*}
\end{theorem}

We will need some notations from \cite{Joh}. For a Hilbert space
$H$ and a von Neumann algebra $\M$ we denote by $H^r\wet L_1(\M)$
the closure of $H\ten L_1(\M)$ with respect to norm
 \[ \noo x\rrm_{H^r\wet L_1(\M)} \lel \noo (x^*,x)\rrm_{\frac12}^{\frac12} \pl. \]
Here the $L_{\frac12}(\M)$-valued scalar product is defined by
 \[ (\summ_k h_k\ten x_k,\summ_j l_j\ten y_j)
  \lel \summ_{k,j}  (h_k,l_j)x_ky_j  \pl \]
and the $^*$ operation  is given by $
 (\sum_k h_k\ten x_k)^*\lel \sum_k h_k\ten x_k^*$.
We refer the reader to \cite[section1]{Joh} for more details on
the operator space structure and to
 \cite{JS} for more information on $L_{\frac p2}$-valued
scalar products. Our motivation for this abstract definition is
the Hilbert space $\ell_2(\mu)$ with scalar product $(h,l)=\sum_k
\mu_k \bar{h}_kl_k$. Let us denote by $(e_k)$ the standard unit
vectors basis. Then
 \begin{equation}\label{muk}
  \noo \summ_k e_k\ten x_k \rrm_{\ell_2^r(\mu)\wet L_1(\M)}
  \lel \noo (\summ_k \mu_k x_k^*x_k)^{\frac12} \rrm_{L_1(\M)} \pl
  .
  \end{equation}
Similarly, we recall that $H^c\wet L_1(\M)$ is the completion of
$H\ten L_1(\M)$ with respect to the norm $\noo x\rrm_{H^c\wet
L_1(\M)}=\noo
 \overline{(x,x^*)}^{\frac12}\rrm$. More precisely, for $x=\sum_k
 h_k\ten x_k$ we have
  \[ \overline{(x,x^*)}\lel \sum_{k,l} \overline{(h_k,h_l)}
  x_kx_l^* \lel  \sum_{k,l} (h_l,h_k) x_kx_l^* \pl .\]
For $H=\ell_2$ this is consistent with $H^c\wet L_1(\M)\cong
\{\sum_k e_{1k}\ten x_k:x_k \in L_1(\M)\}\subset L_1(B(\ell_2)\ten
\M)$. In particular, we have
 \begin{equation}\label{muk2}
  \noo \summ_k e_k\ten x_k \rrm_{\ell_2^c(\mu)\wet L_1(\M)}
  \lel \noo (\summ_k \mu_k x_kx_k^*)^{\frac12} \rrm_{L_1(\M)} \pl
  .
  \end{equation}
The proof of Theorem \ref{carkh} will use the central limit
procedure and also work for the $q$-commuta- tion relations for
$-1\le q\le 1$. By approximation, we shall  assume that $\Mz$ is
$\si$-finite. Again by approximation, it suffices to prove Theorem
\ref{carkh} for a finite number $m$ of generators. This leads  to
the following data: $N=\ell_{\infty}^m(\Mz_2)$ with weight
 \[ \psi(x)\lel \summ_{k=1}^m [(1-\mu_k)x_{11}(k)+
 \mu_kx_{22}(k)] \pl .\]
We define $T=m$ and
 \[ u_{n,T}(x) \lel \sqrt{\frac{T}{n}} \summ_{i=1}^n v_i\ten
 \pi_i(x) \]
where $v_i\in \Mz_{2^n}$ are Clifford matrices  satisfying
$v_iv_j\lel-v_jv_i$. We use the normal faithful state
$\phi_n=\tau_{n}\ten (\frac{\psi}{T})^{\ten_n}$ on
$N_n=\Mz_{2^n}\ten N^{\ten_n}$. We recall that on the Haagerup
$L_2$ space $L_2(N,tr)$ we use $(x,y)=tr(x^*y)$.

\begin{samepage} \begin{lemma}\label{ne1} Let $\phi_{\Mz}$ be a normal
faithful state on $\Mz$.  Let $D$  be the density of
$\phi_{\Mz}\ten \frac{\psi}{T}$, $D_{\Mz}$ be the density of
$\phi_{\Mz}$,  and $D_n$ be the density of $\phi_n\ten
\phi_{\Mz}$. Let $y\in N\ten \Mz$. Then
\begin{align*}
  &\noo  D_n^{\frac12}(
  u_{n,T}\ten id)(y)D_n^{\frac12}\rrm_{L_1(N_n\wet \Mz)}
   \pl \sim_{200} \inf_{y=y_1+y_2+y_3}
  \sqrt{Tn} \noo D^{\frac12}y_1D^{\frac12}\rrm_1 \\
  &\quad +
  \noo (D_{\psi}^{\frac12}\ten D_{\phi_{\Mz}}^{\frac12})y_2(1\ten D_{\phi_{\Mz}}^{\frac12})\rrm_{L_2^r(N,tr)\wet L_1(\Mz)}
  + \noo (1\ten D_{\phi_{\Mz}}^{\frac12})y_3(D_{\psi}^{\frac12}\ten D_{\phi_{\Mz}}^{\frac12})\rrm_{L_2^c(N,tr)\wet L_1(\Mz)}
  \pl.
  \end{align*}
\end{lemma} \end{samepage}

\begin{proof} Let us first mention the well-known absorption principle
 \[ \ez \noo \summ_k \eps_k \ten v_k \ten y_k\rrm_1 \lel
 \noo \summ_k v_k \ten y_k\rrm_1 \pl .\]
Indeed, the unitaries $w_k\lel \eps_k\ten v_k\in
L_{\infty}(\{-1,1\}^n)\ten \Mz_{2^n}$ also satisfy the CAR
relations  $w_kw_j\lel -w_jw_k$, $w_k=w_k^*$  and $w_k^2=1$.
Therefore $\pi(v_k)=w_k$ extends to a trace preserving
homomorphism from $\Mz_{2^n}$ onto the algebra $B_n$ generated by
the $w_k$. Using the trace preserving conditional expectation onto
$B_n$, we deduce
 \begin{align*}
 \ez \noo \summ_k \eps_k \ten v_k \ten y_k\rrm_{L_1(N_n)\wet L_1(\Mz)}\!\!\! &=
 \noo \summ_k w_k  \ten y_k\rrm_{L_1(B_n)\wet
 L_1(\M)} \!\!\!\lel
 \noo \summ_k v_k \ten y_k\rrm_{L_1(\Mz_{2^n})\wet L_1(\M)} \pl
\!\!\! .
 \end{align*}
By approximation it suffices to consider analytic elements $y \in
N\ten \Mz$. We deduce from Corollary \ref{mainten} applied to
$z=\si_{\frac{-i}{2}}(y)$ that
 \begin{align*}
 &\noo  D_n^{\frac12}(id\ten
  u_{n,T})(y)D_n^{\frac12}\rrm_{L_1(\Mz\wet N_n)} \lel
  \sqrt{\frac{T}{n}} \noo \summ_{i=1}^n v_i.\tau \ten D_n^{\frac12}(1\ten
  \pi_i(y))D_n^{\frac12} \rrm \\
  &\pl = \sqrt{\frac{T}{n}} \ez \noo \summ_{i=1}^n \eps_i v_i.\tau \ten (1\ten
  \pi_i(z))D_n \rrm\\
%  \end{align*}
%   \begin{align*}
  &  \sim_{200} \sqrt{\frac{T}{n}}\kla  \inf_{zD=z_1D+z_2D+Dz_3}
   n \noo z_1D\rrm_1
 % \right. \\&\quad \quad \quad  \left.
  +  \sqrt{n} \noo
   (DE(z_2^*z_2)^{\frac12}D)^{\frac12}\rrm_{1} + \sqrt{n}
   \noo
   (DE(z_3z_3^*)^{\frac12}D)^{\frac12}\rrm_1 \mer  \pl .
  \end{align*}
Here $E(a\ten z_2)\lel T^{-1}\psi(a)z_2$ is the conditional
expectation onto $\Mz$. Given $z=\sum_k b_k\ten z_k$ we observe
that
 \begin{align*}
  &\noo z(D_{\psi}^{\frac12}\ten D_{\phi_{\Mz}})\rrm_{L_2^r(N,tr)\wet L_1(\Mz)}
   \lel \noo
  \big((D_{\psi}^{\frac12}\ten D_{\phi_{\Mz}})z^*,z(D_{\psi}^{\frac12}\ten
  D_{\phi_{\Mz}})\big)\rrm_{\frac12}^{\frac12}\\
  &=  \noo \sum_{k,l} \psi(b_k^*b_l) D_{\phi_{\Mz}}z_k^*z_lD_{\phi_{\Mz}}
 \rrm_{\frac12}^{\frac12} \lel
   \sqrt{T} \noo D E(z^*z)
 D \rrm_{\frac12}^{\frac12} \lel
 \sqrt{T} \noo (D E(z^*z)D)^{\frac12}\rrm_1
 \pl .
  \end{align*}
A similar calculation shows that for $z=\sum_k b_k\ten z_k$ we
have
 \begin{align*}
 &\sqrt{T} \noo D E(zz^*) D \rrm_{\frac12}^{\frac12}
 \lel \|\sum_k
 \psi(b_kb_l^*)D_{\phi_{\Mz}}z_kz_l^*D_{\phi_{\Mz}}\|_{\frac12}^{\frac12}
 \\
 &= \|\sum_k
 (D_{\psi}^{\frac12}b_l,D_{\psi}^{\frac12}b_k)_{L_2(N,tr)}
 D_{\phi_{\Mz}}z_kz_l^*D_{\phi_{\Mz}}\|_{\frac12}^{\frac12}
 \lel \| (D_{\psi}^{\frac12}\ten D_{\phi_{\Mz}})z\|_{L_2^c(N,tr)\wet
 L_1(\Mz)} \pl .
 \end{align*}
Since the infimum is taken over analytic elements, we may define
$y_1=\si_{i/2}(z_1),y_2=\si_{i/2}(z_2)$ and $y_3=\si_{-i/2}(z_3)$
such that
 \[ D^{\frac12}yD^{\frac12}\lel zD \lel
  z_1D+z_2D+Dz_3
 \lel D^{\frac12}y_1D^{\frac12} + D^{\frac12}y_2D^{\frac12} +
 D^{\frac12}y_3D^{\frac12}  \pl . \]
This implies $y=y_1+y_2+y_3$. For arbitrary elements the assertion
follows by approximation. \qd

The following Lemma is an easy adaptation of \cite[Proposition
6.2]{Joh}. Indeed, in the proof we may use that $N$ is finite
dimensional and that the  modular operator $S(x)=x^*$ is bounded
on $L_2(N,\psi)$. Then the proof given for $\Mz=\Mz_m$ in
\cite{Joh} is easily adapted to this setting.

\begin{lemma}\label{limitn} Keep the notations from the previous
Lemma. Let $\U$ be an ultrafilter on the integers. Then
 \begin{align*}
  &\lim_{n,\U}
  \noo  D_n^{\frac12}(
  u_{n,T}\ten id)(y)D_n^{\frac12}\rrm_{L_1(N_n\wet \Mz)} \sim_{200} \inf_{y=y_2+y_3}\\
  & \pl
  \noo (D_{\psi}^{\frac12}\ten D_{\phi_{\Mz}}^{\frac12})y_2(1\ten D_{\phi_{\Mz}}^{\frac12})\rrm_{L_2^r(N,tr)\wet L_1(\Mz)}
  + \noo(1\ten  D_{\phi_{\Mz}}^{\frac12})y_3(D_{\psi}^{\frac12}\ten D_{\phi_{\Mz}}^{\frac12})\rrm_{L_2^c(N,tr)\wet L_1(\Mz)}
 \! \! \! \!  \pl.
  \end{align*}
\end{lemma}

\begin{proof}[Proof of Theorem \ref{carkh}.] We fix $m\in \nz$. We
denote by $(\delta_k)$ the unit vector basis in $\ell_{\infty}^m$
and $N=\ell_{\infty}^m(\Mz_2)$. We consider the selfadjoint
subspace $S\subset N$ generated by the elements $x_k=\delta_k\ten
e_{12}$. We apply the central limit procedure Theorem
\ref{limitth} to $u_{n,T}(x)\lel T^{1/2}n^{-1/2}\sum_{j=1}^n
v_j\ten \pi_j(x)$ where $v_jv_l=-v_jv_j$ are anticommuting
unitaries. We obtain an ultraproduct state $\phi_\U$ and a
subalgebra $\N(S)\subset (\prod_{n,\U} (\Ma_{2^n}\ten
N^{\ten_n})_*)^*$ and a map $u_{-1}:S\to N_\U$ such that
 \[ \si_t^{\phi_\U}(u_{-1}(x))\lel u_{-1}(\si_t^{\psi}(x))\]
holds for all $t\in \rz$ and $x\in S$. In particular, we find a
$\phi_\U$-preserving conditional expectation $E:N_\U\to \N(S)$.
This yields a completely isometric embedding $\iota:L_1(\N(S))\to
L_1(\N_\U)$. According to \cite{JF},  Lemma \ref{ne1} and Lemma
\ref{limitn} we deduce for $y=\sum_l f_l\ten x_l$ that
 \begin{align*}
 &\noo
 D_{\phi_\U\ten \phi_{\Mz}}^{\frac12}(u_{-1}\ten id)(y)D_{\phi_\U\ten \phi_{\Mz}}^{\frac12}
 \rrm_{L_1(\N(S))\wet L_1(\Mz)}
 \lel \noo
 D_{\phi_\U\ten \phi_{\Mz}}^{\frac12}(u_{-1}\ten id)(y)D_{\phi_\U\ten \phi_{\Mz}}^{\frac12}
 \rrm_{L_1(N_\U)\wet L_1(\Mz)} \\
 &=   \lim_{n,\U} \|
 D_n^{\frac12}(u_{n,T}\ten id) (y)D_n^{\frac12}\| \\
 &\sim_{200} \inf_{y=y_2+y_3}
 \noo (D_{\psi}^{\frac12}\ten D_{\phi_{\Mz}}^{\frac12})y_2(1\ten
 D_{\phi_{\Mz}}^{\frac12})\rrm_{L_2^r(N,tr)\wet L_1(\Mz)}
   + \noo (1\ten D_{\phi_{\Mz}}^{\frac12})y_3(D_{\psi}^{\frac12}\ten D_{\phi_{\Mz}}^{\frac12})\rrm_{L_2^c(N,tr)\wet
  L_1(\Mz)}\!\!\!  \pl .
  \end{align*}
Now, we consider the subspace $K_{12}={\rm span}\{ \delta_k\ten
e_{12}:k=1,...,m\}\subset L_2(N,tr)$. It is very easily shown that
the orthogonal projection
 \[ P_{12}\kla \begin{array}{cc} f_{11}&f_{12}\\f_{21}&
 f_{22}\end{array}\mer \lel \kla \begin{array}{cc} 0 &f_{12}\\0&
 0 \end{array}\mer \]
satisfies
 \[ \noo P_{12}\ten id_{L_1(\Mz)}:L_2^r(N,tr)\wet L_1(\Mz) \to L_2^r(N,tr)\wet L_1(\Mz) \rrm
 \kl 1 \pl \]
and
 \[ \noo P_{12}\ten id_{L_1(\Mz)}:L_2^c(N,tr)\wet L_1(\Mz) \to L_2^c(N,tr)\wet L_1(\Mz) \rrm
 \kl 1 \pl .\]
Thus given $y=\sum_k \delta_k \ten e_{12}\ten y_k$ and any
decomposition $y=\tilde{y}_2+\tilde{y}_3$, we may define
 \[ y_2\lel (P\ten id)(\tilde{y}_2)\lel \sum_k \delta_k \ten e_{12}\ten
 v_k \quad \mbox{and} \quad
 y_3\lel (P\ten id)(\tilde{y}_3)\lel \sum_k \delta_k \ten e_{12}\ten
 w_k\pl .\]
Then we have
 \begin{align*}
 &\|(D_{\psi}^{\frac12}\ten D_{\phi_{\Mz}}^{\frac12})y_2(1\ten
 D_{\phi_{\Mz}}^{\frac12})\|_{L_2^r(N,tr)\wet L_1(\Mz)}
 \lel \|(P_{12}\ten id)(D_{\psi}^{\frac12}\ten D_{\phi_{\Mz}}^{\frac12})\tilde{y}_2(1\ten
 D_{\phi_{\Mz}}^{\frac12})\|_{L_2^r(N,tr)\wet L_1(\Mz)}\\
 & \kl \|(D_{\psi}^{\frac12} \ten D_{\phi_{\Mz}}^{\frac12})\tilde{y}_2(1\ten
 D_{\phi_{\Mz}}^{\frac12})\|_{L_2^r(N,tr)\wet L_1(\Mz)} \pl .
 \end{align*}
The same argument works for $y_3$. Hence it suffices to consider
decompositions of the form $y_2$ and $y_3$ above. We note that
 \[ (D_{\psi}^{\frac12}(\delta_k\ten e_{12}),D_{\psi}^{\frac12}(\delta_k\ten
 e_{12}))_{tr}
 \lel tr((\delta_l \ten e_{21})D_{\psi}(\delta_k\ten e_{12}))
 \lel \delta_{k,l} (1-\mu_k) \pl .\]
Thus
 \begin{align*}
 \|(D_{\psi}^{\frac12} \ten D_{\phi_{\Mz}}^{\frac12})y_2(1\ten
 D_{\phi_{\Mz}}^{\frac12})\|_{L_2^r(N,tr)\wet L_1(\Mz)}
 &\lel \| (\sum_k (1-\mu_k)
 (D_{\phi_{\Mz}}^{\frac12}v_kD_{\phi_{\Mz}}^{\frac12})^*(D_{\phi_{\Mz}}^{\frac12}v_kD_{\phi_{\Mz}}^{\frac12}))^{\frac12}\|_{L_1(\Mz)}
 \pl.
 \end{align*}
For the $y_3$-term we find
  \[ ((\delta_l\ten e_{12})D_{\psi}^{\frac12},(\delta_k\ten e_{12})D_{\psi}^{\frac12})
 \lel tr(D_{\psi}(\delta_k\delta_l\ten e_{21}e_{12}))\lel
 \delta_{kl} \mu_k \pl .\]
Thus by approximation of $x_k$ with analytic elements of the form
$D_{\phi_{\Mz}}^{1/2}y_kD_{\phi_{\Mz}}^{1/2}$, we obtain
 \begin{align}
 &\noo \sum_k
 D_{\phi_\U}^{\frac12}(u_{-1}(\delta_k\ten
 e_{12}))D_{\phi_\U}^{\frac12} \ten x_k
 \rrm_{L_1(\N(S))\wet L_1(\Mz)} \nonumber\\
 &\quad \quad \sim_{200} \inf_{x_k=c_k+d_k} \|(\summ_k
 (1-\mu_k) c_k^*c_k)^{\frac12} \|+\|(\summ_k
 \mu_k d_k^*d_k)^{\frac12} \| \pl . \label{thatsit}
 \end{align}
Now, we have to identify $u_{-1}(\delta_k \ten e_{12})$  with the
standard generators satisfying the CAR relations. Indeed, we
recall from Theorem \ref{uniq} that we may assume that the map
$\al:A(S)\to L_2(\N(S))$ has  dense range. This enables us  to
apply Lemma \ref{carrel} and to conclude that for
$b_k=\delta_k\ten e_{12}$ we have
 \[ u_{-1}(b_k)u_{-1}(b_j)+u_{-1}(b_j)u_{-1}(b_k)\lel 0 \quad ,\quad
 u_{-1}(b_k)u_{-1}(b_j)^*+u_{-1}(b_j)^*u_{-1}(b_k) \lel \delta_{kj} \]
and
 \[ \phi_\U(u_{-1}(b_{i_r}^*)\cdots u_{-1}(b_{i_1}^*)u_{-1}(b_{j_1})\cdots
 u_{-1}(b_{j_s}))
 \lel  \delta_{r,s} \prod_{l=1}^r \delta_{i_l,j_l}\mu_{i_l} \]
holds for all increasing sequences  $i_1<i_2<\cdots <i_r$ and
$j_1<i_2<\cdots <j_s$. Therefore $\N(S)$ is indeed an isomorphic
copy of the CAR algebra $\Mz_{2^m}$ and $\phi_\U$ induces the
usual quasi-free state $\phi_\mu$. Thus we may replace
$u_{-1}(b_k)$ by $a_k$ in \eqref{thatsit} and the proof is
completed. \qd

\begin{rem}\label{discq} {\rm Only minor modifications of this proof provide the
Khintchine inequality for $q$-gaussian variables in the discrete
setting. We use $N=\ell_{\infty}^m(\Mz_2)$ with the weight
$\psi_{\mu}$ and the inclusion map $j:\ell_{\infty}^m\to
L_2(N,\psi)$. We have seen in Corollary \ref{identify} that the
von Neumann algebra $\Gamma_q(K,U_t)$ and the von Neumann algebra
$\N(S)$ generated by $\{u_q(j(f)):f\in \ell_{\infty}^m\}$ are
isomorphic and $\pi(s_q(f))=u_q(f)$. By construction,
$\N(S)\subset (\prod_n (L_{\infty}(\Om_n;\Mz_n)\ten N^{\ten
n})_*)^*$. Moreover, we have a $\phi_\U$-preserving conditional
expectation and hence a natural inclusion
 \[ L_1(\N(S))\subset \prodd_{n,\U} L_1(\Om_n;\Mz_n\ten N^{\ten n})
 \pl . \]
We deduce from Corollary \ref{identify} and analyticity that
 \begin{align*}
 &\|\summ_l D_{vac}^{\frac12}s_q(f_l)D_{vac}^{\frac12}\ten
 x_l\|_{L_1(\N(S))\wet L_1(\Mz)} \lel  \|\summ_l
 D_{\phi_\U}^{\frac12}u_q(j(f_l))D_{\phi_\U}^{\frac12}\ten
 x_l\|_{L_1(\N(S))\wet L_1(\Mz)}  \\
 &=  \lim_{n,\U} \|\summ_l
 D_n^{\frac12}u_{n,T}(f_l)D_n^{\frac12}\ten
 x_l\|_{L_1(\Om_n;L_1(\Mz_n\ten N^{\ten_n}))\wet L_1(\Mz)}  \pl
 \end{align*}
holds for all $e_j\in \ell_{\infty}^m$ and  $x_l\in L_1(\Mz)$.
Here $u_{n,T}=\sqrt{\frac{T}{n}}\sum_{j=1}^n v_{j}(n)\ten
\pi_j(x)$ and the $v_j(n)$ are  Speicher's  random matrices (see
Corollary \ref{prob}). We consider the special case $f_{2l}=
\delta_l$ and $f_{2l+1}=-i\delta_l$. The $q$-gaussian analogue of
the elements $a_k$ are then given by
 \[ a_k(q) \lel \frac{1}{2}[s_q(\delta_k)-is_q(i\delta_k)]  \cong
  \frac{1}{2}(u_q(j(\delta_k))-iu_q(j(i\delta_k))) \lel
 u_q(\delta_k \ten e_{12})\pl. \]
Therefore the argument in the proof of Theorem \ref{carkh} shows
that
 \begin{align*}
  &\noo \summ_k
  D_{vac}^{\frac12}a_k(q)D_{vac}^{\frac12}\ten x_k\rrm_{L_1(N_{\mu})\wet
  L_1(\Mz)}
\!\!  \!\! \! \sim \!
   \inf
   \noo (\summ_k (1-\mu_k)c_k^*c_k)^{\frac12}\rrm_{L_1(\Mz)}
   + \noo (\summ_k \mu_k d_kd_k^*)^{\frac12}\rrm_{L_1(\Mz)}
   \end{align*}
where the infimum is taken over $x_k=c_k+d_k$.  For $q=1$ we shall
understand $s_1(f)D_{vac}^{1/2}$ as $\frac{d}{idt}
W(\sqrt{2}f)D_{vac}^{1/2}|_{t=0}$. Here $\{W(f):f\in
\ell_{\infty}^m(\cz)\}$ is a family of  unitaries in the Segal
representation (see the proof of Corollary \ref{identify}). We may
then extend $s_1$ by linearity.}\end{rem}

\begin{rem}{\rm In the continuous case we may assume that
$K=L_2(\mu;\cz)$, the inclusion map is given by  $j:K\to
L_2(N,\psi)$, $j(f)=f\ten e_{12}+\bar{f}\ten e_{21}$ and the
weight is given by $\psi(x)=\int_{\Om} [f_1x_{11}+f_2x_{22}]d\mu$.
As usual we assume $f_1+f_2=1$. For $f\in L_2(\mu;\rz)$ we shall
now define
 \[ a_q(f)\lel \frac12(s_q(j(f))-is_q(j(if))) \pl .\]
Then we have %\label{cont-kh}
 \begin{align}
 & \|\sum_l D_{vac}^{\frac12}a_q(f_l)D_{vac}^{\frac12} \ten
  x_l\|_{L_1(\Gamma_q(K,U_t))\wet L_1(\Mz)}
  \sim_{200} \inf_{\sum_l f_l\ten x_l=c+d}  \nonumber \\
 &\pll  \|(\int c(\om)^*c(\om) f_1(\om) d\mu(\om))^{\frac12}\|_{L_1(\Mz)}
 + \|(\int d(\om)d(\om)^* f_2(\om)
 d\mu(\om))^{\frac12}\|_{L_1(\Mz)}
 \pl . \label{cont-kh}
 \end{align}
Let us first assume that $f_2=\sum_k f_2(k)1_{A_k}$ is a finite
simple function. Then \eqref{cont-kh} follows by approximation of
the $f_l$'s by simple functions $f_l \in L_2(\Om_m,\Si,\mu)$ where
$\Si$ is a finite $\si$-algebra on $\Om_m=A_1\cup \cdots \cup A_m$
generated by by $A_1$,...,$A_m$. In that case the unitary group
$U_t=(f_2/f_1)^{it}$ leaves the subspace $L_2(\Om_m,\Si,\mu;\cz)$
invariant. This in turn implies that we have a conditional
expectation $E:\Gamma_q(K,U_t)\to
\Gamma_q(L_2(\Om,\Si,\mu;\cz),U_t)$ and therefore the norm
estimates in $L_1$ are preserved. By density we then obtain
\eqref{cont-kh} for infinite simple functions. For general $f_2$
we consider the sequence
 \[ f_2^r \lel \summ_{k+1 \le -\frac{r}{2}}  2^{\frac{k}{r}} 1_{\{2^{\frac kr}\le f_2<
 2^{\frac{k+1}{r}}\}} +  \summ_{k+1 \le -\frac{r}{2}}  (1-2^{\frac{k}{r}}) 1_{\{2^{\frac kr}\le 1-f_2\le
 2^{\frac{k+1}{r}}\}} \pl
  \]
which converges to $f_2$ everywhere. Then we may consider the
density $\psi_r$ given by $(1-f_2^r,f_2^r)$ and the vacuum state
$\om_r$ on $\Gamma_q(K,U_t(r))$. We use  the ultraproduct
 \[ M_\U\subset e_\U(\prodd_{r,\U} \Gamma_q(K,U_t(r))_*)^*e_\U \]
formed with respect to the support $e_\U$ of the ultraproduct
state $(\om_r)^{\bullet}$. We define a generating system $S\subset
L_2(\mu;\cz)$ of all bounded functions with support contained in
one of the  sets $\Om_m=\{\om : 2^{-m}\le f_2(\om)\le 1-2^{-m}\}$.
The advantage of this generating set is that for $f\in S$ the
family $\si_t^{\psi_r}(f)$ is uniformly bounded and
$\si_t^{\psi_r}(S)\subset S$ for all $r$. Moreover, the family
$(s_{q,r}(f))_{r\in \nz}$ satisfies the moment conditions in
Remark \ref{amom} and therefore we find operators $s_{q,\U}(f)$
affiliated to $M_\U$ such that
 \begin{align*}
 &(D_\U^{\frac12},s_{q,\U}(f_1)\cdots
 s_{q,\U}(f_m)D_{\U}^{\frac12})
 \lel \lim_{r,\U} (D_{\om_r}^{\frac12},s_{q,r}(f_1)\cdots
 s_{q,r}(f_m)D_{\om_r}^{\frac12}) \\
 &= \lim_{r,\U}
    \summ_{\si=\{\{i_1,j_1\},...,\{i_{\frac{m}{2}},j_{\frac{m}{2}}\}\}
 \in P_2(m)}
 q^{I(\si)} \prod_{l=1}^{\frac{m}{2}} \psi_r(j(f_{i_l})j(f_{j_l}))
 \\
 &=   \summ_{\si=\{\{i_1,j_1\},...,\{i_{\frac{m}{2}},j_{\frac{m}{2}}\}\}
 \in P_2(m)}
 q^{I(\si)} \prod_{l=1}^{\frac{m}{2}} \psi(j(f_{i_l})j(f_{j_l}))
 \pl .
 \end{align*}
According to Theorem \ref{uniq} the subalgebra $M(S)\subset M_\U$
generated by the spectral projections of the elements
$\{s_{q,\U}(f):f\in S\}$ is isomorphic to $\Gamma_q(S,U_t)$. Our
choice of $S$ also guarantees that $M(S)$ is invariant under the
modular group of the ultraproduct state. Therefore, we have a
completely isometric embedding of $L_1(\Gamma_q(S,U_t))\subset
L_1(M_\U)$. However, for every $r$ the estimate \eqref{cont-kh}
holds up to a factor $c(r)$ with $\lim_r c(r)=1$. Thus the
estimates also hold in the limit and hence for $\Gamma_q(S,U_t)$.
By density, we may then extend it to $\Gamma_q(K,U_t)$.
 }\end{rem}

\begin{rem}{\rm
Let us note that for fixed $0<\la<1$ we may find a function
$f_2^{\la}$ with $1/c(\la)f_1\le f_1^{\la}\le c(\la)f_1$,
$1/c(\la)f_2\le f_2^{\la}\le c(\la)f_2$ such that
$f_1^{\la}/f_2^{\la}\in \{\la^{n}: n \in \zz\}$. Thus we obtain
the same estimates in a factor of type III$_{\la}$, see Theorem
\ref{hia} for $-1<q<1$. For the border case $q=\pm 1$ we may
perform a similar construction using Remark \ref{typecar} or
Theorem \ref{typccr}.}
\end{rem}

\section{Applications to operator spaces}

In this section we assume the reader to be familiar with operator
spaces theory.  We are now well-prepared for the proof of
Corollary \ref{emb0}.

\begin{theorem}\label{emb} Let $Q$ be a quotient of  $R\oplus C$. Then $Q$
embeds into the predual of the  hyperfinite {\rm III}$_{\la}$,
$0<\la\le 1$.\end{theorem}

\begin{proof} We will start with a general characterization of
quotients of $R\oplus C$. Such a quotient has a direct
decomposition
 \[ Q \approx  R_n\oplus C_m \oplus (R\oplus C/\gr(D_\la))\]
where $n,m\in \nz_0\cup\{\infty\}$ and
$\gr(D_\la)=\{(e_{k,1},\la_k e_{k,1}): k\in \nz\}$  is the graph
of a diagonal operator. Here $\approx$ stand for a completely
isomorphism $u$ with $\|u\|_{cb}\|u^{-1}\|_{cb}\le 4$. Note that
the sequence $(\la_k)$ might be finite as well. Let us briefly
sketch the proof of this decomposition. Indeed, we consider the
linear subspace $S=Q^{\perp}\subset C\oplus R$ and the projections
$\pi_C$ and $\pi_R$ on the column and row component. By splitting
off $R_m$ and $C_m$, we may assume that $\pi_C(S^{\perp})$ and
$\pi_R(S^{\perp})$ are dense. Then, we mod out by $S\cap
\{0\}\times R$ and $S\cap C\times \{0\}$. By homogeneity this does
not change the operator space structure. We obtain
$\tilde{S}\subset H^r\oplus K^c$ such that $\tilde{S}$ is the
graph of an operator with domain in  a separable Hilbert space
$H$. Using homogeneity again we may assume that $\tilde{S}\subset
H^r\oplus H^c$ is the graph of a positive operator $T$ on $H^r$.
By the spectral theorem (see e.g. \cite[Theorem 5.6.2]{Kad}) $T$
is unitarily equivalent to a multiplication operator $M_f(f)=fg$
on some $L_2(\Om,\mu)$. Using a small perturbation, we may assume
that $f$ is an infinite sum of characteristic functions (changing
the operator space structure of $L_2^r(\mu)\oplus L_2^c(\mu)/{\rm
graph}(M_f)$ only by a constant $(1+\eps)$). By assumption
$L_2(\mu)$ is separable and hence a suitable choice of the basis
yields a diagonal operator. For more details on this argument,
known to Xu and the author for quite some time, see also
\cite{Pl}.

Thus in the following, we have to find an embedding of
$Q(\la)=R\oplus C/{\rm graph}(D_{\la})$ such that $(\la_k)$ are
positive integers. We observe that $Q(\la)$ has a basis $(f_k)$
such that the quotient mapping $q:R\oplus C\to Q(\la)$ is given
by $q((e_{1k},0))=-\la_kf_k$ and $q((0,e_{k1}))=f_k$. This implies
that
 \begin{align*}
  \noo \summ_k x_k \ten f_k\rrm_{Q(\la)\wet L_1(\Mz_m)}
  &= \inf_{x_k=-\la_k c_k+d_k} \noo (\summ_k
  c_k^*c_k)^{\frac12}\rrm_{L_1(\Mz_m)}+
   \noo (\summ_k d_kd_k^*)^{\frac12}\rrm_{L_1(\Mz_m)} \!\!\!.
  \end{align*}
Let $\mu_k$ to be determined later. The change of variables
$\hat{c}_k=-(1-\mu_k)^{-1/2}c_k$ and  $\hat{d}_k=\mu_k^{-1/2}d_k$
yields
 \begin{align*}
  \noo \summ_k x_k \ten f_k\rrm_{Q(\la)\wet L_1(\Mz_m)}
  &= \inf
  \noo (\summ_k (1-\mu_k) \hat{c}_k^*\hat{c}_k)^{\frac12}\rrm_{L_1(\Mz_m)}
 + \noo (\summ_k \mu_k
 \hat{d}_k\hat{d}_k^*)^{\frac12}\rrm_{L_1(\Mz_m)} \pl .
  \end{align*}
Here the infimum is taken over $\mu_k^{-1/2}x_k=-\mu_k^{-1/2}\la_k
c_k+\mu_k^{-1/2}c_k=\hat{c}_k+\hat{d}_k$. For the last equality to
hold we choose $\mu_k=(1+\la_k^{-2})^{-1}$.
 Then Theorem
\ref{carkh} implies that
$w(f_k)=\mu_k^{-1/2}D_{\mu}^{1/2}a_kD_{\mu}^{1/2}$ extends to a
complete isomorphism $w:Q(\la)\to L_1(\N(\mu))=\N(\mu)^{op}_*$.
However, $\N(\mu)$ and $\N(\mu)^{op}$ are hyperfinite von Neumann
algebras. If we want to accommodate the additional pieces $R_m$
and $C_m$, we can use $N=B(\ell_2)\oplus \N(\mu)^{op}\subset
B(\ell_2)\ten \N(\mu)^{op} $. Using a conditional expectation, we
can replace $\mu$ by a sequence $\mu'$ such that for every
rational $0<\la<1$ there are infinitely may $\mu'_k$'s with
$\mu_k=\la/1+\la$. According to \cite[section8]{AW}, we deduce
that $\N(\mu')$ and $\N(\mu')^{op}$ are  of type III$_1$. Then
$B(\ell_2)\ten \N(\mu')^{op}\cong \N(\mu')^{op}$ and we find a
complete  embedding in the predual of a hyperfinite III$_1$
factor. Due to the results of Haagerup, Rosenthal and Sukochev
\cite{HRS}, the predual of the hyperfinite III$_1$ factor also
embeds into the predual of the hyperfinite III$_{\la}$ factor for
all $0<\la<1$. \qd

\begin{rem}{\rm We  learned from Haagerup that a von Neumann
algebra admits a normal conditional expectation onto a copy of the
hyperfinite III$_1$ factor if and only if the flow of weights,
i.e. the restriction of the dual automorphism group on the center
of the core, admits a normal invariant measure. In particular, for
those von Neumann algebras $N$ we have an embedding of the predual
of the hyperfinite factor in the predual of $N$. This implies that
$R\oplus C/Q$ completely embeds in $N_*$. It is open whether every
quotient of   $R\oplus C$ completely embeds into the predual of
every type III factor. }\end{rem}

\begin{rem}\label{concr} {\rm In our applications, we will often consider a concrete
quotient $Q=Q(\la,\nu)$  of $R\oplus C$ with basis $(f_k)$
satisfying
 \begin{align*}
 \noo \summ_k f_k\ten x_k\rrm_{Q(\la,\nu)\wet L_1(\Mz_m)} &=
  \inf_{x_k=c_k+d_k} \|(\sum_k \la_k c_k^*c_k)^{\frac12}\|_1+
 \|(\sum_k \nu_k d_kd_k^*)^{\frac12}\|_1 \\
 &= \inf_{\sqrt{\la_k+\nu_k}\p x_k=c_k+d_k} \|(\sum_k \frac{\la_k}{\la_k+\nu_k} c_k^*c_k)^{\frac12}\|_1+
 \|(\sum_k \frac{\nu_k}{\la_k+\nu_k}d_kd_k^*)^{\frac12}\|_1 \\
 &\sim_{200} \|\summ_k \sqrt{\la_k+\nu_k} \p D_{\mu}^{\frac12}a_kD_{\mu}^{\frac12}
 \ten  x_k \|_{L_1(\N(\mu))\wet L_1(\Mz_m)}
 \end{align*}
for all $x_k\in L_1(\Mz)$. Here we used
$\mu_k=\frac{\nu_k}{\la_k+\nu_k}$ and
$\phi_{\mu_k}(x)=(1-\mu_k)x_{11}+\mu_kx_{22}$. Then $D_{\mu}$ is
the density of the quasi-free state $\phi_{\mu}$ on $\N(\mu)$. We
obtain an isomorphism $v^{op}:Q\to L_1(\N(\mu))$ defined by
 \[ v^{op}(f_k)\lel \sqrt{\la_k+\nu_k}
 D_{\mu}^{\frac12}a_kD_{\mu}^{\frac12}\pl .\]
Let us  make this even more explicit. We may use the transposition
map $\pi^t(x_1\ten \cdots x_m\ten 1\cdots)= x_1^t \ten \cdots \ten
x_m^{t}\ten 1\cdots$ and obtain an isomorphism $\pi^t:\ten_k
\Mz_2\to \ten_k \Mz_2^{op}$ such that $\phi_{\mu}\pi^t\lel
\phi_{\mu}$. Therefore $\pi^t$ extends to an isomorphism
$\pi^t:\N(\mu)\to \N(\mu)^{op}$. Note that $a_k^t=a_k^*$  and
hence
 \begin{align*}
 &v^{op}(f_k)(\pi^t(x))  \lel  (\la_k+\nu_k)^{\frac12} \pl
tr(D_{\mu}^{\frac12}a_kD_{\mu}^{\frac12}\pi^t(x))
 \lel [(\la_k+\nu_k)\mu_k(1-\mu_k)^{-1}]^{\frac12}
\pl  tr(a_kD_{\mu}\pi^t(x)) \\
 &= [(\la_k+\nu_k)\mu_k(1-\mu_k)^{-1}]^{\frac12} \pl
 \phi_{\mu}(\pi^t(x)a_k) \lel [(\la_k+\nu_k)\mu_k(1-\mu_k)^{-1}
]^{\frac12}
 \pl  \phi_{\mu}(a_k^tx) \\
 &= [(\la_k+\nu_k)\mu_k(1-\mu_k)^{-1}]^{\frac12} \pl
 \phi_{\mu}.a_k^*(x)  \lel
 [(\la_k+\nu_k)(1-\mu_k)\mu_k^{-1}]^{\frac12} \pl
 a_k^*.\phi_{\mu}(x)
 \pl .
 \end{align*}
This implies that $v^{op}\pi^t(f_k)=
[(\la_k+\nu_k)\nu_k/\la_k]^{1/2}
 \phi_{\mu}.a_k^*=[(\la_k+\nu_k)\la_k/\nu_k]^{1/2}
 a_k^*.\phi_{\mu}$ is a complete isomorphism in the predual of $\N(\mu)$.
 }\end{rem}

\begin{rem}{\rm We refer to \cite{Pl} for a characterization of
operator spaces in $Q(R\oplus C)$ which embed in the predual of a
semifinite hyperfinite factor.}
\end{rem}

Before studying further applications to operator spaces, we
provide an application of operator space theory to  norm
inequalities for linear functionals.

\begin{cor} Let $\phi_{\mu}$ and $\phi_{\nu}$ be quasi free states
on the CAR algebra. Let $(b_{ij})$ be a matrix with finitely many
non-zero entries.  Then
\begin{align*}
 &\noo \summ_{ij} b_{ij} (a_i.\phi_{\mu}\ten
 a_j.\phi_{\nu})\rrm_{(\N(\mu)\bar{\ten}\N(\nu))_*}
 \sim_{c} \inf_{b_{ij}=f_{ij}+g_{ij}+h_{ij}+k_{ij}}
 \kla \summ_{ij}  |f_{ij}|^2 \mu_i\nu_j \mer^{\frac12}\\
 &\pl  + \noo \left[g_{ij}\frac{\sqrt{\mu_i}\nu_j}{\sqrt{1-\nu_j}}\right]
 \rrm_{S_1}+
 \noo \left [h_{ij} \frac{\mu_i \sqrt{\nu_j}}{
 \sqrt{1-\mu_i}}
 \right] \rrm_{S_1}
 +\kla \summ_{ij} |k_{ij}|^2
 \frac{\mu_i^2\nu_j^2}{(1-\mu_i)(1-\nu_j)}
 \mer^{\frac12}
 \pl .
 \end{align*}\end{cor}

\begin{proof} From \eqref{erdual} and \eqref{wetL1} we deduce that
 \begin{align*}
 &\noo \summ_{ij} b_{ij} a_i.\phi_{\mu}\ten
 a_j.\phi_{\nu}\rrm_{(\N(\mu)\bar{\ten}\N(\nu))_*}
 \lel  \noo \summ_{ij} b_{ij} a_iD_\mu \ten
 a_jD_{\nu}\rrm_{L_1(\N(\mu)\bar{\ten}\N(\nu))} \\
 &\pl = \noo \summ_{ij} b_{ij}(\frac{\mu_i}{1-\mu_i}\frac{\nu_j}{1-\nu_j})^{\frac12}
  D_{\mu}^{\frac12}a_iD_\mu^{\frac12} \ten
 D_{\nu}^{\frac12}a_jD_{\nu}^{\frac12} \rrm_{L_1(\N(\mu))\wet L_1(\N(\nu))}
 \pl .
 \end{align*}
Now we consider the operator spaces
 \[ X(\mu) \lel \overline{{\rm
 span}}\{D_{\mu}^{1/2}a_kD_\mu^{1/2}:k\in \nz\} \pll \mbox{and}\pll
  X(\nu)\lel \overline{{\rm span}}\{D_{\nu}^{1/2}a_kD_\nu^{1/2}:k\in \nz\} \pl .\]
with canonical basis $f_k(\mu)=D_{\mu}^{1/2}a_kD_\mu^{1/2}$ and
$f_k(\nu)=D_{\nu}^{1/2}a_kD_\nu^{1/2}$. According to Theorem
\ref{carkh} we have
 \[ X(\mu) \approx_{200} \ell_2^r(1-\mu)\oplus
 \ell_2^c(\mu)/\Delta \pl =:\pl  K(\mu)   \pl ,\]
where $\Delta=\{(x,x):x\in \ell_2(1-\mu)\cap \ell_2(\mu)\}$ is the
diagonal, and the isomorphism is given by
$u(f_k)=(e_k,e_k)+\Delta$. Let us denote the operator space
defined by the right hand side by $K(\mu)$. In \cite[Corollary
7.12]{Joh} it is shown that $K(\mu)$ is completely contractively
complemented in $L_1(M(\mu))$ for some von Neumann algebra with
QWEP.  We denote the corresponding embeddings by
$w_{\mu}:K(\mu)\to L_1(M(\mu)))$  and $w_{\nu}:K(\mu)\to
L_1(M(\nu)))$, respectively. According to \cite[Lemma 4.4, Lemma
4.5]{Joh} and \cite[Corollary 7.12]{Joh}, we deduce that
 \begin{align*}
 &\|\summ_{ij} c_{ij} f_i(\mu)\ten f_j(\nu)\|_1
  \lel   \pi_1^o(T_c:X(\mu)^*\to X(\nu)) \sim_{200^2}
 \pi_1(T_c:K(\mu)^*\to K(\nu))\\
 & \sim_{9}
 \|\summ_{ij} c_{ij} w_{\mu}u(f_i)\ten w_{\nu}u(f_j)\|_{L_1(M(\mu))\wet
 L_1(M(\nu))} \sim_{9} \|\summ_{ij} c_{ij} u(f_i)\ten u(f_j)\|_{K(\mu)\wet
 K(\nu)} \pl .
 \end{align*}
We may now apply \cite[Lemma 5.1]{Joh} and deduce that
 \begin{align*}
 &K(\mu) \wet K(\nu) \lel  \ell_2^r(1-\mu)\oplus \ell_2^c(\mu)/\Delta \wet   \ell_2^r(1-\nu)\oplus
  \ell_2^c(\nu)/\Delta \\
  &= \ell_2(\nz^2;(1-\mu)\ten (1-\nu)) \oplus
  \ell_2(1-\mu)\ten_{\pi} \ell_2(\nu) \oplus
  \ell_2(\mu)\ten_{\pi} \ell_2(1-\nu)\oplus \ell_2(\nz^2;\mu \ten
  \nu)/\Delta\ten \Delta \pl .
 \end{align*}
Here $H\ten_{\pi}K=S_1(H,K)$ is the Banach space projective tensor
product and can be calculated using the norm in the Schatten
class. The assertion follows from a change of variables. \qd

\begin{rem}{\rm Following Remark \ref{discq}, we know that $X(\mu)$ and
\[ X_{free}(\mu)=\overline{{\rm span}}\{D_{\mu}^{1/2}a_k(0)D_\mu^{1/2}:k\in \nz\}\]
are completely isomorphic because they both satisfy the formula in
Theorem \ref{carkh}. In the case of free random variables this
result follows by easily duality from Pisier's  estimates in
\cite{Pl}. Pisier also proves the complementation result in
$\N_{free}(\mu)$. This yields slightly better constants than using
the results from \cite{Joh}. }\end{rem}

\begin{prob} Describe the operator space structure of $\overline{\rm
span}\{a_ia_j.\phi_{\mu}: i,j\in  \nz\}$.
\end{prob}

\begin{cor}\label{caracG} Let $X$ be a separable operator space. Then $X$ and $X^*$
embeds into the predual of a hyperfinite factor if and only if
there exists subspaces $S \subset R\oplus C$  such that $X$ is
completely isomorphic to a subspace of $R\oplus C/S$.
\end{cor}

\begin{proof} Let $X\subset Q=R\oplus C/S$ be a subspace of a
quotient of $R\oplus C$. Then Theorem \ref{emb} implies that $Q$
and hence $X$ embeds in the predual of the hyperfinite III$_1$
factor. However, $(R\oplus C)^*=C\oplus R$ is completely
isomorphic to $R\oplus C$. The duality between subspaces and
quotients implies $X^*\in QS((R\oplus C)^*)\lel QS(C\oplus R)$. It
is an elementary fact in Banach space theory that $QS(X)=SQ(X)$.
This also holds true in the category  of operator spaces. Thus
$X\in SQ(R\oplus C)$ implies $X^*\in SQ(R\oplus C)$ and hence both
embed in the predual of a hyperfinite III$_1$ factor. Since the
predual of a hyperfinite factor has the completely bounded
approximation property, the converse follows from
Pisier/Shlyahtenko's Grothendieck Theorem \cite[Theorem
0.5ii)]{PS}.
 \qd

The following result motivated our approach.
\begin{cor} The operator space $OH$ embeds into the hyperfinite
{\rm III}$_\la$ factor, $0<\la\le 1$.
\end{cor}

\begin{proof} In \cite{Joh}, we found  the  densities $f_1=1/t$ and $f_2=1/1-t$
with respect to $\mu \lel dt/\pi\sqrt{t(1-t)}\ten m$, $m$ the
counting measure. In this case the spectrum $f_2/f_1$ is
continuous. Both embeddings using the CAR and CCR relations lead
to the hyperfinite III$_1$ factor. The results of of Haagerup,
Rosenthal and Sukochev \cite{HRS} complete the proof. \qd

We want to discuss two concrete embeddings for the interpolation
spaces $R_p=[R,C]_{\frac1p}$, see \cite{Po} for a precise
definition.  We will need the results from \cite{JX3}. The space
$R_p$ has a basis $(g_k)$ such that
 \begin{align}\label{jx1}
 \noo \summ_k g_k\ten x_k\rrm_{R_p\wet L_1(\Mz_m)} \sim_{c_p} \inf_{x_k=c_{kj}+d_{kj}}
 &\pll \noo (\summ_{kj} (1-\si_j)c_{jk}^*c_{jk})^{\frac12}\rrm_1+
  \noo (\summ_{kj} \si_jd_{jk}d_{jk}^*)^{\frac12}\rrm_1 \pl.
  \end{align}
Here $j \in \zz$ and the coefficients satisfy
 \begin{equation}\label{coef}  \si_j \lel \begin{cases} |j|^{-p'} & j\ge 1\\
                             \frac12      & j=0\\
                             1-|j|^{-p} & j\le 1
                             \end{cases}
                             \pl . \end{equation}

\begin{cor}\label{concr5}  Let $(\mu_{kj})$ be a double indexed sequence
such that $\mu_{kj}=\si_j$.  Let $(a_{kj})$ be a  double indexed
sequence of generators of the CAR algebra satisfying
\eqref{carrk}. Then the map
 \[ v(g_k) \lel \summ_{j<0}
 (1+|j|)^{-\frac{p}{2}}a_{kj}^*.\phi_{\mu} +
 \summ_{j\ge 0}  (1+|j|)^{-\frac{p'}{2}} \phi_\mu.a_{kj}^*.
 \]
defines a complete embedding of $R_p$ in $\N(\mu)_*$.
\end{cor}

\begin{proof} This  follows immediately
from \eqref{jx1} and Remark \ref{concr}. In order to obtain
absolutely summable coefficients, we observe  that for $j\ge 0$ we
have $[\si_j(1-\si_j)^{-1}]^{1/2}\sim_c (1+|j|)^{-p'/2}$. In that
case we use $\phi_{\mu}.a_{kj}^*$. For $j<0$ we prefer
$a_{kj}^*.\phi_{\mu}$ and find $[(1-\si_j)\si_j^{-1}]^{1/2}\sim_c
(1+|j|)^{-p/2}$.\qd

In our next application we want to find an embedding of $R_p^n$ in
$S_1^m$ with some control of $m=m(p,n)$. We refer to \cite{JX} for
the following result. Let $\la>1$. Then there exists a constant
$c(p,\la)$ such that
\begin{align}\label{concrt1}
 \noo \summ_k g_k\ten x_k\rrm_{R_p\wet L_1(\Mz_m)}\!\!\! \sim_{c(p,\la)} \inf_{x_k=c_{kj}+d_{kj}}
 &\pll \noo (\summ_{k\in \nz, j\in \zz} \la^{\frac{j}{p'}} c_{jk}^*c_{jk})^{\frac12}\rrm_1+
  \noo (\summ_{k\in \nz,j\in\zz}  \la^{-\frac{j}{p}}
  d_{jk}d_{jk}^*)^{\frac12}\rrm_1\pl .
  \end{align}

\begin{lemma}\label{reduc} There  exists a constant $c(p,\la)$ and  $c(p)$ such that for all $\nen$
\begin{samepage} \begin{align*}
 &\noo \summ_{k\le n} g_k\ten x_k\rrm_{R_p\wet L_1(\Mz_m)}  \\
 &\pll \sim_{c(p,\la)}
 \inf_{x_k=c_{kj}+d_{kj}} \noo \bigg(\summ_{k\le n,|j|\le c(p)\frac{\log n}{\log \la} } \la^{\frac{j}{p'}} c_{jk}^*c_{jk}\bigg)^{\frac12}\rrm+
  \noo \bigg(\summ_{k\le n,|j|\le c(p)\frac{\log n}{\log \la} }  \la^{-\frac{j}{p}}
  d_{jk}d_{jk}^*\bigg)^{\frac12}\rrm\pl
  \end{align*} \end{samepage}
\noindent holds for all $m\in \nz$  and  sequences $x_k\in
L_1(\Mz_m)$.
\end{lemma}

\begin{proof} Note that the lower estimate
 \begin{align*}
 & \inf_{x_k=c_{kj}+d_{kj}}
 \noo \bigg(\summ_{k\le n,|j|\le c(p)\frac{\log n}{\log \la}
 } \la^{\frac{j}{p'}} c_{jk}^*c_{jk}\bigg)^{\frac12}\rrm+
  \noo \bigg(\summ_{k\le n, |j|\le c(p)\frac{\log n}{\log \la}}  \la^{-\frac{j}{p}}
  d_{jk}d_{jk}^*\bigg)^{\frac12}\rrm \\
  &\pl \le  \inf_{x_k=c_{kj}+d_{kj}}
 \noo \bigg(\summ_{k\le n j\in \zz, } \la^{\frac{j}{p}} c_{jk}^*c_{jk}\bigg)^{\frac12}\rrm+
  \noo \bigg(\summ_{k\le n,j\in \zz}  \la^{-\frac{j}{p'}}
  d_{jk}d_{jk}^*\bigg)^{\frac12}\rrm \pl
  \end{align*}
is obvious. Now, we assume that the right hand side is $\le 1$. In
particular, we may find $c_k=c_{k,0}$ and $d_k=d_{k,0}$ such that
 \[ \noo (\summ_{k=1}^n  c_k^*c_k)^{\frac12}\rrm+
 \noo (\summ_{k=1}^n  d_kd_k^*)^{\frac12}\rrm \le 1 \pl .\]
This implies
 \[ \max\{\| (\summ_k x_kx_k^*)^{\frac12} \|,\|(\summ_k x_k^*x_k)^{\frac12} \| \}
 \kl 1+\sqrt{n} \pl .\]
Then, we see that for the smallest integer $j_0\ge c(p,\la)\log n$
we have
 \begin{align*}
  \noo (\summ_{k\le n, j\ge j_0}  \la^{-\frac{j}{p}}
  x_kx_k^*)^{\frac12}\rrm &\le (1+\sqrt{n})
  (\summ_{j\ge j_0} \la^{-\frac{j}{p}})^{\frac12}
  \kl
   (1-\la^{-\frac{1}{p}})^{-1}
  \la^{-\frac{j_0}{2p}} 2\sqrt{n}\pl.
  \end{align*}
Thus for this part we need $j_0\ge \frac{p}{\log \la}\log n$ in
order to obtain a constant independent of $n$. If we require $c(p)
\lel 2\max\{p,p'\}$, i.e. $c(p,\la)=c(p)/\log \la$,  we can
estimates both tails by a constant $C(p,\la)$. \qd

\begin{cor}\label{eps} Let $1<p<\infty$, $\eps>0$ and  $\nen$. Then $[R_n,C_n]_{\frac1p}$ completely embeds into $S_1^{n^{\eps n}}$.
The constant depends only on $p$ and $\eps$.
\end{cor}

\begin{proof} We use
the index set $I=\{1,...,n\}\times \{j\in \zz: |j|\le c(p) (\log
\la)^{-1}\log n\}$ which has cardinality $m\le 2c(p)\frac{n\log
n}{\log \la}$. We may chose $\la$ large enough such that $2^m\le
n^{\eps n}$. Let us note that for a finite sequence $(\mu_{kj})$
the $L_1$ space $L_1(\N(\mu))$ and $S_1^{2^m}$ are canonically
isomorphic by sending $D_\mu\in L_1(\N(\mu))$ to
 \[ \hat{D}_{\mu} \lel \ten_{j,k} \kla \begin{array}{cc} 1-\mu_j& 0 \\0&
 \mu_j \end{array}\mer  \in S_1^{2^m} \pl \]
and extending this map to a $\N(\mu)=\Mz_{2^m}$--bimodule map.
Following Remark \ref{concr} we shall define
 \[ \mu_{kj}\lel \frac{\la^{-j/p}}{\la^{j/p'}+\la^{-j/p}} \lel
 (1+\la^{j})^{-1}\pl .\]
We recall $\nu_{kj}+\la_{kj}=\la^{j/p'}+\la^{-j/p}$ and
\[ \hat{D}_{\mu}^{1/2}a_{kj}\hat{D}_{\mu}^{1/2} \lel
(1-\mu_{kj})^{\frac12}\mu_{kj}^{-\frac12}a_{kj}\hat{D}_{\mu} \lel
 \la^{\frac j2} a_{kj}\hat{D}_{\mu} \pl . \] Thus the embedding
$v^{op}$ from Remark \ref{concr} yields
 \[ v^{op}_{\la}(g_k) \lel \summ_{|j|\le c(p) (\log \la)^{-1}\log n}
 \sqrt{\la^{j/p'}+\la^{-j/p}} \pl  \la^{\frac{j}{2}}
  a_{kj}\hat{D}_{\mu} \in S_1^{2^m}  \pl .\]
The cb-norm of  this  isomorphism $v^{op}_{\la}$ and its inverse
depend only on $\la$ and $p$.\qd

At the end of this section, we  will describe a completely
isomorphic embedding of $OH$ in $B(H)$ in terms of the Brown
algebra $U_m^{nc}$ introduced in \cite{Bro}, see also
\cite{McCl}. The algebra $U_m^{nc}$ is the universal algebra
spanned by elements $(u_{ij})_{i,j=1}^m$ such that every unitary
$U=(U_{ij})\in \Ma_m(B(H))$ induces a representation
$\pi_U:U_{m}^{nc}\to B(H)$ satisfying
 \[ \pi_U(u_{ij}) \lel U_{ij} \pl .\]
Equivalently $U_{m}^{nc}$ is the universal algebra spanned by
contractions $u_{ij}$ satisfying
 \[ \summ_{k=1}^m u_{ki}^*u_{kj}\lel \delta_{ij} \lel
  \summ_{k=1}^m u_{ik}u_{jk}^*
   \pl .\]
Let us state some elementary facts. Every complete contraction
$T:S_1^m\to B(H)$ admits a complete positive extension
$\hat{T}:U_{m}^{nc}\to B(H)$. Indeed, consider then  matrix
$v=[T(e_{ij})]\in \Mz_m(B(H))$ which is   a contraction. Then, we
may define the unitary
 \[ U\lel \kla \begin{array}{cc} v &
 \sqrt{1-vv^*}\\-\sqrt{1-v^*v} & v\end{array} \mer \in
 \Mz_2(\Mz_{m}(B(H))) \pl .\]
Let $f_{ij}$ be the matrix units in $\Mz_2$. Then $\hat{T}(a)\lel
f_{11}\pi_U(a)f_{11}$ defines  the completely positive extension
of $T$. Let us now show that the map $\iota(e_{ij})=u_{ij}$ is a
complete isometry.  Let $a_{ij}\in \Mz_k$ be matrices. Then
 \begin{align}
  &\|\summ_{ij} a_{ij}\ten e_{ij}\|_{\Mz_k\ten_{\min}S_1^m}
  \lel  \sup_{\|T\|_{cb}\le 1} \|\summ a_{ij}\ten T(e_{ij})\|\lel
    \sup_{\|U\|\le 1} \|\summ a_{ij}\ten f_{11}U_{ij}f_{11}\| \nonumber \\
  &\le
  \sup_{\|U\|\le 1} \|\summ a_{ij}\ten U_{ij}\| \lel   \|\summ_{ij} a_{ij}\ten u_{ij}
  \|_{\Mz_k\ten_{\min}U_{m}^{nc}} \kl   \sup_{\|T\|_{cb}\le 1} \|\summ a_{ij}\ten T(e_{ij})\|
  \label{brown}
  \pl .
 \end{align}
The last line follows from the fact that a unitary $U=[U_{ij}]$
yields a complete contraction $T(e_{ij})=U_{ij}$. Thus we have
equality and $\iota$ is completely isometric.

\begin{cor} Let $1< p<\infty$ and $\eps>0$. Let $m\ge n^{\eps n}$.
Then $[C_n,R_n]_{\frac1p}$ embeds completely isomorphically {\rm
(}uniformly in $n${\rm )} into $U_{m}^{nc}$ using a linear
combination of the generators.
\end{cor}

\begin{proof} It is sufficient to note that the map $u:S_1^m\to
U_{m}^{nc}$ given by $u(e_{ij}) \lel u_{ij}$ is a complete
isometry. Then Corollary \ref{eps} implies the assertion.\qd

\begin{rem} {\rm Alternatively, we may use the completely isometric embedding
$u:S_1^m\to M_m\dot{\ast}M_m$ in the full free product of matrix
algebras (amalgamated over $1$) given by
 \[ u(e_{ij}) \lel e_{1j}\ast e_{i1} \pl .\]
We refer to \cite{Har} (see also \cite{Po}) for this complete
isometry.}
\end{rem}

We will now describe an infinite dimensional analogue of the Brown
algebra. We fix a sequence $(\mu_k)$ and the state
$\phi_{\mu}=\ten_{k\in \nz} \phi_{\mu_k}$ on $\ten_{k\in \nz}
\Mz_2$. Then we have conditional expectations
 \[ E_n:\Mz_{2^{\infty}}\to \Mz_{2^n}\pl ,\pl E_n(x_1\ten \cdots
 ) \lel x_1\ten \cdots \ten x_n \prodd_{k>n}\phi_{\mu_k}(x_k) \pl
 .\]
We use the notation $\E_n\lel E_n\ten id$.  A sequence $y=(y_n)$
of operators in $\Mz_{2^{\infty}}\ten B(\ell_2)$ is called
\emph{adapted} if $y_n\in \Mz_{2^n}\ten B(\ell_2)$. An adapted
sequence $(y_n)$ is a \emph{martingale} if in addition
$\E_n(y_m)=y_n$ holds for all $n\le m$. We define the set
 \begin{align*}
  U(\mu) &=  \{ (u_n) \subset B(\ell_2) \pl: \pl (u_n) \mbox{ a martingale of
  contractions} \pl ,\pl \\
  &\quad \quad \quad  \quad \quad \quad \quad \quad \quad \pl   \forall_{n}\pl:  w-\lim_{m} \E_n(u_m^*u_m)=1=w-\lim_{m} \E_n(u_m^*u_m)\} \pl.
  \end{align*}
We observe that $U(\mu)$ is in one-to-one correspondence with the
set of unitaries in $\N(\mu)\bar{\ten}B(\ell_2)$. Indeed, if $u\in
\N(\mu)\bar{\ten} B(\ell_2)$ is a unitary, then $u_n\lel \E_n(u)$
satisfies the condition listed above. Conversely, let $(u_n)$ be a
martingale of contractions. Then we may define $u$ as a weak$^*$
limit of the $u_n$'s. Since the $E_m$'s have finite rank, we find
$u_m=\E_m(u)$. Being a martingale, we obtain from this strong and
strong$^*$ convergence of $u_m$ to $u$. Hence the conditions above
imply $u^*u=1=uu^*$. We may now define the Brown algebra as the
universal $^*$-algebra generated by coefficients
$(u_{i_1,...,i_n;j_1,...,j_n})$. To be more precise, let
$S=\bigcup_n 2^{\{1,..,n\}}\times 2^{\{1,..,n\}}$ and $q$ a
noncommutative polynomial in $(g_s)_{s\in S}$ variables.  As usual
we define
\[ \noo q\rrm_{U_{\mu}^{nc}} \lel \sup_{(u_n)\in U(\mu)}
\noo
 q(\{u_n(\vec{i}_n,\vec{j}_n): (\vec{i}_n,\vec{j}_n)\in S\})
 \rrm_{B(\ell_2)} \pl .\]
Here $\vec{i}_n=(i_1,...,i_n)$ stand for an $n$-tuple and
$u_n(\vec{i}_n,\vec{j}_n)$ are the matrix coefficients of
$u_n=\sum_{\vec{i},\vec{j}}e_{\vec{i},\vec{j}}\ten
u_n(\vec{i}_n,\vec{j}_n)$ of $u_n$.
\begin{lemma}\label{einfach} \begin{enumerate}
\item[i)] $U_{\mu}^{nc}$ is a direct limit of  subalgebras $A_n$
generated by the $g_{i_1,....,i_n,j_1,...,j_n}$'s. \item[ii)] The
map $\iota:S_1^{2^n}\to A_n$ given by
 \[
 \iota(e_{i_1,....,i_n;j_1,....,j_n}) \lel g_{i_1,....,i_n,j_1,...,j_n}\]
is a complete isometry.
 \item[iii)] There is a natural inclusion map $i_{n,n+1}:A_n\to
 A_{n+1}$ such that
  \[ i_{n,n+1}(g_{i_1,....,i_n;j_1,...,j_n})\lel (1-\mu_{n+1})
 g_{i_1,....,i_n,0;j_1,....,j_n,0}+ \mu_{n+1}
 g_{i_1,....,i_n,1;j_1,....,j_n,1}  \pl .\]
\end{enumerate}
\end{lemma}

\begin{proof} The last assertion iii) follows
immediately from the martingale property of the $(u_n)$'s. Thus by
definition of the $A_n$'s we deduce $\iota_{n,n+1}(A_n)\subset
A_{n+1}$ completely isometrically. Then assertion i) follows
immediately from the definition of $U_{\mu}^{nc}$ because
polynomials with finitely many entries are dense. For the proof of
ii) we consider a martingale sequence $(u_k)$ of contractions. In
particular, the element $u_n= \sum_{\vec{i},\vec{j}}
e_{\vec{i},\vec{j}} \ten
 u_n(\vec{i},\vec{j})$ is  a contraction. This allows us to
apply \eqref{brown}. Thus $\iota_n:S_1^{2^n}\to A_n$ is a complete
contraction. For the converse, we consider $a_{\vec{i},\vec{j}}\in
\Mz_{2^n}$ and a unitary $U\in \Mz_{2^n}(B(H))$.  We define the
martingale $u_k=\E_k(U)$ for $k\le n$ and $U_k=u$ for $k>n$. By
definition of $A_n$ we have a homomorphism $\pi_{(u_k)}:A_n\to
B(H)$ given by
$\pi_{(u_k)}(g_{i_1,...,i_k,j_1,...,j_k})=u_k(i_1,...,i_k,j_1,...,j_k)$.
Then we have
 \begin{align*}
 &\|\summ_{\vec{i}_n,\vec{j}_n} a_{\vec{i},\vec{j}}\ten
 U_{\vec{i},\vec{j}}
\|_{\Mz_m(B(H))} \! \! \lel\!\!
  \|(id \ten \pi_{(u_k)})
 (\summ_{\vec{i}_n,\vec{j}_n} a_{\vec{i}_n,\vec{j}_n}\ten
 g_{\vec{i}_n,\vec{j}_n})\|\!\! \kl\!\!
  \|\summ_{\vec{i}_n,\vec{j}_n} a_{\vec{i}_n,\vec{j}_n}\ten
 g_{\vec{i}_n,\vec{j}_n}\|_{\Mz_m\ten_{\min}A_n}\! .
 \end{align*}
This shows that $\iota_n^{-1}$ is also a complete contraction.\qd

\begin{lemma}\label{pred}  $L_1(\N(\mu))$ embeds completely isometrically in
$U_{\mu}^{nc}$.
\end{lemma}

\begin{proof} $L_1(\N(\mu))$ is the direct limit of $S_1^{2^n}$
with respect to the inclusion  map $\iota_{n,n+1}:S_1^{2^n}\to
S_1^{2^{n+1}}$ is given by
 \[ \iota_{n,n+1}(e_{i_1,...,i_n;j_1,...,j_n}) \lel (1-\mu_{n+1})
 e_{i_1,...,i_n,0;j_1,...,j_n,0}+\mu_{n+1}e_{i_1,...,i_n,1;j_1,...,j_n,1}
 \pl .\]
According to Lemma \ref{einfach}, we deduce that
 \[ v(e_{i_1,...,i_n;j_1,...,j_n})\lel g_{i_1,...,i_n;j_1,...,j_n}
 \]
extends to a complete isometry of $v:L_1(\N(\mu))\to
B_{\mu}^{nc}$. \qd

\begin{cor}  Let $(\mu_{kj})_{k\in \nz,j\in \zz}$ be defined as $\mu_{kj}=\frac{(1+|j|)^{-2}}{2}$ for $j\ge 0$ and
$\mu_{kj}=1-\frac{(1+|j|)^{-2}}{2}$ for $j<0$. Then $OH$ embeds
into $U_{\mu}^{nc}$ using a  linear combination of the generators.
\end{cor}

\begin{proof} This follows immediately from Corollary
\ref{concr5}, more precisely Remark \ref{concr} and Lemma
\ref{pred}.  \qd

{AMS Subject classification 2000:} 46L53, 47L25.

\begin{thebibliography}{GvW78}

\bibitem[AB]{AB}
Luigi Accardi and Marek Bo{\.z}ejko.
\newblock Interacting {F}ock spaces and {G}aussianization of probability
  measures.
\newblock {\em Infin. Dimens. Anal. Quantum Probab. Relat. Top.},
  1(4):663--670, 1998.

\bibitem[Ara]{Ar}
Huzihiro Araki.
\newblock A lattice of von {N}eumann algebras associated with the quantum
  theory of a free {B}ose field.
\newblock {\em J. Mathematical Phys.}, 4:1343--1362, 1963.

\bibitem[AW63]{AW1}
H.~Araki and E.~J. Woods.
\newblock Representations of the canonical commutation relations describing a
  nonrelativistic infinite free {B}ose gas.
\newblock {\em J. Mathematical Phys.}, 4:637--662, 1963.

\bibitem[AW69]{AW}
Huzihiro Araki and E.~J. Woods.
\newblock A classification of factors.
\newblock {\em Publ. Res. Inst. Math. Sci. Ser. A}, 4:51--130, 1968/1969.

\bibitem[Bla]{Bla}
Etienne Blanchard.
\newblock A few remarks on exact {$C(X)$}-algebras.
\newblock {\em Rev. Roumaine Math. Pures Appl.}, 45(4):565--576 (2001), 2000.


\bibitem[BS]{BS}
Philippe Biane and Roland Speicher.
\newblock Stochastic calculus with respect to free {B}rownian motion and
  analysis on {W}igner space.
\newblock {\em Probab. Theory Related Fields}, 112(3):373--409, 1998.

\bibitem[Bro]{Bro}
Lawrence~G. Brown.
\newblock Ext of certain free product {$C\sp{\ast} $}-algebras.
\newblock {\em J. Operator Theory}, 6(1):135--141, 1981.

\bibitem[BGS]{Sch}
Anis Ben~Ghorbal and Michael Sch{\"u}rmann.
\newblock Non-commutative notions of stochastic independence.
\newblock {\em Math. Proc. Cambridge Philos. Soc.}, 133(3):531--561, 2002.

\bibitem[CGH]{cock}
A.~M. Cockroft, S.~P. Gudder, and R.~L. Hudson.
\newblock A quantum-mechanical functional central limit theorem.
\newblock {\em J. Multivariate Anal.}, 7(1):125--148, 1977.

\bibitem[CH]{Cu-Hud}
C.~D. Cushen and R.~L. Hudson.
\newblock A quantum-mechanical central limit theorem.
\newblock {\em J. Appl. Probability}, 8:454--469, 1971.

\bibitem[Con]{Co}
A.~Connes.
\newblock Classification of injective factors. {C}ases {$II\sb{1},$}
  {$II\sb{\infty },$} {$III\sb{\lambda },$} {$\lambda \not=1$}.
\newblock {\em Ann. of Math. (2)}, 104(1):73--115, 1976.

\bibitem[Dix]{Dix}
Jacques Dixmier.
\newblock {\em von {N}eumann algebras}, volume~27 of {\em North-Holland
  Mathematical Library}.
\newblock North-Holland Publishing Co., Amsterdam, 1981.
\newblock With a preface by E. C. Lance, Translated from the second French
  edition by F. Jellett.

\bibitem[ER]{ER}
Ed. Effros and Z-J. Ruan.
\newblock {\em Operator spaces}, volume~23 of {\em London Mathematical Society
  Monographs. New Series}.
\newblock The Clarendon Press Oxford University Press, New York, 2000.

\bibitem[Gro]{Groh}
Ulrich Groh.
\newblock Uniform ergodic theorems for identity preserving {S}chwarz maps on
  {$W\sp{\ast} $}-algebras.
\newblock {\em J. Operator Theory}, 11(2):395--404, 1984.

\bibitem[GvW]{GiW}
N.~Giri and W.~von Waldenfels.
\newblock An algebraic version of the central limit theorem.
\newblock {\em Z. Wahrscheinlichkeitstheorie und Verw. Gebiete},
  42(2):129--134, 1978.

\bibitem[Har]{Har}
Asma Harcharras.
\newblock On some ``stability'' properties of the full {$C\sp *$}-algebra
  associated to the free group {$F\sb \infty$}.
\newblock {\em Proc. Edinburgh Math. Soc. (2)}, 41(1):93--116, 1998.

\bibitem[Heg]{Heg}
Gerhard~C. Hegerfeldt.
\newblock Noncommutative version of the central limit theorem and of
  {C}ram\'er's theorem.
\newblock In {\em Stochastic processes and their applications in mathematics
  and physics (Bielefeld, 1985)}, volume~61 of {\em Math. Appl.}, pages
  187--202. Kluwer Acad. Publ., Dordrecht, 1990.

\bibitem[Hia]{Hiai}
F.~Hiai.
\newblock {$q$}-deformed {A}raki-{W}oods algebras.
\newblock In {\em Operator algebras and mathematical physics (Constan\c ta,
  2001)}, pages 169--202. Theta, Bucharest, 2003.

\bibitem[HP]{HP}
Uffe Haagerup and Gilles Pisier.
\newblock Bounded linear operators between {$C\sp *$}-algebras.
\newblock {\em Duke Math. J.}, 71(3):889--925, 1993.

\bibitem[HRS]{HRS}
U.~Haagerup, H.~P. Rosenthal, and F.~A. Sukochev.
\newblock Banach embedding properties of non-commutative {$L\sp p$}-spaces.
\newblock {\em Mem. Amer. Math. Soc.}, 163(776):vi+68, 2003.

\bibitem[Hud]{Hud0}
R.~L. Hudson.
\newblock A quantum-mechanical central limit theorem for anti-commuting
  observables.
\newblock {\em J. Appl. Probability}, 10:502--509, 1973.

\bibitem[J1]{JD}
M.~Junge.
\newblock Doob's inequality for non-commutative martingales.
\newblock {\em J. Reine Angew. Math.}, 549:149--190, 2002.

\bibitem[J2]{JF}
M.~Junge.
\newblock Fubini's theorem for ultraproducts of noncommutative {$L_p$} spaces.
\newblock {\em Canad. J. Math.}, 56 (5), 2004, p 983-1021.

\bibitem[J3]{Joh}
M.~Junge.
\newblock Embedding of the operators space {$OH$} and the logarithmic `little
  {G}rothendieck' inequality.
\newblock {\em Inventiones Mathematica\!e} {\bf  161} (2005),
p 225-286.


\bibitem[JNRX]{4-author}
M.~Junge, N.~Nielsen, Z.-J. Ruan, and Q.~Xu.
\newblock {${\mathcal C}{\mathcal O}{\mathcal L}_p$} spaces - the local
  structure of non-commutative {$L_p$} spaces.
\newblock {\em Advances in Mathematics}, {\bf
117} (2003), no.\!\! 2, p 313-341.


\bibitem[JR]{JR-ap}
M.~Junge and Z-J. Ruan.
\newblock {A}pproximation properties for noncommutative {$L\sb p$}-spaces
  associated with discrete groups.
\newblock {\em Duke Math. J.}, 117(2):313--341, 2003.

\bibitem[JSh]{JS}
M.~Junge and D.~Sherman.
\newblock Noncommutative {$L_p$} modules.
\newblock {\em Journal of Operator Theory}, 53 (2005), 3-34.

\bibitem[JX1]{JX}
M.~Junge and Q.~Xu.
\newblock {N}oncommutative {B}urkholder/{R}osenthal inequalities.
\newblock {\em Ann. Probab.}, 31(2):948--995, 2003.


\bibitem[JX2]{JX-ros}
M.~Junge and Q.~Xu.
\newblock {N}oncommutative {B}urkholder/{R}osenthal inequalities: Applications.
\newblock In preparation.


\bibitem[JX3]{JX3}
M.~Junge and Q.~Xu.
\newblock {Q}untum probablistic methods in operator spaces theory:
\newblock In preparation.

\bibitem[JR]{JRo}
M.~Junge and H.~Rosenthal.
\newblock Noncommutative {$L_p$} spaces are asymptotially stable.
\newblock In preparation.

\bibitem[Kad]{Kad}
Richard~V. Kadison.
\newblock A generalized {S}chwarz inequality and algebraic invariants for
  operator algebras.
\newblock {\em Ann. of Math. (2)}, 56:494--503, 1952.

\bibitem[K{\"o}s]{ko2}
Claus K{\"o}stler.
\newblock {P}rivate communication.

\bibitem[KS]{KS}
C.~K{\"o}stler and R.~Speicher.
\newblock {O}n the structure of non-commutative white noises.
\newblock To appear in {\it Transactions of the AMS}.

\bibitem[Kos]{Kos}
H.~Kosaki.
\newblock Applications of the complex interpolation method to a von {N}eumann
  algebra: noncommutative {$L\sp{p}$}-spaces.
\newblock {\em J. Funct. Anal.}, 56(1):29--78, 1984.

\bibitem[KR1]{kar-I}
R.~V. Kadison and J.~R. Ringrose.
\newblock {\em Fundamentals of the theory of operator algebras. {V}ol. {I}},
  volume~15 of {\em Graduate Studies in Mathematics}.
\newblock American Mathematical Society, Providence, RI, 1997.
\newblock Elementary theory, Reprint of the 1983 original.

\bibitem[KR2]{kar-II}
R.~V. Kadison and J.~R. Ringrose.
\newblock {\em Fundamentals of the theory of operator algebras. {V}ol. {II}},
  volume~16 of {\em Graduate Studies in Mathematics}.
\newblock American Mathematical Society, Providence, RI, 1997.
\newblock Advanced theory, Corrected reprint of the 1986 original.


\bibitem[K{\"u}m]{Ku}
Burkhard K{\"u}mmerer.
\newblock Survey on a theory of noncommutative stationary {M}arkov processes.
\newblock In {\em Quantum probability and applications, III (Oberwolfach,
  1987)}, volume 1303 of {\em Lecture Notes in Math.}, pages 154--182.
  Springer, Berlin, 1988.

\bibitem[LP]{LP}
F.~Lust-Piquard.
\newblock In\'egalit\'es de {K}hintchine dans {$C\sb p\;(1<p<\infty)$}.
\newblock {\em C. R. Acad. Sci. Paris S\'er. I Math.}, 303(7):289--292, 1986.

\bibitem[LPP]{LPP}
F.~Lust-Piquard and G.~Pisier.
\newblock Noncommutative {K}hintchine and {P}aley inequalities.
\newblock {\em Ark. Mat.}, 29(2):241--260, 1991.

\bibitem[LTI]{LT-I}
Joram Lindenstrauss and Lior Tzafriri.
\newblock {\em Classical {B}anach spaces. {I}}.
\newblock Springer-Verlag, Berlin, 1977.
\newblock Sequence spaces, {\em  Ergebnisse der Mathematik und ihrer Grenzgebiete},
  Vol. 92.

\bibitem[LTII]{LT-II}
J.~Lindenstrauss and L.~Tzafriri.
\newblock {\em Classical {B}anach spaces. {II}}, {\em Ergebnisse
  der Mathematik und ihrer Grenzgebiete}, Vol 97.
\newblock Springer-Verlag, Berlin, 1979.
\newblock Function spaces.

\bibitem[McC]{McCl}
Kevin McClanahan.
\newblock {$C\sp *$}-algebras generated by elements of a unitary matrix.
\newblock {\em J. Funct. Anal.}, 107(2):439--457, 1992.

\bibitem[Nou]{Nou}
Alexandre Nou.
\newblock Asymptotic matricial models and {QWEP} property for
  {$q$}-{A}raki--{W}oods algebras.
\newblock {\em J. Funct. Anal.}, 232(2):295--327, 2006.

\bibitem[OL]{OL}
{\em Lie groups and {L}ie algebras. {I}}, volume~20 of {\em
Encyclopaedia of
  Mathematical Sciences}.
\newblock Springer-Verlag, Berlin, 1993.
\newblock Foundations of Lie theory. Lie transformation groups, A translation
  of {\it Current problems in mathematics. Fundamental directions. Vol. 20}
  (Russian), Akad.\ Nauk SSSR, Vsesoyuz.\ Inst.\ Nauchn.\ i Tekhn.\ Inform.,
  Moscow, 1988 , Translation by A. Kozlowski, Translation edited by A. L.
  Onishchik.

\bibitem[Pau]{Pau}
V.~Paulsen.
\newblock {\em Completely bounded maps and operator algebras}, volume~78 of
  {\em Cambridge Studies in Advanced Mathematics}.
\newblock Cambridge University Press, Cambridge, 2002.


\bibitem[Ps1]{Po}
G.~Pisier.
\newblock {\em Introduction to operator space theory}, volume 294 of {\em
  London Mathematical Society Lecture Note Series}.
\newblock Cambridge University Press, Cambridge, 2003.

\bibitem[Ps2]{Pl}
G.~Pisier.
\newblock Completely bounded maps into certain hilbertian operator spaces.
\newblock {\em Internat. Math. Res. Notices}, 74  3983-4018, 2004.

\bibitem[PS]{PS}
G.~Pisier and D.~Shlyakhtenko.
\newblock Grothendieck's theorem for operator spaces.
\newblock {\em Invent. Math.}, 150(1):185--217, 2002.

\bibitem[Qua]{Quaeg}
J.~Quaegebeur.
\newblock A noncommutative central limit theorem for {CCR}-algebras.
\newblock {\em J. Funct. Anal.}, 57(1):1--20, 1984.

\bibitem[Ray]{Ra}
Y.~Raynaud.
\newblock On ultrapowers of non commutative {$L\sb p$} spaces.
\newblock {\em J. Operator Theory}, 48(1):41--68, 2002.

\bibitem[Rid]{Rid}
G.~Rideau.
\newblock On some representations of the anticommutations relations.
\newblock {\em Comm. Math. Phys.}, 9:229--241, 1968.

\bibitem[Shl]{S}
D.~Shlyakhtenko.
\newblock Free quasi-free states.
\newblock {\em Pacific J. Math.}, 177(2):329--368, 1997.

\bibitem[Spe]{Sp}
Roland Speicher.
\newblock A noncommutative central limit theorem.
\newblock {\em Math. Z.}, 209(1):55--66, 1992.

\bibitem[Str]{Strat}
{\c{S}}erban Str{\u{a}}til{\u{a}}.
\newblock {\em Modular theory in operator algebras}.
\newblock Editura Academiei Republicii Socialiste Rom\^ania, Bucharest, 1981.
\newblock Translated from the Romanian by the author.

\bibitem[SZ]{Strat-Z}
{\c{S}}erban Str{\u{a}}til{\u{a}} and L{\'a}szl{\'o} Zsid{\'o}.
\newblock {\em Lectures on von {N}eumann algebras}.
\newblock Editura Academiei, Bucharest, 1979.
\newblock Revision of the 1975 original, Translated from the Romanian by Silviu
  Teleman.

\bibitem[Tak1]{Tak}
M.~Takesaki.
\newblock {\em Theory of operator algebras. {I}}.
\newblock Springer-Verlag, New York, 1979.

\bibitem[Tak2]{TK-II}
M.~Takesaki.
\newblock {\em Theory of operator algebras. {II}}, volume 125 of {\em
  Encyclopaedia of Mathematical Sciences}.
\newblock Springer-Verlag, Berlin, 2003.
\newblock Operator Algebras and Non-commutative Geometry 6.

\bibitem[Tak3]{TK-III}
M.~Takesaki.
\newblock {\em Theory of operator algebras. {III}}, volume 127 of {\em
  Encyclopaedia of Mathematical Sciences}.
\newblock Springer-Verlag, Berlin, 2003.
\newblock Operator Algebras and Non-commutative Geometry 8.

\bibitem[Ter]{Te}
M.~Terp.
\newblock {$L_p$} spaces associated with von {N}eumann algebras.
\newblock Notes~~, Math. Institute, Copenhagen Univ. ~ 1981.
\newblock

\bibitem[Voi]{Voi2}
Dan Voiculescu.
\newblock A strengthened asymptotic freeness result for random matrices with
  applications to free entropy.
\newblock {\em Internat. Math. Res. Notices}, (1):41--63, 1998.

\bibitem[vW]{vWa}
W.~von Waldenfels.
\newblock An algebraic central limit theorem in the anti-commuting case.
\newblock {\em Z. Wahrscheinlichkeitstheorie und Verw. Gebiete},
  42(2):135--140, 1978.

\bibitem[Wid]{Wid}
David~Vernon Widder.
\newblock {\em The {L}aplace {T}ransform}.
\newblock Princeton Mathematical Series, v. 6. Princeton University Press,
  Princeton, N. J., 1941.

\end{thebibliography}
\end{document}